\documentclass[11pt,a4paper]{article}
\usepackage[T1]{fontenc}
\usepackage{lmodern}
\usepackage[margin=26mm]{geometry}
\usepackage{amsmath,amssymb,amsthm,mathtools,microtype,cite}
\usepackage[hidelinks]{hyperref}
\mathtoolsset{showonlyrefs}
\numberwithin{equation}{section}
\newtheorem{theorem}{Theorem}[section]
\newtheorem{proposition}[theorem]{Proposition}
\newtheorem{lemma}[theorem]{Lemma}
\newtheorem{corollary}[theorem]{Corollary}
\theoremstyle{remark}

\newcommand{\E}{\mathbb E}
\newcommand{\Prob}{\mathbb P}
\newcommand{\Tr}{\operatorname{Tr}}
\newcommand{\Cov}{\operatorname{Cov}}
\newcommand{\Var}{\operatorname{Var}}
\newcommand{\Law}{\mathcal L}
\newcommand{\Poi}{\operatorname{Poisson}}
\newcommand{\Wtwo}{\mathsf W_2}
\newcommand{\TV}{d_{\mathrm{TV}}}
\newcommand{\BL}{d_{\mathrm{BL}}}
\newcommand{\ind}{\mathbf 1}
\newcommand{\cA}{\mathcal A}
\newcommand{\cR}{\mathcal R}

\newcommand{\cF}{\widehat{\mathcal F}}
\newcommand{\cZ}{\mathcal Z}
\newcommand{\op}{\mathrm{op}}
\newcommand{\dd}{\,\mathrm d}

\newcommand{\Ga}{\mathcal C}
\newcommand{\new}{\mathrm{new}}

\newcommand{\perm}{\mathfrak S}

\title{Growing-time trace fluctuations and self-averaging form factors\\of Grover walks on random lifts}
\author{
\normalsize \href{https://orcid.org/0009-0009-5064-0499}{Nino~Nagata}$^{1}$,\quad \href{https://orcid.org/0000-0001-9979-1882}{Ryoichi~Suzuki}$^{2,*}$\\[10pt]
\parbox{0.96\textwidth}{\small\centering
$^1$ Department of Mathematical Sciences, College of Science and Engineering,\\
Ritsumeikan University, Kusatsu, Shiga 525-8577, Japan\\[3pt]
$^2$ School of Management, Tokyo University of Science,\\
Tokyo 102-0071, Japan\\[5pt]
$^*$ Corresponding author: Ryoichi Suzuki;\\
\href{mailto:rsuzukimath@gmail.com}{\texttt{rsuzukimath@gmail.com}}
}}
\date{}
\hypersetup{pdftitle={Growing-time trace fluctuations and self-averaging form factors of Grover walks on random lifts},pdfauthor={Nino Nagata; Ryoichi Suzuki}}

\begin{document}
\maketitle
\begin{abstract}
We prove a functional central limit theorem for trace statistics of Grover
walks on random $N$-lifts of a fixed finite connected simple $d$-regular
graph, $d\geq3$. All observations are computed from the same random lift.
The squared centred traces, divided by their observation times, have
averages converging in probability to a degree-dependent constant. Their
partial sums, centred at this constant and divided by the square root
of the window length, converge to Brownian motion. For a window of
length $M$ beginning at $t\geq1$, it suffices that $M\to\infty$ and
$R^4(d-1)^R/N\to0$, where $R=t+M$. The initial time may remain fixed,
with no further restriction on $M/t$. Convergence holds jointly with
fixed low-degree polynomial trace statistics; their limiting vector is
independent of the Brownian motion. For non-bipartite bases, we also
identify Gaussian trace limits at fixed time offsets. Under common
initial-time normalisation, subtracting a deterministic growth profile
gives an explicit Brownian time change when $M/t$ has a finite limit.
Bipartite bases retain the same averaged constant and have twice the
Brownian variance rate. The proof combines a growing-dimensional Poisson
approximation, a uniform bilinear Chebyshev-energy estimate, a
time-decaying Poisson--Gaussian coupling, and a Gaussian quadratic-array
theorem. Separate moment estimates give expectation-centred versions
under an additional sufficient condition.
\end{abstract}


\section{Introduction and main results}\label{sec:main}
Can averaging a dependent sequence of trace statistics on one random
graph produce a law of large numbers and a functional central limit
theorem? We study this question for Grover walks on random lifts of a
fixed regular base. An $N$-lift replaces each base vertex by a fibre of
$N$ copies and joins adjacent fibres by independent uniformly random
perfect matchings. The graph is sampled once; every observation is then
computed from the same graph by varying the power of its Grover operator.
Thus increasing the observation window does not provide independent
samples of graphs.

The trace of a power sums contributions from all eigenphases. The
expectation of its squared mean-centred trace is the connected spectral
form factor. For a non-bipartite base, the squared centred trace divided
by its observation time has a nondegenerate limit at a single growing
time. Averaging the observation-time-normalised squared traces instead
recovers an explicit deterministic value. This question of temporal
self-averaging has also been studied for quantum graphs and random-matrix
ensembles \cite{Prange1997,KottosSmilansky1999,BSW2002,GnutzmannSmilansky2006,HSW,CES}.

Our main result describes the cumulative averaging error by a Brownian
limit while the number of sheets and the observation window grow
simultaneously. The initial time may remain fixed: within the stated
logarithmic cutoff range, the window need not be short relative to its
initial time. The limit is joint with fixed low-degree polynomial trace
statistics and is independent of their limiting Poisson structure.
Bipartite and non-bipartite bases have the same averaged limit, but the
Brownian variance rate doubles in the bipartite case. Dividing all
observations by a common initial time instead leaves a deterministic
growth profile and, on proportional windows, a nonconstant variance rate.

Theorem~\ref{S:thm:natural} gives the functional limit and
Theorem~\ref{S:thm:lln} the law of large numbers.
Corollaries~\ref{S:cor:macro} and~\ref{S:cor:drift} describe the common-time
normalisation. The subsequent fixed-offset theorem identifies the
fluctuations that remain at a single time. Both exact mean centring and
a centring computable from the base and the number of sheets are
available. All trace observables use the full Grover spectrum.

\paragraph{Main contributions.}
The proof requires uniform estimates because both the number of cycle
lengths and the number of observations increase. A fixed-dimensional
cycle limit or a stationary Gaussian approximation at fixed time offsets
alone would not control the cumulative error. Theorem~\ref{L:thm:poisson}
provides a total-variation error of order $R^4q^R/N$, where $q=d-1$,
for the entire cycle-count vector up to the cutoff $R$.
Theorem~\ref{S:thm:energy} evaluates the bilinear Chebyshev energy with
a bounded remainder at all positive observation times.
Proposition~\ref{F:prop:alltime} retains time decay in a simultaneous
Poisson--Gaussian coupling. Theorem~\ref{F:thm:gaussian} then gives an
invariance principle for Gaussian quadratic arrays with a product bound
on the covariance error, rather than exact stationarity. These estimates
allow the initial observations to have nonvanishing covariance errors
while their cumulative contribution becomes negligible. The exact trace
identity in Section~\ref{sec:inputs} and separate finite-lift moment
bounds connect these arguments to the stated trace statistics.

\subsection{Model and observables}

Let $B=(V,E)$ be a fixed finite connected simple $d$-regular graph, where $d\geq3$. Set
\[
 q=d-1,\qquad n_0=|V(B)|,\qquad m_0=|E(B)|,
 \qquad r_d=\frac{2\sqrt q}{d},\qquad \alpha_d=\arccos r_d.
\]
Thus $0<r_d<1$ and $0<\alpha_d<\pi/2$. Choose an orientation for each base edge and assign independent uniform permutations of $[N]=\{1,\ldots,N\}$ to the oriented edges. Traversal in the reverse direction uses the inverse permutation. For an oriented edge $e=(u,v)$ with permutation $\sigma_e$, join $(u,a)$ to $(v,\sigma_e(a))$ for every $a\in[N]$. This defines the random $N$-lift $G_N$; we do not condition on connectedness. We write $\perm_N$ for the symmetric group on $[N]$. Constants $C_B$ may change from line to line and depend only on the fixed base, unless further dependence is indicated by subscripts.

For a finite simple $d$-regular graph $G$, let $A_G$ be its adjacency matrix and $D(G)$ its directed-edge set. Write $o(e),t(e)$ for the initial and terminal vertices of a directed edge, and $\bar e$ for its reversal. We use $I_n$ for the $n\times n$ identity matrix and $\mathbf1$ for the all-ones column vector of the dimension indicated by the surrounding formula; $M^{\mathsf T}$ and $M^*$ denote transpose and conjugate transpose. The untwisted Grover matrix is
\[
 (U_G)_{e,f}=\frac{2}{d}\ind_{\{t(f)=o(e)\}}-\ind_{\{e=\bar f\}}.
\]
It is real and orthogonal: it is the product of the arc-reversal involution and the direct sum of the local reflections $(2/d)\mathbf1\mathbf1^{\!*}-I_d$, each acting on the $d$ arcs terminating at a vertex. We use the centred real trace
\begin{equation}\label{eq:At}
 \cA_N(s)=\Tr U_{G_N}^{\,s}-\E[\Tr U_{G_N}^{\,s}],\qquad s\in\mathbb Z_{\geq0},
\end{equation}
and the connected form factor with no matrix-dimension normalisation,
\begin{equation}\label{eq:formfactor}
 \mathcal F_N^{\mathrm c}(s)
   =\E[|\cA_N(s)|^2]
   =\E[|\Tr U_{G_N}^{\,s}|^2]-|\E[\Tr U_{G_N}^{\,s}]|^2.
\end{equation}
The symbols $\Prob$ and $\E$ refer to the probability law or coupling under consideration. For quantities defined from $G_N$, expectations are over the random covering permutations. The auxiliary Poisson and Gaussian laws are specified when they are introduced. Throughout, the full Grover spectrum is retained.

For any graph $G$, let $H_G$ denote its directed-edge non-backtracking matrix. For the base,
\[
 (H_B)_{e,f}=\ind_{\{t(e)=o(f),\ f\ne\bar e\}}.
\]
The trace $\Tr H_G^k$ counts rooted, oriented cyclically non-backtracking closed walks: the non-backtracking condition is imposed also on the transition from the last edge to the first. A matrix is primitive if all entries of some positive integer power are positive. Under our standing regularity assumptions, $H_B$ is primitive exactly when $B$ is non-bipartite; see Lemma~\ref{S:lem:period}. The time-averaging results below cover both bipartite and non-bipartite bases. Define
\begin{equation}\label{eq:kernel}
 \Gamma_d(h)=\frac{2}{\pi}\int_{\alpha_d}^{\pi-\alpha_d}\cos(h\theta)\dd\theta,
 \qquad h\in\mathbb Z.
\end{equation}

The variance coefficient and the period are
\[
 c_d=\Gamma_d(0)=\frac4\pi\arcsin\left(\frac{2\sqrt{d-1}}d\right),
 \qquad
 \mathfrak p_B=\begin{cases}1,&B\text{ is not bipartite},\\2,&B\text{ is bipartite}.\end{cases}
\]
The set $\mathcal I_d=[\alpha_d,\pi-\alpha_d]\cup[-\pi+\alpha_d,-\alpha_d]$ is the limiting dispersive Grover band, and $\Gamma_d$ is the Fourier transform of $\pi^{-1}\mathbf 1_{\mathcal I_d}$. Hence $c_d=|\mathcal I_d|/\pi=2\,|\mathcal I_d|/(2\pi)$: it is the orthogonal-class ramp slope $2$ multiplied by the fraction of the quasi-energy circle occupied by the band. By Lemma~\ref{S:lem:period}, $\mathfrak p_B$ is the period of $H_B$. For a cutoff $R$, write
\begin{equation}\label{eq:cutoff-errors}
 \tau_{N,R}=R^4q^R/N,\qquad \delta_{N,R}=R^5q^R/N.
\end{equation}
The quantity $\tau_{N,R}$ controls the cycle approximation in total variation; $\delta_{N,R}$ controls covariance errors. The fixed polynomial trace statistics $Y_{N,k}$ and their limits $Y_k$ are defined in Section~\ref{sec:inputs}. We use $W$ for standard Brownian motion. Probability metrics and path-space conventions are specified at the end of this section.

\subsection{Time averaging on one lift}
We distinguish two choices of centring and two choices of time normalisation. Set
\[
 \cA_N^0(s)=\cA_N(s),\qquad \cA_N^*(s)=\Tr U_{G_N}^s-m_N(s),
\]
where $m_N(s)$ is the explicit finite sum defined in Section~\ref{TA:sec:centering}. The superscripts $0$ and $*$ distinguish exact mean centring and this explicit centring. The superscripts $0$ and $*$ are labels, not powers or adjoints. For $\star\in\{0,*\}$, define
\[
 \cF_{N,t,M}^{\natural,\star}=\frac1M\sum_{h=0}^{M-1}\frac{\cA_N^\star(t+h)^2}{t+h},
 \qquad
 \cF_{N,t,M}^{\star}=\frac1{Mt}\sum_{h=0}^{M-1}\cA_N^\star(t+h)^2.
\]
Every summand is evaluated on the same random lift. The symbol $\natural$ means that each squared trace is divided by its own observation time. The cumulative error for this normalisation is
\begin{equation}\label{S:eq:naturalprocess}
 \cZ_{N,t,M}^{\natural,\star}(u)=\frac1{\sqrt M}
  \sum_{h=0}^{\lfloor Mu\rfloor-1}
   \left\{\frac{\cA_N^\star(t+h)^2}{t+h}-c_d\right\},
 \qquad 0\leq u\leq1,
\end{equation}

Empty sums are zero. Each process is right-continuous and constant between consecutive points $u=j/M$; its polygonal interpolation is linear between consecutive grid values.

\begin{theorem}[Fluctuations with observation-time normalisation]\label{S:thm:natural}
Let $B$ be any fixed finite connected simple $d$-regular graph with $d\geq3$. Let $t_N,M_N$ be positive integers such that
\begin{equation}\label{S:eq:conditions}
 t_N\geq1,\qquad M_N\longrightarrow\infty,\qquad
 R_N=t_N+M_N,\qquad \tau_{N,R_N}\longrightarrow0.
\end{equation}
Then, for either trace centring,
\[
 \cZ_{N,t_N,M_N}^{\natural,\star}
   \Longrightarrow2\sqrt{\mathfrak p_Bc_d}\,W
 \quad\text{in }D([0,1])\text{ with }J_1.
\]
The polygonal interpolations converge in $C([0,1])$. The processes obtained from the two centrings are uniformly close in probability and converge jointly to the same Brownian motion. For every fixed $L\geq3$, the convergence holds jointly with $(Y_{N,k})_{k=3}^L$, whose limit $(Y_k)_{k=3}^L$ is independent of $W$. In particular,
\[
 \sqrt{M_N}\bigl(\cF_{N,t_N,M_N}^{\natural,\star}-c_d\bigr)
       \Longrightarrow\mathcal N(0,4\mathfrak p_Bc_d).
\]
\end{theorem}

Apart from the common cutoff condition, Theorem~\ref{S:thm:natural}
imposes no restriction on the ratio $M_N/t_N$. In particular, $t_N$ may
remain fixed. A sufficient
graph-size restriction is $t_N+M_N\leq(1-\varepsilon)\log_qN$ for a fixed
$0<\varepsilon<1$. For example, $t_N=1$ and
$M_N=\lfloor(1-\varepsilon)\log_qN\rfloor-1$ satisfy this restriction
for all sufficiently large $N$. Windows comparable to, or longer
than, their initial time are also allowed. For a bipartite base, odd-time traces vanish and the leading even-time trace variance doubles; after averaging over consecutive times, the Brownian variance rate is twice the non-bipartite value.

The law of large numbers distinguishes the two normalisations. Dividing
each summand by its own observation time recovers the same constant on
all the preceding windows. With a common initial time, the limiting average
also records deterministic growth across the window.

\begin{theorem}[Self-averaging for all regular bases]\label{S:thm:lln}
Under \eqref{S:eq:conditions}, for either trace centring,
\[
 \cF_{N,t_N,M_N}^{\natural,\star}\longrightarrow c_d
 \quad\text{in probability}.
\]
If, in addition, $M_N/t_N\to\rho\in[0,\infty)$, then
\[
 \cF_{N,t_N,M_N}^{\star}\longrightarrow c_d(1+\rho/2)
 \quad\text{in probability}.
\]
The respective limits hold in $L^1$ under the additional condition
\begin{equation}\label{S:eq:averagedmoment}
 \frac{\delta_{N,R_N}}{M_N}
       =\frac{R_N^5q^{R_N}}{NM_N}\longrightarrow0.
\end{equation}
\end{theorem}

If $R_N/M_N$ is bounded, the extra moment condition follows from
$\tau_{N,R_N}\to0$. This includes windows with a fixed initial time and
windows for which $M_N/t_N$ converges to a positive constant.
The proof in Section~\ref{S:sec:llnproof} uses direct ideal second moments
and finite-lift expectation estimates, independently of the functional
limit theorem. Section~\ref{sec:timeaverages} also gives a finite-size
probability bound for the common-time statistic on shorter windows when
the base is non-bipartite.

\subsection{The drift caused by a common time normalisation}
Dividing all squared traces by the initial time leaves a deterministic growth term across the window. To compare with the preceding process, define
\begin{equation}\label{F:eq:process}
 \cZ_{N,t,M}^{\star}(u)=\frac1{\sqrt M}
 \sum_{h=0}^{\lfloor Mu\rfloor-1}
 \left\{\frac{\cA_N^\star(t+h)^2}{t}-c_d\right\},
 \qquad 0\leq u\leq1.
\end{equation}
We denote its polygonal interpolation by $\overline{\cZ}_{N,t,M}^{\star}$.

For windows comparable to the initial time, subtracting the time-dependent
profile gives the following functional limit.

\begin{corollary}[Common-time profile and Brownian time change]\label{S:cor:macro}
Assume \eqref{S:eq:conditions} and $M_N/t_N\to\rho\in[0,\infty)$, and put
$\sigma_B=2\sqrt{\mathfrak p_Bc_d}$. For either trace centring,
\begin{equation}\label{S:eq:macro}
 \frac1{\sqrt{M_N}}\sum_{h<\lfloor M_Nu\rfloor}
 \left\{\frac{\cA_N^\star(t_N+h)^2}{t_N}
                -c_d\frac{t_N+h}{t_N}\right\}
 \Longrightarrow \sigma_B\int_0^u(1+\rho v)\dd W(v).
\end{equation}
Convergence holds in $D([0,1])$ with $J_1$, and the polygonal
interpolations converge in $C([0,1])$. Jointly, the observation-time
process in Theorem~\ref{S:thm:natural} converges to $\sigma_BW$ with the
same Brownian motion. The two trace centrings have the same joint limit,
and the limiting processes are independent of each fixed low-mode vector.
The limit in \eqref{S:eq:macro} has the same law as
\[
 u\longmapsto\sigma_B\widetilde W
       \left(u+\rho u^2+\frac{\rho^2u^3}{3}\right),
\]
where $\widetilde W$ is standard Brownian motion. In particular,
\[
 \sqrt{M_N}\left\{\cF_{N,t_N,M_N}^\star
       -c_d\left(1+\frac{M_N-1}{2t_N}\right)\right\}
 \Longrightarrow
 \mathcal N\left(0,4\mathfrak p_Bc_d
                  \left(1+\rho+\frac{\rho^2}{3}\right)\right).
\]
\end{corollary}

On windows with $M_N/t_N\to0$, centring each squared observation at $c_d$
instead of the profile produces a transition at the scale
$M_N^{3/2}/t_N$.

\begin{corollary}[Drift under the common initial-time normalisation]\label{S:cor:drift}
Under \eqref{S:eq:conditions} and $M_N/t_N\to0$, for either centring,
\begin{equation}\label{S:eq:drift}
 \left\|\cZ_{N,t_N,M_N}^{\star}
   -\cZ_{N,t_N,M_N}^{\natural,\star}
   -\frac{c_d M_N^{3/2}}{2t_N}(u\mapsto u^2)\right\|_\infty
       \longrightarrow0\quad\text{in probability}.
\end{equation}
If $M_N^{3/2}/t_N\to\kappa\in[0,\infty)$, the common-time process and its interpolation converge to
\[
  u\longmapsto2\sqrt{\mathfrak p_Bc_d}\,W(u)
                           +\frac{\kappa c_d}{2}u^2.
\]
If $M_N^{3/2}/t_N\to\infty$, its terminal value tends to $+\infty$ in probability, so the sequence is not tight. The finite-$\kappa$ limit holds jointly with the same independent low-mode vector as in Theorem~\ref{S:thm:natural}.
\end{corollary}

For power windows $M_N=\lfloor t_N^\gamma\rfloor$ with $0<\gamma<1$,
the drift vanishes for $\gamma<2/3$, equals $c_du^2/2$ at $\gamma=2/3$,
and makes the common-time process non-tight for $\gamma>2/3$, under the
stated cutoff condition.

Table~\ref{tab:normalisations} compares the three cumulative statistics.
For either trace centring, the process in each row is the sum of the
indicated summands over
$0\leq h<\lfloor M_Nu\rfloor$, divided by $\sqrt{M_N}$, and
\eqref{S:eq:conditions} is assumed. Write
$Q_{N,h}=\cA_N^\star(t_N+h)^2$ and
$\sigma_B=2\sqrt{\mathfrak p_Bc_d}$.
\begin{table}[htbp]
\centering
\small
\caption{Time normalisation, centring, and the resulting process limit.}
\label{tab:normalisations}
\setlength{\tabcolsep}{4pt}
\renewcommand{\arraystretch}{1.45}
\begin{tabular}{@{}p{43mm}p{43mm}p{64mm}@{}}
\hline
Summand & Additional condition & Process limit\\
\hline
$\displaystyle\frac{Q_{N,h}}{t_N+h}-c_d$
 & None
 & $\sigma_BW(u)$\\[2pt]
$\displaystyle\frac{Q_{N,h}}{t_N}-c_d\frac{t_N+h}{t_N}$
 & $M_N/t_N\to\rho\in[0,\infty)$
 & $\displaystyle\sigma_B\int_0^u(1+\rho v)\dd W(v)$\\[2pt]
$\displaystyle\frac{Q_{N,h}}{t_N}-c_d$
 & $M_N/t_N\to0$,\newline
   $M_N^{3/2}/t_N\to\kappa\in[0,\infty)$
 & $\displaystyle\sigma_BW(u)+\frac{\kappa c_d}{2}u^2$\\
\hline
\end{tabular}
\end{table}

Centring each squared observation at its own expectation is a separate
operation from either choice of trace centring. Corollary~\ref{F:cor:expectation}
gives expectation-centred versions of the first two limits in
Table~\ref{tab:normalisations} under the additional sufficient condition
$\delta_{N,R_N}/\sqrt{M_N}\to0$. For common-time normalisation,
$M_N/t_N\to\rho<\infty$ is still assumed, but no condition forcing
$M_N^{3/2}/t_N$ to vanish is needed after expectation centring.

\subsection{Fluctuations at nearby times}
We now turn from time averages to observations at fixed time differences. The next theorem describes their joint fluctuations for a non-bipartite base. The covariance depends on the base only through its degree. The corresponding parity-dependent statement for bipartite bases is Theorem~\ref{S:thm:bipoffset}.

\begin{theorem}[Grover trace fluctuations at growing times]\label{thm:main}
Suppose $H_B$ is primitive. Fix distinct integers $h_1,\ldots,h_p$, put $H=\max_a|h_a|$, and let $t_N$ be a sequence of positive integers tending to infinity. Set $R_N=t_N+H$. If
\begin{equation}\label{eq:weakcondition}
 \frac{R_N^4q^{R_N}}{N}\longrightarrow0,
\end{equation}
then
\begin{equation}\label{eq:mainclt}
 \left(\frac{\cA_N(t_N+h_a)}{\sqrt{t_N}}\right)_{a=1}^p
 \ \xrightarrow{\mathrm d}\ (G_{h_a})_{a=1}^p,
 \qquad \E[G_hG_{h'}]=\Gamma_d(h-h'),
\end{equation}
where $(G_h)_{h\in\mathbb Z}$ is a centred stationary real Gaussian process. If in addition
\begin{equation}\label{eq:momentcondition}
 \frac{R_N^5q^{R_N}}{N}\longrightarrow0,
\end{equation}
then the convergence in \eqref{eq:mainclt} is in quadratic Wasserstein distance, and every covariance converges. For each fixed $L\geq3$, the joint vector formed by appending $(Y_{N,k})_{k=3}^{L}$ converges to $((G_{h_a})_{a=1}^{p},(Y_k)_{k=3}^{L})$, where the Gaussian field is independent of the limiting low-mode vector and $Y_{N,k}$ and $Y_k$ are defined in Section~\ref{sec:inputs}. This joint convergence is weak under \eqref{eq:weakcondition} and holds in quadratic Wasserstein distance under \eqref{eq:momentcondition}.
\end{theorem}
For every fixed $0<\varepsilon<1$, the condition
\[
 R_N\leq(1-\varepsilon)\log_qN
\]
is sufficient for both \eqref{eq:weakcondition} and \eqref{eq:momentcondition}. These sufficient conditions differ near the logarithmic cutoff; neither is claimed to be optimal. Theorem~\ref{thm:main} involves no time averaging.

\begin{corollary}[Linear connected form factor]\label{cor:ff}
Under the assumptions of Theorem~\ref{thm:main} and $t_N^5q^{t_N}/N\to0$,
\begin{equation}\label{eq:ramplaw}
 \frac{\mathcal F_N^{\mathrm c}(t_N)}{t_N}\longrightarrow
 c_d=\Gamma_d(0)
 =\frac{4}{\pi}\arcsin\left(\frac{2\sqrt{d-1}}{d}\right).
\end{equation}
More generally, the normalised connected trace covariance at fixed offsets converges to the kernel in \eqref{eq:kernel}.
\end{corollary}

The parity formulas are
\begin{equation}\label{eq:paritykernel}
 \begin{aligned}
 \Gamma_d(2m)&=-\frac{2\sin(2m\alpha_d)}{\pi m}
       &&(m\in\mathbb Z\setminus\{0\}),\\
 \Gamma_d(2m+1)&=0 &&(m\in\mathbb Z).
 \end{aligned}
\end{equation}
The even and odd subsequences of the limiting Gaussian process are therefore independent, while successive coordinates within either subsequence need not be independent. In particular, $\Gamma_d(2)<0$. These statements concern the limiting process; they do not assert exact decorrelation for finite lifts.

\subsection{Outline of the proof}
Trace identities express the observable in terms of cyclically
non-backtracking closed walks. Simple cycles supply the main random
contribution; the remaining walks are controlled by direct moment
bounds. Since the number of relevant cycle lengths grows with time,
Theorem~\ref{L:thm:poisson} approximates the whole vector up to the
terminal cutoff rather than only finitely many fixed lengths.

The analytic step is uniform in both observation times.
Lemma~\ref{S:lem:absoluteD} gives an $L^1$ first-order approximation for
$|D|e^{in\varphi}$, where $D=-i\partial_x$ and the smooth real periodic
phase has nondegenerate critical points. Theorem~\ref{S:thm:energy}
turns it into the required bilinear energy estimate. For the ideal
Poisson trace field, dividing each trace by the square root of its
observation time gives a covariance error bounded by $C_B/\sqrt{su}$
at times $s,u$; for bipartite bases the stationary comparison is made
on the even-time subsequence. The error need not be small at the first
few observations, but its cumulative effect is negligible.

Rio's coupling \cite{Rio}, with its time decay retained, replaces the
ideal Poisson trace field by a Gaussian field with the same covariance.
For the observation-time squared process, the comparison error is
$O((1+\log M)/\sqrt M)$ in expected uniform norm, uniformly over
positive initial times.
The quadratic-array theorem then gives the Brownian limit.
Cycle approximation and a first-moment remainder estimate transfer it
to the finite lift. Discrete summation by parts yields the common-time
profile and drift, while separate moment estimates give the covariance
and $L^1$ conclusions.

\subsection{Relation to earlier work}
\paragraph{Grover spectra and spectral fluctuations.}
The determinant correspondence of Konno and Sato \cite{KS}, its covering
decomposition \cite{SKKM}, and Grover trace formulae \cite{KMMS} provide
the algebraic starting point. We use the regular-graph specialisation in
\cite[Corollary 3.1]{AKM}, together with the old/new decomposition.
The finite-size study \cite{AKM} examined local spacings of the new Grover
eigenangles after removing the inherited base spectrum. Here the
observable is the full centred trace at increasing powers of the
evolution operator, followed by time averaging under the unconditioned
lift law. Oren and Smilansky \cite{OS} connect closed-walk statistics to
form factors of random regular graphs, using random-matrix assumptions
for higher-order predictions. Gnutzmann and Altland \cite{GA} study
spectral correlations of quantum graphs under energy averaging.
Our randomness is instead supplied by the covering permutations.

\paragraph{Relation to quantum graphs and periodic-orbit theory.}
The local Grover reflection $(2/d)\mathbf1\mathbf1^*-I_d$ is the Neumann--Kirchhoff vertex-scattering matrix, so $U_G$ coincides with the bond-scattering map of the equilateral Neumann quantum graph on $G$ after removal of the common propagation phase; see, for example, \cite{KottosSmilansky1999,Tanner2006,GnutzmannSmilansky2006}. Periodic-orbit expansions for quantum graphs relate form factors to closed trajectories, with the diagonal approximation producing the leading orthogonal-class ramp and correlated orbit pairs producing its first off-diagonal correction \cite{KottosSmilansky1999,BSW2002}. Trace formulae for discrete regular graphs likewise organise spectral statistics by periodic non-backtracking walks \cite{OrenGodelSmilansky2009}. The question whether a spectral form factor is self-averaging goes back, in a broader semiclassical setting, to Prange \cite{Prange1997}. The present fixed-base lift results concern logarithmic times, far below
the Heisenberg scale. They describe both trace fluctuations and temporal
averaging in a simultaneous graph-size and time limit.

\paragraph{Time averaging and other ensembles.}
Temporal self-averaging has precedents in circular unitary ensembles \cite[Section VI]{HSW} and, for a centred energy variance, in random surface covers \cite[Theorem 1.6]{Maoz}. Naud \cite[Theorem 1.1]{Naud} and Maoz \cite[Theorem 1.3 and Corollary 1.4]{Maoz} obtain variance limits and Gaussian limits for smoothed surface-cover statistics by taking the covering-degree limit before the smoothing limit. Rudnick and Wigman \cite[Theorem 1.1]{RW} show that the smooth energy variance on a typical random hyperbolic surface approaches the GOE ensemble variance, taking the large-genus limit before the energy limits. These surface models impose relations absent from the independent edge permutations used here. Here the base graph is fixed, the edge permutations are independent, and the graph size, time and averaging width vary simultaneously.

Cipolloni, Erd\H os and Schr\"oder \cite[Section 1.3.1 and Theorem 2.5]{CES} study time-dependent linear statistics for Wigner and monoparametric ensembles. Our limiting covariance is described by the same leading-order orthogonal covariance form; the precise correspondence is given in Section~\ref{sec:interpretation}. Maltsev and Osman \cite[Definitions 2.2--2.4 and Theorems 2.1--2.2]{MO} prove mesoscopic variance asymptotics for two quantum-graph ensembles on polynomial mesoscopic scales. Their first ensemble combines a uniformly sampled $d$-regular graph
with independent random edge lengths and symmetric equi-transmitting
vertex-scattering matrices. Their second uses independent Haar-distributed
vertex-scattering matrices on a complete graph. These choices of
randomisation and their polynomial mesoscopic scales differ from the
fixed-base lift with deterministic Grover scattering and logarithmic
observation times considered here.

\paragraph{Random lifts and cycle counts.}
Joint Poisson limits for fixed cycle lengths on fixed-base lifts are
due to Greenhill, Janson and Ruci\'nski \cite[Lemma 5.1]{GJR}.
Quantitative cycle approximations and their links to polynomial traces
also appear in the permutation and regular-graph models of Dumitriu,
Johnson, Pal and Paquette \cite{DJPP}, Johnson and Pal \cite{JP}, and
Johnson \cite{Johnson}. In particular, Johnson
\cite[Proposition 10]{Johnson} groups cycle indicators before applying
mean-dependent Stein factors. Our estimate uses classical Poisson
approximation \cite{Bar88,Rollin} and retains this benefit of aggregation.
The fixed-base argument combines an exact conditional forcing coupling
with a symmetric discrepancy estimate to obtain the explicit
cutoff-dependent error in Theorem~\ref{L:thm:poisson}.

A close comparison is provided by Dumitriu and Zhu
\cite[Theorems 2.10, 4.4 and 4.7]{DZ}, who use Poisson approximation of
cycle counts to study global eigenvalue fluctuations in uniformly sampled
biregular bipartite graphs, for fixed and growing degrees. Here the degree
and the base remain fixed, while the trace power and averaging window
grow with the number of sheets. Our main observable is a cumulative sum
of squared traces on a single lift. Thus the cycle-count method is shared,
but the ensemble, the varying observables, and the functional
time-averaging limit are different.

Yu \cite[Definitions 2.7--2.8 and Theorem 3.6]{Yu} proves logarithmic-length
Poisson approximation for cycle-count vectors in noisy fixed-base lifts,
including bipartite bases. At zero noise, its cycle-count laws on simple
bases agree with those considered here, with even sheet numbers in its
non-bipartite definition. Its non-bipartite statement allows cutoffs $D\leq c(H)\log n$
for some $c(H)>0$, where $n$ is the total number of vertices.
Theorem~\ref{L:thm:poisson} gives the error bound $C_BR^4q^R/N$, allowing
$R\leq(1-\varepsilon)\log_qN$ for every fixed $\varepsilon\in(0,1)$.
The comparison at the end of Section~\ref{sec:cyclepoisson} identifies the
limitation of the particular unweighted bound used in Yu's proof.
Yu's model also permits noise and half-loops, which are not treated here.

\paragraph{Moment estimates and functional limits.}
The moment analysis builds on fixed-word and permutation-cycle results \cite{Nica,LP,PZ}, remainder estimates \cite{FK}, and methods for random permutation matrices \cite{BC,CN}. The covariance scale can also be obtained from the general coefficient method of Hanany and Puder \cite[Propositions 3.7, 5.7 and 6.8]{HP}, using their multiword coefficient estimates. We give a direct proof that retains the length dependence needed for the time-dependent observables. The empirical spectral limits of Gong, Li and Liu \cite{GLL} concern a different object from these fluctuation laws.

The analysis of cumulative errors draws on Rio's Poisson--normal coupling
\cite[Corollary 2.2]{Rio}, functional Breuer--Major theory \cite{NN},
and its Poisson-noise counterpart \cite[Theorems 1.2 and 1.5]{KW}.
Theorem~\ref{F:thm:gaussian} handles the product-form covariance error of
the normalised field instead of assuming exact stationarity. Combined
with the uniform cycle and energy estimates, it yields temporal
self-averaging and Brownian cumulative errors on a single lift while
the graph size and observation window grow simultaneously.

\paragraph{Probability and path-space conventions.}
For probability laws on a common measurable space, let $d_{\mathrm{TV}}(\mu,\nu)=\sup_A|\mu(A)-\nu(A)|$, where the supremum is over measurable sets. On Euclidean and path spaces we use the Borel sigma-algebra of the stated topology. For laws on a finite-dimensional Euclidean space with finite second moments, write
\[
 \Wtwo(\mu,\nu)=\inf\{(\E[\|V-W\|_2^2])^{1/2}:\Law(V)=\mu,\ \Law(W)=\nu\},
\]
where the infimum is over couplings. On a finite-dimensional Euclidean space, the bounded Lipschitz distance $\BL$ is defined using test functions bounded in absolute value by one and with Lipschitz constant at most one with respect to the $\ell^1$ norm. Here and below, $\Law(V)$ denotes the probability law of a random variable $V$, and $\mathcal N(m,v)$ denotes the real normal law with mean $m$ and variance $v$.

We write $D([0,1])$ for the space of real-valued right-continuous functions with left limits, equipped with the Skorokhod $J_1$ topology, and $C([0,1])$ for the space of real-valued continuous functions, equipped with the uniform norm. For a path $x$, $\|x\|_\infty=\sup_{u\in[0,1]}|x(u)|$. All times and degrees indexing powers of $U_G$, Chebyshev polynomials or cycle counts are integers. The auxiliary immigration--death chain uses a continuous time parameter. Logarithms without a subscript are natural logarithms.

\paragraph{Organisation.}
Section~\ref{sec:inputs} gives the spectral representation.
Sections~\ref{sec:cycles}--\ref{sec:cyclepoisson} establish the cycle approximation,
and Section~\ref{sec:moments} supplies moment bounds. Section~\ref{sec:coefficients}
proves the bilinear energy estimate. Section~\ref{sec:poisson} organises the period,
covariance correction and ideal field; Section~\ref{sec:transfer} transfers
fixed-offset limits. Section~\ref{sec:timeaverages} gives explicit centring,
a common transfer estimate and finite-size bounds. Section~\ref{sec:functional}
proves the temporal limit theorems and their common-time consequences.
Section~\ref{sec:interpretation} discusses their spectral meaning and scope.

\section{Spectral representation and notation}\label{sec:inputs}
We use the counting-measure inner product on $\ell^2(V\times[N])$. The old subspace consists of functions constant on each fibre. Its orthogonal complement is
\[
 \mathcal H_N^{\new}=\{f:\sum_{a=1}^N f(v,a)=0\text{ for every }v\in V\}.
\]
The subspace of fibre-constant functions is invariant under the symmetric adjacency matrix, and hence so is its orthogonal complement. Thus
\[
 A_{G_N}\simeq A_B\oplus A_N^{\new},
\]
where $A_N^{\new}$ is the restriction to $\mathcal H_N^{\new}$; see \cite[Proposition 2.2]{AKM}. The old block is deterministic, so the centred polynomial traces of the full matrix and its new restriction agree.

The Chebyshev polynomials of the first kind satisfy $T_0(x)=1$, $T_1(x)=x$ and $T_{k+1}(x)=2xT_k(x)-T_{k-1}(x)$. Define
\begin{equation}\label{eq:Ymodes}
 Y_{N,k}=2\left\{\Tr T_k\!\left(\frac{A_N^{\new}}{2\sqrt q}\right)
       -\E\left[\Tr T_k\!\left(\frac{A_N^{\new}}{2\sqrt q}\right)\right]\right\}.
\end{equation}
The modes of orders zero, one and two are identically zero because the corresponding new-sector quantities are deterministic: $\dim\mathcal H_N^{\new}=n_0(N-1)$, $\Tr A_N^{\new}=0$, and $\Tr\![(A_N^{\new})^2]=n_0(N-1)d$. All eigenvalues of the new restriction are retained.

Write
\begin{equation}\label{L:eq:basecounts}
 b_j=\Tr H_B^j,\qquad \lambda_j=\frac{b_j}{2j},\qquad
 \mu_{N,j}=\E[C_{N,j}],
\end{equation}
where $C_{N,j}$ is the number of unrooted, unoriented simple $j$-cycles in $G_N$. In particular $b_1=b_2=0$. Let $Z_j$, $j\geq3$, be independent Poisson variables of means $\lambda_j$; a zero mean corresponds to a variable that is identically zero. Set
\begin{equation}\label{eq:idealY}
\begin{aligned}
 Y_k&=2q^{-k/2}\sum_{j\mid k,\ j\geq3}j(Z_j-\lambda_j),\qquad Y_1=Y_2=0,\\
 \Lambda_B(k,\ell)&=\Cov(Y_k,Y_\ell)
     =2q^{-(k+\ell)/2}\sum_{j\mid\gcd(k,\ell)}j b_j.
\end{aligned}
\end{equation}
The covariance formula follows from independence and $\Var(Z_j)=\lambda_j$. The finite-lift approximation of these variables is proved in Sections~\ref{sec:cyclepoisson} and \ref{sec:moments}.

The determinant relation of \cite{KS} and the covering decomposition of \cite{SKKM}, in the specialisation to $d$-regular graphs recorded in \cite[Corollary 3.1]{AKM}, give, for every integer $s\geq1$,
\[
 \Tr U_{G_N}^{\,s}-\Tr U_B^{\,s}
  =(N-1)(m_0-n_0)(1+(-1)^s)
     +2\Tr T_s(A_N^{\mathrm{new}}/d).
\]
The identity also holds when the lift is disconnected, by applying the determinant relation componentwise. Define the real coefficients $a_{s,k}$ by
\begin{equation}\label{eq:chebcoeff}
 T_s(r_dx)=\sum_{k=0}^{s}a_{s,k}T_k(x),
 \qquad a_{s,k}=0\quad(k>s).
\end{equation}
Centring this identity gives the exact formula
\begin{equation}\label{eq:tracecoef}
 \cA_N(s)=\sum_{k=3}^{s}a_{s,k}Y_{N,k}.
\end{equation}
This centred identity retains all eigenvalues; only deterministic trace contributions cancel. We will compare \eqref{eq:tracecoef} with
\begin{equation}\label{eq:idealsignal}
 P_s=\sum_{k=3}^{s}a_{s,k}Y_k.
\end{equation}

\section{Short cycles and conditional exploration}\label{sec:cycles}
\subsection{Admissible cycles and exact means}
Consider the deterministic graph on $V\times[N]$ containing all edges between fibres over adjacent vertices of $B$. A simple cycle in this graph is \emph{admissible} if, for every base edge, the pairs it requires form a partial bijection between the two incident fibres. For such a simple cycle, admissibility is equivalent to its projection being cyclically non-backtracking. Indeed, two cycle edges at a lifted vertex cannot use the same base edge, and distinct lifted vertices in a given fibre give distinct rows or columns of the corresponding partial permutation.

For a cycle $\alpha$, let $|\alpha|$ be its number of edges and $V(\alpha)$ its vertex set. Let $\Ga_{N,j}$ be the finite set of admissible $j$-cycles, with cyclic shifts and reversal identified. Write
\[
 I_\alpha=\ind_{\{\alpha\subseteq G_N\}},\qquad
 p_\alpha=\Prob(I_\alpha=1),\qquad
 C_{N,j}=\sum_{\alpha\in\Ga_{N,j}}I_\alpha.
\]
If $r_e(\alpha)$ is the number of edges of $\alpha$ projecting onto the unoriented edge $e$, independence of the base-edge permutations gives
\begin{equation}\label{L:eq:cycleprob}
 p_\alpha=\prod_{e\in E}\frac1{(N)_{r_e(\alpha)}},
 \qquad (N)_r=N(N-1)\cdots(N-r+1),\quad (N)_0=1.
\end{equation}

\begin{proposition}[Exact mean and one-sided error bound]\label{L:prop:mean}
Let $\mathcal W_j(B)$ be the set of rooted oriented cyclically non-backtracking closed base walks of length $j$. For $w\in\mathcal W_j(B)$, let $a_v(w)$ count occurrences of $v$ among its $j$ vertex positions, and let $a_e(w)$ count traversals of the unoriented edge $e$. If $N\geq j$, then
\begin{equation}\label{L:eq:exactmean}
 \mu_{N,j}=\frac1{2j}\sum_{w\in\mathcal W_j(B)}
              \frac{\prod_v(N)_{a_v(w)}}{\prod_e(N)_{a_e(w)}}.
\end{equation}
Moreover,
\begin{equation}\label{L:eq:meanerror}
 0\leq\lambda_j-\mu_{N,j}
   \leq\frac{j(j-1)}{2N}\lambda_j.
\end{equation}
\end{proposition}
\begin{proof}
For a fixed base walk, distinct sheet labels must be assigned to the vertex occurrences in each fibre. There are $\prod_v(N)_{a_v}$ such assignments. The specified edges then occur with probability $\prod_e(N)_{a_e}^{-1}$. Every lifted simple cycle is counted exactly $2j$ times, including cycles with a periodic base projection. This proves \eqref{L:eq:exactmean}.

Put $F(r)=\log(N)_r$ for $0\leq r\leq N$. Its successive increments decrease, so $F$ is discretely concave and $F(0)=0$. At every vertex with $a_v>0$, cyclic non-backtracking gives
\[
 0\leq a_e\leq a_v\quad(e\ni v),\qquad
 \sum_{e\ni v}a_e=2a_v.
\]
Consequently, $F(a_e)\geq(a_e/a_v)F(a_v)$, and hence
\[
 \sum_{e\ni v}F(a_e)\geq2F(a_v).
\]
Sum over $v$, noting that each base edge is counted twice. It follows that $\sum_eF(a_e)\geq\sum_vF(a_v)$, so each ratio in \eqref{L:eq:exactmean} is at most one.

For the lower bound, use $\sum_ea_e=\sum_va_v=j$ to obtain
\[
 \frac{\prod_v(N)_{a_v}}{\prod_e(N)_{a_e}}
 \geq\prod_v\prod_{r=0}^{a_v-1}\left(1-\frac rN\right)
 \geq1-\frac1{2N}\sum_va_v(a_v-1)
 \geq1-\frac{j(j-1)}{2N}.
\]
Summing both bounds over the base walks, applying the normalising factor, and using $|\mathcal W_j(B)|=b_j$ proves \eqref{L:eq:meanerror}.
\end{proof}

Since each row of $H_B$ has sum $q$, we have $b_j\leq2m_0q^j$. Thus, when $R\leq N$,
\begin{equation}\label{L:eq:geometricmeans}
 \sum_{j=3}^R j\mu_{N,j}\leq C_Bq^R,\qquad
 \sum_{j=3}^R j^2\mu_{N,j}\leq C_BRq^R,\qquad
 \sum_{j=3}^R j q^j\leq C_BRq^R.
\end{equation}
The first two estimates also hold with $\lambda_j$ in place of $\mu_{N,j}$.

\subsection{Forcing a cycle without destroying compatible cycles}
Two admissible cycles are \emph{compatible} if their combined requirements form a partial bijection over every base edge. Distinct admissible cycles need not be compatible; vertex-disjoint cycles are always compatible.

\begin{lemma}[Forcing a partial permutation]\label{L:lem:forceperm}
Let $\pi$ be uniformly distributed over the permutations satisfying $r-1$ prescribed pairs with distinct sources and targets. Let $a$ be a new source and $b$ a new target. Then
\begin{equation}\label{L:eq:force}
 \pi'=(b\ \pi(a))\circ\pi
\end{equation}
is uniformly distributed over the permutations satisfying all $r$ pairs, where the transposition is interpreted as the identity if $\pi(a)=b$. Previously prescribed pairs remain unchanged.
\end{lemma}
\begin{proof}
Neither $b$ nor $\pi(a)$ is a previously prescribed target, so the earlier constraints are preserved. Fix a permutation $\tau$ satisfying all $r$ pairs. Its preimages under \eqref{L:eq:force} are exactly
\[
 \pi=(b\ c)\circ\tau,
\]
where $c$ ranges over the $N-r+1$ targets not used by the first $r-1$ pairs. Each such $\pi$ satisfies the earlier constraints and maps to $\tau$, and the preimages are distinct. The number of preimages is independent of $\tau$, proving uniformity.
\end{proof}

\begin{lemma}[Cycle-forcing coupling]\label{L:lem:forcecycle}
Fix $\alpha\in\Ga_{N,i}$. There is a graph $G_N^\alpha$ on the same probability space as $G_N$ such that
\begin{equation}\label{L:eq:conditionallaw}
 \Law(G_N^\alpha)=\Law(G_N\mid I_\alpha=1).
\end{equation}
Writing $J_{\beta,\alpha}=\ind_{\{\beta\subseteq G_N^\alpha\}}$, one has
\begin{equation}\label{L:eq:compatiblepreserve}
 J_{\beta,\alpha}\geq I_\beta
 \quad\hbox{whenever $\beta$ is compatible with $\alpha$}.
\end{equation}
For incompatible $\beta$, $J_{\beta,\alpha}=0$.
\end{lemma}
\begin{proof}
For each base edge, order the finitely many pairs required by $\alpha$ and impose them successively by Lemma~\ref{L:lem:forceperm}. The edge permutations are treated separately, so they remain independent, with their respective uniform conditional laws. Since $\{I_\alpha=1\}$ is exactly the intersection of those edge-wise constraints, this gives \eqref{L:eq:conditionallaw}.

To prove \eqref{L:eq:compatiblepreserve}, suppose that $\beta$ is present at a given stage and consider one forcing step $a\mapsto b$. A nontrivial step changes only row $a$ and row $c=\pi^{-1}(b)$. If $\beta$ prescribed a value at $a$, compatibility would force that value to be $b$, and the step would be trivial. If $\beta$ prescribed the current value $b$ at $c\ne a$, its constraints would conflict with the required pair $a\mapsto b$. Thus no edge of a present compatible $\beta$ is removed by a nontrivial step. Induction over the forcing steps proves preservation. Incompatibility makes simultaneous occurrence impossible, giving the final assertion.
\end{proof}

Lemma~\ref{L:lem:forcecycle} will be used separately for compatible and incompatible cycle pairs.

\begin{lemma}[Vertex-disjoint cycles]\label{L:lem:disjoint}
Suppose that $|\alpha|=i$, $|\beta|=j$, and $V(\alpha)\cap V(\beta)=\varnothing$. For $i,j\leq R$ and $N\geq8R^2$,
\begin{equation}\label{L:eq:disjointchange}
 \E[|I_\beta-J_{\beta,\alpha}|]
 \leq \frac{4ij}{N}p_\beta.
\end{equation}
\end{lemma}
\begin{proof}
Put $a_e=r_e(\alpha)$ and $c_e=r_e(\beta)$. The row sets and target sets used by the two cycles are disjoint over every base edge. Therefore
\begin{equation}\label{L:eq:ratio}
 \frac{\Prob(I_\beta=1\mid I_\alpha=1)}{p_\beta}
 =\prod_e\frac{(N)_{a_e}(N)_{c_e}}{(N)_{a_e+c_e}}.
\end{equation}
Every factor is at least one. Its logarithm is
\[
 \sum_{r=0}^{c_e-1}\log\left(1+\frac{a_e}{N-a_e-r}\right)
 \leq\frac{a_ec_e}{N-i-j}.
\]
The logarithm of \eqref{L:eq:ratio} is thus at most $2ij/N\leq1/4$. Using $e^x-1\leq2x$ on this range and the pathwise inequality \eqref{L:eq:compatiblepreserve} gives \eqref{L:eq:disjointchange}.
\end{proof}

\subsection{Cycles intersecting a forced cycle}\label{sec:deferred}
\paragraph{Deferred decisions.}
Fix the deterministic input of an exploration algorithm whose next query and stopping decision depend only on the answers received so far. Queries may evaluate a permutation or its inverse. For a transcript $h$ that can occur, let $F_e(h)$ be the partial bijection revealed over base edge $e$. Then
\[
 \{\text{transcript}=h\}
   =\bigcap_{e\in E}\{\sigma_e\supseteq F_e(h)\}.
\]
Induction over the successive queries shows that every completion on the right reproduces the transcript. Thus, conditional on $h$, the edge permutations are independent, uniformly distributed completions of the respective partial bijections. The same conditional-completion description holds after exposing any initial compatible set of pairs. We also use it for the pair exploration in Section~\ref{sec:moments}.

\begin{proposition}[Conditional overlap bound]\label{L:prop:overlap}
Fix an admissible cycle $\alpha$ of length $i$. If $N>i+j$ and $j\geq3$, then
\begin{equation}\label{L:eq:overlapexplicit}
 \E\!\left[
  \sum_{\substack{\beta\in\Ga_{N,j},\ \beta\ne\alpha\\
                   V(\beta)\cap V(\alpha)\ne\varnothing}} I_\beta
       \,\middle|\, I_\alpha=1\right]
 \leq\frac{i^2j\,d q^{j-1}}{N-i-j}.
\end{equation}
In particular, if $i,j\leq R$ and $N\geq8R^2$, the right-hand side is at most $C_Bi^2j q^j/N$.
\end{proposition}
\begin{proof}
Expose all pairs belonging to $\alpha$. Under the conditional law, each edge permutation is uniformly distributed among the completions of its exposed partial bijection.

Choose a starting vertex $x\in V(\alpha)$ and a prescribed non-backtracking base walk of length $j$ beginning at the projection of $x$. The base walk has a unique lifted trajectory, which we reveal one step at a time. If the required matching pair is already known, follow the determined step. When it is unknown, expose it; its other endpoint is uniform among the still unused labels in the appropriate fibre. This statement holds for traversal in either orientation, by using the partial inverse permutation when necessary.

Throughout this exploration, at most $i+j$ pairs of any one base-edge permutation have been specified. Thus the conditional probability that a fresh exposure lands in $V(\alpha)$ is at most
\begin{equation}\label{L:eq:hitprob}
 \frac{i}{N-i-j}.
\end{equation}
Define the event at a step to be empty when no fresh pair is exposed. Summing the conditional probabilities over at most $j$ steps bounds the probability of at least one such fresh hit by $ji/(N-i-j)$.

If the lifted trajectory is a simple cycle $\beta\ne\alpha$ meeting $\alpha$, it necessarily makes a fresh hit on $V(\alpha)$. To see this, a distinct simple cycle cannot consist solely of edges of $\alpha$. If it traverses an edge outside $\alpha$ with both endpoints on $\alpha$, that traversal is already a fresh hit. Otherwise, it leaves $V(\alpha)$ and later returns for the first time. The return edge is not an edge of $\alpha$, and the simplicity of $\beta$ prevents the return edge from having been traversed earlier. It is therefore newly exposed at the return step and lands in $V(\alpha)$. This use of a fresh return edge relies on the simplicity of the second cycle.

There are $i$ choices of $x$ and at most $dq^{j-1}$ choices of the base walk. Every cycle counted on the left-hand side of \eqref{L:eq:overlapexplicit} gives at least one such rooted trajectory (in fact $2|V(\alpha)\cap V(\beta)|$ of them). Summing their probability bounds gives the asserted upper bound
\[
 i\,dq^{j-1}\frac{ji}{N-i-j}.\qedhere
\]
\end{proof}

\begin{lemma}[Unconditional overlap bound]\label{L:lem:unconditional}
For a fixed set $A\subset V\times[N]$ of size $i$,
\begin{equation}\label{L:eq:overlapunconditional}
 \E\left[\sum_{\substack{\beta\in\Ga_{N,j}\\V(\beta)\cap A\ne\varnothing}}I_\beta\right]
 \leq\frac{ij}{N}\mu_{N,j}.
\end{equation}
\end{lemma}
\begin{proof}
Let $\mathcal T_{N,j}(x)$ count the present $j$-cycles containing $x$. The lift law is invariant under independent relabellings of the fibres: an edge permutation is then multiplied on the left and right by fixed permutations, which preserves the independent uniform law. Consequently $\E[\mathcal T_{N,j}(x)]$ is constant within each fibre. For $x$ above $v$,
\[
 \E[\mathcal T_{N,j}(x)]=\frac1N\E\left[\sum_{y\in\{v\}\times[N]}\mathcal T_{N,j}(y)\right]
 \leq\frac jN\mu_{N,j}.
\]
Summing over $x\in A$ proves \eqref{L:eq:overlapunconditional}.
\end{proof}

\section{Length-aggregated Poisson approximation}\label{sec:cyclepoisson}
We use the classical immigration--death construction for the Poisson Stein equation \cite{Bar88,DJPP}, together with a coordinatewise bound in terms of the aggregated means. The second-difference estimate below includes diagonal coordinate pairs.

\begin{theorem}[Uniform cycle approximation]\label{L:thm:poisson}
There is a constant $C_B<\infty$ such that, for integers $R\geq3$ and $N\geq8R^2$,
\begin{equation}\label{L:eq:poissonmain}
 \TV\!\left(\Law(C_{N,3},\ldots,C_{N,R}),
       \bigotimes_{j=3}^R\Poi(\lambda_j)\right)
 \leq C_B\frac{R^4q^R}{N}.
\end{equation}
Primitivity of $H_B$ is not required for this assertion.
\end{theorem}

\subsection{Mean-sensitive second differences}
Let $Z$ have independent coordinates $Z_i\sim\Poi(\mu_i)$, for $i$ in a finite set $J$, with all $\mu_i>0$. The state space is $\mathbb Z_{\geq0}^{J}$, and $e_i$ is its $i$th coordinate unit vector. Set $\Delta_i h(x)=h(x+e_i)-h(x)$. Suprema of functions on this state space are taken over all $x\in\mathbb Z_{\geq0}^{J}$. Consider the generator
\[
 \mathcal A h(x)=\sum_i\mu_i\Delta_i h(x)
       +\sum_i x_i\{h(x-e_i)-h(x)\}.
\]
Terms with $x_i=0$ in the second sum are zero. For an indicator test function $f$, a solution of $\mathcal A h=f-\E[f(Z)]$ is
\begin{equation}\label{L:eq:semigroup}
 h(x)=-\int_0^\infty\{\E[f(X_t^x)]-\E[f(Z)]\}\,dt,
\end{equation}
where $X_t^x$ is the immigration--death chain started at $X_0^x=x$, with coordinate immigration rates $\mu_i$ and individual death rate one. A coupling with the stationary chain bounds the integrand by $(\|x\|_1+\sum_i\mu_i)e^{-t}$, which justifies the integral.

Write $\log^+a=\max(0,\log a)$ and define, for $a>0$,
\begin{equation}\label{L:eq:rho}
 \rho(a)=\min\left\{1,\frac{2+\log^+a}{a}\right\}.
\end{equation}
\begin{lemma}[Coordinatewise second differences]\label{L:lem:hessian}
For the solution \eqref{L:eq:semigroup} and every pair $i,j$, including $i=j$,
\begin{equation}\label{L:eq:hessian}
 \|\Delta_i\Delta_j h\|_\infty
       \leq \rho(\sqrt{\mu_i\mu_j}).
\end{equation}
The bound is independent of $|J|$ and of $x$.
\end{lemma}
\begin{proof}
Couple chains started at $x$, $x+e_i$, $x+e_j$, and $x+e_i+e_j$ using two extra labelled particles with independent unit-rate lifetimes. The same coupling applies when both particles are in the same coordinate. The bound $2e^{-2t}$ justifies taking the second difference under the integral. The terms cancel unless both extra particles survive, giving
\[
 \Delta_i\Delta_j h(x)
   =-\int_0^\infty e^{-2t}\E[\Delta_i\Delta_j f(X_t^x)]\,dt.
\]
At time $t$, the immigration part consists of independent Poisson variables with means $\mu_i u$, $u=1-e^{-t}$; it is independent of the surviving initial particles.

For $v>0$, let $p_v$ be the Poisson probability mass function on $\mathbb Z_{\geq0}$, and let $\delta p_v(n)=p_v(n-1)-p_v(n)$, with $p_v(-1)=0$. Its total mass is zero and
\begin{equation}\label{L:eq:shift}
 \|\delta p_v\|_{\ell^1}
 =\E\left[\left|\frac{Z_v}{v}-1\right|\right]\leq v^{-1/2},
 \qquad Z_v\sim\Poi(v).
\end{equation}
For $i\ne j$, the double difference is an integral against the tensor product of these two signed measures. Conditioning on the survivors and other coordinates does not increase its norm. The integral of an indicator against a zero-mass signed measure is bounded in absolute value by half its $\ell^1$ norm. Therefore
\[
 |\E[\Delta_i\Delta_j f(X_t^x)]|
       \leq\frac{1}{2u\sqrt{\mu_i\mu_j}}.
\]
For $i=j$, split the immigration variable of mean $\mu_i u$ into the sum of two independent Poisson variables of mean $\mu_i u/2$. Its second signed difference is
$\delta p_{\mu_i u/2}*\delta p_{\mu_i u/2}$.
The convolution inequality and \eqref{L:eq:shift} bound its $\ell^1$ norm by $2/(\mu_i u)$, giving the bound $1/(\mu_i u)$ for indicator test functions.

In both cases, the elementary bound $|\Delta_i\Delta_j f|\leq2$ also holds. Put $a=\sqrt{\mu_i\mu_j}$. We obtain
\begin{align}\label{L:eq:integral}
 \|\Delta_i\Delta_jh\|_\infty
 &\leq\int_0^\infty e^{-2t}
             \min\left\{2,\frac{1}{a(1-e^{-t})}\right\}\,dt\notag\\
 &=\int_0^1(1-u)\min\left\{2,\frac1{au}\right\}\,du
 =\begin{cases}
 1,&a\leq1/2,\\
 (\log(2a))/a+1/(4a^2),&a>1/2.
 \end{cases}
\end{align}
This is at most one, and at most $(2+\log^+a)/a$. Equation \eqref{L:eq:hessian} follows. All intermediate bounds are uniform over the conditioned values of the survivors, so the bound is uniform in the number of coordinates.
\end{proof}

Logarithmic dependence on the means is already present in \cite[Lemma 3]{Bar88}, and \cite[Section 3]{Rollin} discusses why logarithmic factors in multivariate Stein estimates cannot be discarded without further structure. We make no claim that the polynomial loss in our graph application is optimal.

\subsection{Aggregation before estimation}
\begin{lemma}[Aggregation before estimation]\label{L:lem:grouped}
Let $\Gamma=\bigsqcup_{i\in J}\Gamma_i$ index a finite Bernoulli family $I=(I_\alpha)_{\alpha\in\Gamma}$, with $p_\alpha=\Prob(I_\alpha=1)>0$. Put $C_i=\sum_{\alpha\in\Gamma_i}I_\alpha$ and $C=(C_i)_{i\in J}$. Discard empty index classes, so that $\mu_i=\E[C_i]>0$ for every retained index. For every $\alpha$, suppose that $J^\alpha=(J_{\beta,\alpha})_\beta$ is coupled with $I$ and has the conditional law of this Bernoulli family given $I_\alpha=1$. Define
\[
 S_{ij}=\sum_{\alpha\in\Gamma_i}p_\alpha
       \sum_{\substack{\beta\in\Gamma_j\\\beta\ne\alpha}}
                    \E[|I_\beta-J_{\beta,\alpha}|],
 \qquad M_{ij}=\rho(\sqrt{\mu_i\mu_j}).
\]
Then
\begin{equation}\label{L:eq:grouped}
 \TV\left(\Law(C),\bigotimes_i\Poi(\mu_i)\right)
 \leq\sum_i M_{ii}\sum_{\alpha\in\Gamma_i}p_\alpha^2
                      +\sum_{i,j}M_{ij}S_{ij}.
\end{equation}
\end{lemma}
\begin{proof}
Set $C_j^\alpha=\sum_{\beta\in\Gamma_j}J_{\beta,\alpha}$. For $\alpha\in\Gamma_i$, the coordinate $C_i^\alpha$ is at least one. The conditional law gives the exact Stein identity
\[
 \E[\mathcal A h(C)]=\sum_i\sum_{\alpha\in\Gamma_i}p_\alpha
          \E[\Delta_i h(C)-\Delta_i h(C^\alpha-e_i)].
\]
Telescoping between the two nonnegative lattice points and applying Lemma~\ref{L:lem:hessian} yields
\[
 |\Delta_i h(x)-\Delta_i h(y)|
       \leq\sum_j M_{ij}|x_j-y_j|.
\]
For $j\ne i$, the coordinate difference is bounded by the sum of the corresponding $|I_\beta-J_{\beta,\alpha}|$. In coordinate $i$, use $J_{\alpha,\alpha}=1$ to obtain
\[
 |C_i-(C_i^\alpha-1)|
 \leq I_\alpha+
       \sum_{\substack{\beta\in\Gamma_i\\\beta\ne\alpha}}
                           |I_\beta-J_{\beta,\alpha}|.
\]
Taking expectations gives the self term $p_\alpha^2$, and the other terms are $S_{ij}$. Taking the supremum over indicator test functions proves \eqref{L:eq:grouped}. The factors $M_{ij}$ depend on the aggregated means $\mu_i$, instead of the individual cycle probabilities $p_\alpha$.
\end{proof}

\subsection{Symmetric cycle discrepancies}
For the remainder of this section, write $\Gamma_i=\Ga_{N,i}$, $C_i=C_{N,i}$ and $\mu_i=\mu_{N,i}$. Use the couplings of Lemma~\ref{L:lem:forcecycle} and let $S_{ij}$ be as in Lemma~\ref{L:lem:grouped}.

\begin{lemma}[Pairwise symmetry]\label{L:lem:symmetry}
For distinct admissible cycles,
\begin{equation}\label{L:eq:absCov}
 p_\alpha\E[|I_\beta-J_{\beta,\alpha}|]
      =|\Cov(I_\alpha,I_\beta)|.
\end{equation}
Consequently $S_{ij}=S_{ji}$.
\end{lemma}
\begin{proof}
For compatible cycles, preservation in Lemma~\ref{L:lem:forcecycle} gives
\[
 p_\alpha\E[|I_\beta-J_{\beta,\alpha}|]
 =p_\alpha\{\Prob(I_\beta=1\mid I_\alpha=1)-p_\beta\}
 =\Cov(I_\alpha,I_\beta)\geq0.
\]
For incompatible cycles, $J_{\beta,\alpha}=0$ and both sides of \eqref{L:eq:absCov} equal $p_\alpha p_\beta$. Symmetry of the absolute covariance and summation over the cycles prove the assertion, including $i=j$. The two cases account for both signs of the covariance.
\end{proof}

\begin{lemma}[Bidirectional discrepancy estimate]\label{L:lem:balance}
If $i,j\leq R$, $N\geq8R^2$ and $\mu_i\mu_j>0$, then
\begin{equation}\label{L:eq:balance}
 \frac{S_{ij}}{\sqrt{\mu_i\mu_j}}
 \leq\left(5m_0+\frac{2d}{q}\right)
       \frac{(ij)^{3/2}q^{(i+j)/2}}{N}.
\end{equation}
\end{lemma}
\begin{proof}
Lemma~\ref{L:lem:disjoint} bounds the vertex-disjoint contribution by $4ij\mu_i\mu_j/N$. For distinct intersecting cycles, use
\[
 \E[|I_\beta-J_{\beta,\alpha}|]
 \leq p_\beta+\Prob(I_\beta=1\mid I_\alpha=1).
\]
Lemma~\ref{L:lem:unconditional} bounds the sum of the unconditional terms. Proposition~\ref{L:prop:overlap}, with $N-i-j\geq N/2$, bounds the conditional terms. Together,
\begin{equation}\label{L:eq:directional}
 S_{ij}\leq\frac{5ij\mu_i\mu_j}{N}
          +\frac{2d}{q}\frac{i^2j\mu_i q^j}{N}.
\end{equation}
By Lemma~\ref{L:lem:symmetry}, \eqref{L:eq:directional} also holds with $i,j$ interchanged. Subtracting the common first term from both upper bounds allows us to replace the numerator of the second term by
\[
 \min\{i^2j\mu_iq^j,j^2i\mu_jq^i\}
 \leq (ij)^{3/2}\sqrt{\mu_i\mu_j}\,q^{(i+j)/2}.
\]
Proposition~\ref{L:prop:mean} and $b_i\leq2m_0q^i$ give $\mu_i\leq m_0q^i/i$. Thus the first numerator, divided by $\sqrt{\mu_i\mu_j}$, is at most $5m_0\sqrt{ij}\,q^{(i+j)/2}$. Replacing $\sqrt{ij}$ by the larger quantity $(ij)^{3/2}$ proves \eqref{L:eq:balance}. No uniform positive lower bound on the nonzero aggregated means is needed.
\end{proof}

\begin{proof}[Proof of Theorem~\ref{L:thm:poisson}]
If $\mu_i=0$, then $C_i=0$ almost surely. Equation \eqref{L:eq:meanerror} and $N\geq8R^2$ imply $\mu_i\geq(15/16)\lambda_i$, and hence $\lambda_i=0$. These coordinates can be removed without error. Let $J$ consist of the remaining indices; if $J$ is empty, the theorem holds with zero distance.

Set
\[
 L_R=2+\log m_0+R\log q.
\]
Since $\mu_i\leq m_0q^R$, the smoothing factor satisfies
\begin{equation}\label{L:eq:factor}
 M_{ij}\leq\frac{L_R}{\sqrt{\mu_i\mu_j}},
 \qquad L_R\leq C_B R.
\end{equation}
Apply Lemma~\ref{L:lem:grouped} with the exact forcing couplings.
For $\alpha\in\Gamma_i$, $p_\alpha\leq(N-R)^{-i}\leq(N-R)^{-3}$, so
\begin{equation}\label{L:eq:self}
 \sum_{i\in J}M_{ii}\sum_{\alpha\in\Gamma_i}p_\alpha^2
       \leq\frac{L_R|J|}{(N-R)^3}
       \leq\frac{L_R R}{(N-R)^3}.
\end{equation}
For the distinct-pair term, Lemma~\ref{L:lem:balance} and \eqref{L:eq:factor} give
\begin{align}\label{L:eq:sum}
 \sum_{i,j\in J}M_{ij}S_{ij}
 &\leq\frac{L_R(5m_0+2d/q)}{N}
             \left(\sum_{i=3}^R i^{3/2}q^{i/2}\right)^2\notag\\
 &\leq\frac{L_R(5m_0+2d/q)R^3q^R}
                  {(1-q^{-1/2})^2N}
 \leq C_B\frac{R^4q^R}{N}.
\end{align}
We have thus compared the vector of cycle counts with independent Poisson variables of means $\mu_i$. To replace these means by $\lambda_i$, add independent Poisson increments with means $\lambda_i-\mu_i$. The probability that at least one increment is nonzero is at most
\begin{equation}\label{L:eq:meanchange}
 \sum_{i=3}^R(\lambda_i-\mu_i)
 \leq\frac{m_0}{2N}\sum_{i=3}^R(i-1)q^i
 \leq\frac{m_0 Rq^R}{2(1-q^{-1})N}.
\end{equation}
Combining \eqref{L:eq:self}--\eqref{L:eq:meanchange} proves \eqref{L:eq:poissonmain}. In particular, an explicit, unoptimised upper bound is
\[
 \frac{L_RR}{(N-R)^3}
 +\frac{L_R(5m_0+2d/q)R^3q^R}{(1-q^{-1/2})^2N}
 +\frac{m_0 Rq^R}{2(1-q^{-1})N}.
\]
All denominators are positive under the stated hypothesis.
\end{proof}

\paragraph{Comparison with an unweighted neighbourhood estimate.}
We compare with \cite[Theorem 6.1 and Section 6, p.~28]{Yu} at zero
noise on a simple non-bipartite base. Its parameters are $H=B$, $m=N$,
$n=n_0N$ and $D=R$, with $N$ even as in its Definition~2.7. The random
fibre partition in its Remark~2.10 has no effect on cycle-count laws.
A positive-probability local configuration is then an admissible cycle.
Its set $S_\alpha$ of stubs consists of its vertex--base-edge incidences
and has cardinality $2|\alpha|$. There are $n_0Nd$ stubs. The neighbourhoods
used there satisfy $\beta\in\mathcal N_\alpha$ exactly when
$S_\alpha\cap S_\beta\ne\varnothing$. Denote the first error term of
\cite[Theorem 6.1]{Yu} by
\[
 \mathfrak b_1=\sum_\alpha\sum_{\beta\in\mathcal N_\alpha}p_\alpha p_\beta.
\]
Retain only length-$R$ cycles and put $u_x=\sum_{\alpha\in\Ga_{N,R}:x\in S_\alpha}p_\alpha$. Double counting and Cauchy--Schwarz give
\[
 \sum_xu_x=2R\mu_{N,R},\qquad
 2R\mathfrak b_1\geq\sum_xu_x^2
       \geq\frac{(2R\mu_{N,R})^2}{n_0Nd}.
\]
The first inequality uses $|S_\alpha\cap S_\beta|\leq2R$ and includes pairs of identical cycles. If $H_B$ is primitive, $R\to\infty$ and $R^2/N\to0$, Proposition~\ref{L:prop:mean} and $b_R/q^R\to1$ imply
\[
 \mu_{N,R}\sim\frac{q^R}{2R},\qquad
 \mathfrak b_1\geq c_B\frac{q^{2R}}{NR}
\]
for some $c_B>0$ and all sufficiently large $N$. For $R=\lfloor\xi\log_qN\rfloor$ with $1/2<\xi<1$, this lower bound diverges, whereas $R^4q^R/N\to0$. Hence that particular unweighted right-hand side $2(2\mathfrak b_1+2\mathfrak b_2+\mathfrak b_3)$ in \cite[Theorem 6.1]{Yu} is not informative in this range with those neighbourhoods. The nonnegative quantities $\mathfrak b_2,\mathfrak b_3$ are the other errors in that cited bound. This obstruction concerns the specified upper bound, not the true total-variation distance or the validity of Yu's theorem; a modified proof may give a different range. Length aggregation retains the mean-dependent Stein factors and yields the stated cutoff range.

\section{Non-backtracking traces and moment estimates}\label{sec:moments}
\subsection{The non-cycle remainder}
Put $W_{N,k}=\Tr H_{G_N}^k$. A simple $j$-cycle contributes exactly $2j$ cyclically non-backtracking closed walks of length $k$ when $j\mid k$. Thus
\begin{equation}\label{L:eq:walkdecomp}
 W_{N,k}=\sum_{\substack{j\mid k\\j\geq3}}2jC_{N,j}+\mathcal R_{N,k},
 \qquad \mathcal R_{N,k}\geq0,
\end{equation}
where $\mathcal R_{N,k}$ counts the walks that are not repetitions of a simple cycle.

The following first-moment bound is also a consequence of \cite[Lemma 1.3.5 and the following paragraph]{FK}. We include an exploration proof, since the same construction is used to estimate pairs of walks.

\begin{lemma}[Two-collision bound]\label{L:lem:remainder}
For $N\geq2k$ and $k\geq3$,
\begin{equation}\label{L:eq:remainder}
 \E[\mathcal R_{N,k}]\leq C_B\frac{k^4q^k}{N}.
\end{equation}
\end{lemma}
\begin{proof}
Fix a lifted starting vertex and a prescribed non-backtracking base walk of length $k$. Reveal its lifted trajectory by exposing a matching pair only when it is not already known. A \emph{collision} is a newly exposed edge whose endpoint is a previously visited vertex. At each step, conditional on the complete exploration history, its probability is at most $k/(N-k)$; declare the collision indicator zero at a step using an already exposed edge.

If $S$ is the connected graph formed by the distinct edges traversed, the number of collisions is exactly $|E(S)|-|V(S)|+1$. A cyclically non-backtracking closed walk uses at least two distinct incident edges at every visited vertex. Therefore $S$ has minimum degree at least two. If its cycle rank is one, then $|E(S)|=|V(S)|$ and $S$ is a simple cycle. A non-backtracking closed walk on that cycle has a fixed orientation and is a repetition of the cycle. Every walk contributing to $\mathcal R_{N,k}$ consequently requires at least two collisions.

For any two step indices, conditioning at the later step and then at the earlier step bounds the probability of two collisions by $(k/(N-k))^2$. There are $n_0N$ starting vertices and $dq^{k-1}$ non-backtracking base walks of length $k$ from each projected starting vertex. It follows that
\[
 \E[\mathcal R_{N,k}]
 \leq n_0Ndq^{k-1}\binom{k}{2}\left(\frac{k}{N-k}\right)^2
 \leq C_B\frac{k^4q^k}{N}.
\]
Including trajectories that are not closed only enlarges the upper bound.
\end{proof}

The classical non-backtracking spectral identity for a finite $d$-regular graph is
\begin{equation}\label{L:eq:traceidentity}
 \Tr H_G^k=2q^{k/2}\Tr T_k\!\left(\frac{A_G}{2\sqrt q}\right)
            +( |E(G)|-|V(G)|)(1+(-1)^k).
\end{equation}
For example, it follows from the non-backtracking eigenvalue description in \cite[Lemma 5.3]{GJR}: each adjacency eigenvalue $a$ gives the two roots of $z^2-az+q=0$, whose $k$th powers sum to $2q^{k/2}T_k(a/(2\sqrt q))$, and the additional eigenvalues are $1$ and $-1$, each repeated $|E(G)|-|V(G)|$ times. The formula holds componentwise for disconnected graphs as well.

Since the old block in Section~\ref{sec:inputs} and the last term of \eqref{L:eq:traceidentity} are deterministic,
\begin{equation}\label{L:eq:centeredtrace}
 Y_{N,k}=q^{-k/2}(W_{N,k}-\E[W_{N,k}]).
\end{equation}
Combining \eqref{L:eq:walkdecomp} and \eqref{L:eq:centeredtrace} gives
\begin{equation}\label{L:eq:modecycle}
 Y_{N,k}=2q^{-k/2}\sum_{\substack{j\mid k\\j\geq3}}j(C_{N,j}-\mu_{N,j})
          +q^{-k/2}(\mathcal R_{N,k}-\E[\mathcal R_{N,k}]).
\end{equation}
The last term has expected absolute value at most $C_B k^4q^{k/2}/N$, by Lemma~\ref{L:lem:remainder}. This estimate controls finite-$N$ centring separately from the total variation approximation.

\subsection{Cycle-count covariances}
We estimate the second moments directly using the forcing couplings.

\begin{proposition}[Cycle covariance error]\label{L:prop:cyclecovariance}
For $3\leq i,j\leq R$ and $N\geq8R^2$,
\begin{equation}\label{L:eq:cyclecovariance}
 \left|\Cov(C_{N,i},C_{N,j})-\ind_{\{i=j\}}\lambda_i\right|
 \leq C_B\frac{ijq^{i+j}}{N}.
\end{equation}
\end{proposition}
\begin{proof}
Separating the identical-cycle term gives
\begin{align*}
 \Cov(C_{N,i},C_{N,j})
 &=\ind_{\{i=j\}}\left(\mu_{N,i}
              -\sum_{\alpha\in\Ga_{N,i}}p_\alpha^2\right)\\
 &\quad+\sum_{\substack{\alpha\in\Ga_{N,i},\ \beta\in\Ga_{N,j}\\
                         \beta\ne\alpha}}
       p_\alpha\left\{\Prob(I_\beta=1\mid I_\alpha=1)-p_\beta\right\}.
\end{align*}
The contribution of vertex-disjoint pairs is bounded in absolute value by
$4ij\mu_{N,i}\mu_{N,j}/N$, by Lemma~\ref{L:lem:disjoint}.
For intersecting pairs, use the sum of the conditional and unconditional probabilities. Lemma~\ref{L:lem:unconditional} and Proposition~\ref{L:prop:overlap} bound this contribution by
\[
 \frac{ij}{N}\mu_{N,i}\mu_{N,j}
       +C_B\frac{i^2j q^j}{N}\mu_{N,i}.
\]
Since $\mu_{N,i}\leq\lambda_i\leq C_Bq^i/i$, the total contribution of distinct cycles is at most $C_Bijq^{i+j}/N$, including when $i=j$.
For the identical-cycle term, Proposition~\ref{L:prop:mean} gives
$|\mu_{N,i}-\lambda_i|\leq C_Biq^i/N$, and
\[
 \sum_{\alpha\in\Ga_{N,i}}p_\alpha^2
       \leq (N-R)^{-3}\mu_{N,i}.
\]
When $i=j$, both errors are absorbed by the right-hand side of
\eqref{L:eq:cyclecovariance}.
\end{proof}

\subsection{A second moment for the non-cycle remainder}
We first record a deterministic counting fact. For a vertex set $A$ in a $d$-regular graph, the number of rooted oriented cyclically non-backtracking closed walks of length $k$ whose vertex set meets $A$ is at most
\begin{equation}\label{L:eq:walkshitset}
 |A|\,k\,d q^{k-1}.
\end{equation}
Indeed, root such a walk at a visit to $A$. There are at most
$|A|dq^{k-1}$ resulting non-backtracking trajectories and at most $k$ choices for restoring the original root. This is an upper bound even for periodic walks, for which some rotations coincide.

\begin{proposition}[Second-moment remainder bound]\label{L:prop:remaindersecond}
For $k\geq3$ and $N\geq4k$,
\begin{equation}\label{L:eq:remaindersecondfull}
 \E[\mathcal R_{N,k}^2]
    \leq C_Bq^{2k}\left(\frac{k^6}{N}+\frac{k^8}{N^2}\right).
\end{equation}
Consequently, for $N\geq8k^2$,
\begin{equation}\label{L:eq:remaindersecond}
 \E[\mathcal R_{N,k}^2]\leq C_B\frac{k^6q^{2k}}{N},
 \qquad
 \left\|\frac{q^{-k/2}}{\sqrt{2k}}
       (\mathcal R_{N,k}-\E[\mathcal R_{N,k}])\right\|_{L^2}
       \leq C_B\left(\frac{k^5q^k}{N}\right)^{1/2}.
\end{equation}
\end{proposition}
\begin{proof}
Expand $\mathcal R_{N,k}^2$ as the number of ordered pairs of walks counted by $\mathcal R_{N,k}$. Split the pairs according to whether their vertex sets intersect.

\emph{Intersecting pairs.}
The first walk visits at most $k$ vertices. By \eqref{L:eq:walkshitset}, at most $k^2dq^{k-1}$ rooted closed walks of length $k$ meet that vertex set. Dropping the requirement that the second walk be counted by the remainder can only increase the count. Thus, pathwise, the number of intersecting pairs is bounded above by
$k^2dq^{k-1}\mathcal R_{N,k}$. Lemma~\ref{L:lem:remainder} bounds its expectation by
\begin{equation}\label{L:eq:intersectremainder}
 C_B\frac{k^6q^{2k}}{N}.
\end{equation}
Two intersecting remainder walks may traverse the same graph of cycle rank two, so four separate collisions need not occur.

\emph{Vertex-disjoint pairs.}
Index the possible first trajectories by a lifted starting vertex and a prescribed non-backtracking base walk of length $k$. For each such index $\tau$, expose the trajectory and let $\mathcal F_\tau$ be the sigma-algebra generated by its full exploration transcript. Write $I_\tau$ for the indicator that this trajectory belongs to the remainder. The event $\{I_\tau=1\}$ is $\mathcal F_\tau$-measurable. On that event, its support $S$ has at most $k$ vertices and $k$ edges. By the deferred-decisions observation in Section~\ref{sec:deferred}, conditional on $\mathcal F_\tau$, the edge permutations are independent, uniformly distributed completions of the exposed partial bijections.

Fix a second starting vertex outside $V(S)$ and a second base walk of length $k$. Explore its lift, stopping if it hits $V(S)$. Before that stopping time, none of the already exposed edges of $S$ can be traversed. Hence a second trajectory that remains disjoint from $S$ and belongs to the remainder must have at least two fresh collisions with its own previously visited vertices, by the cycle-rank argument of Lemma~\ref{L:lem:remainder}.

Let $\mathcal H_{\tau,t}$ be the sigma-algebra generated by $\mathcal F_\tau$ and the second exploration through step $t$. Let $\xi_t$ indicate a fresh collision with the second trajectory's own previously visited vertices, and set $\xi_t=0$ after stopping and on steps with no fresh exposure. At most $2k$ pairs of the relevant permutation have been exposed, and at most $k$ candidate target vertices have already been visited by the second trajectory. Consequently,
\[
 \E[\xi_t\mid\mathcal H_{\tau,t-1}]\leq p,
 \qquad p=\frac{k}{N-2k}.
\]
For $a<b$, the tower property and the measurability of $\xi_a$ give
\begin{align*}
 \E[\xi_a\xi_b\mid\mathcal F_\tau]
 &=\E\!\left[\xi_a\E[\xi_b\mid\mathcal H_{\tau,b-1}]
                                \,\middle|\,\mathcal F_\tau\right]\\
 &\leq p\,\E[\xi_a\mid\mathcal F_\tau]\leq p^2.
\end{align*}
A disjoint remainder trajectory has $\sum_{t=1}^k\xi_t\geq2$. Since
$\ind_{\{\sum_t\xi_t\geq2\}}\leq\sum_{a<b}\xi_a\xi_b$, summing over second starting vertices and base walks gives the conditional expectation bound
\begin{equation}\label{L:eq:conditionaldisjointremainder}
 n_0N d q^{k-1}\binom{k}{2}
        \left(\frac{k}{N-2k}\right)^2
       \leq C_B\frac{k^4q^k}{N}.
\end{equation}
This conditional bound is uniform over all transcripts with $I_\tau=1$. Multiply by $I_\tau$, take expectations, and sum over the indices $\tau$ of the first trajectory. Since $\sum_\tau I_\tau=\mathcal R_{N,k}$, Lemma~\ref{L:lem:remainder} shows that the expected number of vertex-disjoint ordered pairs is at most
\[
 C_B\frac{k^4q^k}{N}\E[\mathcal R_{N,k}]
       \leq C_B\frac{k^8q^{2k}}{N^2}.
\]
Together with \eqref{L:eq:intersectremainder}, this proves
\eqref{L:eq:remaindersecondfull}. Under $N\geq8k^2$ the second term is absorbed by the first. Finally use
$\Var(\mathcal R_{N,k})\leq\E[\mathcal R_{N,k}^2]$ and the specified normalisation to obtain \eqref{L:eq:remaindersecond}.
\end{proof}

\subsection{Mixed cycle--remainder moments}
The mixed terms can be estimated directly using nonnegativity and conditional cycle counts.

\begin{proposition}[Mixed cycle--remainder bound]\label{L:prop:mixedremainder}
For integers $i,\ell\geq3$ and $N\geq2(i+\ell)$,
\begin{equation}\label{L:eq:mixedcovariance}
 \left|\Cov(C_{N,i},\mathcal R_{N,\ell})\right|
       \leq C_B\frac{\ell^5q^{i+\ell}}{N}.
\end{equation}
The same upper bound holds for $\E[C_{N,i}\mathcal R_{N,\ell}]$.
\end{proposition}
\begin{proof}
Index the possible remainder trajectories of length $\ell$ by their lifted starting vertex and prescribed non-backtracking base walk. For an index $\tau$, let $\mathcal F_\tau$ be its full exploration sigma-algebra and $I_\tau$ its remainder indicator. On $\{I_\tau=1\}$, let $S$ be its support, so $|V(S)|,|E(S)|\leq\ell$. Split the present simple $i$-cycles into those meeting $V(S)$ and those disjoint from it.

\emph{Cycles meeting the support.}
In every $d$-regular graph and for every vertex set $A$, the number of simple $i$-cycles meeting $A$ is at most $|A|dq^{i-1}/2$. Indeed, a simple cycle meeting $A$ gives at least two rooted oriented non-backtracking trajectories of length $i$ starting in $A$, whereas the total number of such trajectories is at most $|A|dq^{i-1}$. Unlike the count of rooted walks, this cycle count requires no factor $i$ to recover an original root. Therefore, pathwise, the number of ordered pairs consisting of a remainder $\ell$-walk and a simple $i$-cycle meeting its support is at most
\[
 \frac{\ell dq^{i-1}}{2}\,\mathcal R_{N,\ell}.
\]

\emph{Cycles disjoint from the support.}
Condition on a transcript with $I_\tau=1$. Its exposed pairs all belong to $S$, and the remaining pairs have the uniform completion law of Section~\ref{sec:deferred}. Fix a starting vertex $x_0$ outside $V(S)$ and a prescribed non-backtracking base walk $w$ of length $i$. Let $\mathcal H_{\tau,w}$ be the sigma-algebra generated by this first transcript and the first $i-1$ steps of the second exploration. Write $E_{\tau,w}$ for the event that $x_0,\ldots,x_{i-1}$ are distinct and outside $V(S)$; set this event to be empty unless $w$ is closed and cyclically non-backtracking. This event is $\mathcal H_{\tau,w}$-measurable.

On $E_{\tau,w}$, the final matching pair is fresh at both endpoints. The last vertex is first reached at step $i-1$, and the non-backtracking condition excludes its incoming edge. At the starting vertex, only the first edge has been used, and the cyclic non-backtracking condition excludes it as the final edge. No exposed edge of $S$ is incident to either vertex. Thus, the endpoint slots for the final matching pair are both unused, including when the inverse permutation is traversed. Let $r_{\tau,w}\leq\ell+i-1$ count its exposed pairs, and let $J_{\tau,w}$ indicate that the full trajectory is a simple $i$-cycle disjoint from $S$. Uniform completion gives
\[
 \E[J_{\tau,w}\mid\mathcal H_{\tau,w}]
   =\frac{\ind_{E_{\tau,w}}}{N-r_{\tau,w}}
   \leq\frac{\ind_{E_{\tau,w}}}{N-\ell-i}.
\]
We condition only on the observed prefix. Each present simple $i$-cycle disjoint from $S$ has exactly $2i$ rooted oriented representations, including cycles with a periodic base projection. Consequently, the conditional expected number of such cycles is at most
\begin{equation}\label{L:eq:conditionaldisjointcycles}
 \frac{n_0N d q^{i-1}}{2i(N-\ell-i)}
       \leq C_B\frac{q^i}{i}.
\end{equation}
This bound is uniform in the first transcript, since $N\geq2(i+\ell)$. Multiplying by $I_\tau$, averaging and summing over $\tau$ bounds the expected disjoint contribution by $C_Bq^i\E[\mathcal R_{N,\ell}]/i$.

Combining both contributions and using Lemma~\ref{L:lem:remainder} yields
\[
 \E[C_{N,i}\mathcal R_{N,\ell}]
 \leq C_B\left(\ell+\frac1i\right)q^i\E[\mathcal R_{N,\ell}]
 \leq C_B\frac{\ell^5q^{i+\ell}}N.
\]
Also, Proposition~\ref{L:prop:mean} gives $\mu_{N,i}\leq C_Bq^i/i$, and hence
\[
 \mu_{N,i}\E[\mathcal R_{N,\ell}]
       \leq C_B\frac{\ell^4q^{i+\ell}}{iN}.
\]
Both variables are nonnegative, so
$|\Cov(C_{N,i},\mathcal R_{N,\ell})|
\leq\E[C_{N,i}\mathcal R_{N,\ell}]+\mu_{N,i}\E[\mathcal R_{N,\ell}]$.
These bounds prove \eqref{L:eq:mixedcovariance}.
\end{proof}

\subsection{Pointwise polynomial-trace covariances}
Define
\begin{equation}\label{L:eq:normalizedmodes}
 X_{N,k}=Y_{N,k}/\sqrt{2k},\qquad
 K_B(k,\ell)=\frac{\Lambda_B(k,\ell)}{2\sqrt{k\ell}}.
\end{equation}
Set
\begin{equation}\label{L:eq:QandD}
 Q_{N,k}=\sqrt{\frac{2}{k}}\,q^{-k/2}
                \sum_{\substack{i\mid k\\i\geq3}}i(C_{N,i}-\mu_{N,i}),
 \qquad D_{N,k}=X_{N,k}-Q_{N,k}.
\end{equation}
Both variables are centred. For $k\geq3$ and $N\geq8k^2$, Proposition~\ref{L:prop:remaindersecond} gives
\begin{equation}\label{L:eq:Dnorm}
 \|D_{N,k}\|_{L^2}\leq C_B\left(\frac{k^5q^k}{N}\right)^{1/2}.
\end{equation}
Independence and $\Var(Z_j)=\lambda_j=b_j/(2j)$ give the normalised covariance:
\begin{align*}
 \Cov\left(\frac{Y_k}{\sqrt{2k}},\frac{Y_\ell}{\sqrt{2\ell}}\right)
 &=\frac{2q^{-(k+\ell)/2}}{\sqrt{k\ell}}
      \sum_{\substack{j\mid\gcd(k,\ell)\\j\geq3}}j^2\lambda_j\\
 &=K_B(k,\ell).
\end{align*}
Moreover, $K_B(k,k)\leq C_B$ for every $k$, since
$\sum_{j\leq k}j b_j\leq C_B kq^k$.

\begin{lemma}[Covariance of the cycle part]\label{L:lem:cyclepartcov}
For $3\leq k,\ell\leq R$ and $N\geq8R^2$,
\begin{equation}\label{L:eq:cyclepartcov}
 \left|\Cov(Q_{N,k},Q_{N,\ell})-K_B(k,\ell)\right|
 \leq C_B\frac{k^{3/2}\ell^{3/2}q^{(k+\ell)/2}}{N}.
\end{equation}
\end{lemma}
\begin{proof}
Apply Proposition~\ref{L:prop:cyclecovariance} term by term in
\eqref{L:eq:QandD}. The independent-Poisson part is exactly $K_B(k,\ell)$, and the error is at most
\[
 \frac{C_Bq^{-(k+\ell)/2}}{N\sqrt{k\ell}}
       \left(\sum_{i\mid k}i^2q^i\right)
       \left(\sum_{j\mid\ell}j^2q^j\right).
\]
The geometric bound
$\sum_{i\leq k}i^2q^i\leq C_q k^2q^k$ and its analogue for $\ell$ prove the result.
\end{proof}

\begin{proposition}[Pointwise covariance estimate]\label{lem:pointcov}
For $3\leq k,\ell\leq R$ and $N\geq8R^2$,
\begin{equation}\label{eq:pointcov}
 \left|\Cov(Y_{N,k},Y_{N,\ell})-\Lambda_B(k,\ell)\right|
 \leq\frac{C_Bq^{(k+\ell)/2}}{N}
       \left(k^2\ell^2+k\ell^5+k^5\ell+k^3\ell^3\right).
\end{equation}
Primitivity is not required.
\end{proposition}
\begin{proof}
In the notation of \eqref{L:eq:QandD}, decompose the normalised covariance as
\begin{align*}
 \Cov(X_{N,k},X_{N,\ell})-K_B(k,\ell)
 &=\Cov(Q_{N,k},Q_{N,\ell})-K_B(k,\ell)\\
 &\quad+\Cov(Q_{N,k},D_{N,\ell})+\Cov(D_{N,k},Q_{N,\ell})\\
 &\quad+\Cov(D_{N,k},D_{N,\ell}).
\end{align*}
Lemma~\ref{L:lem:cyclepartcov} controls the first term by
$C_B k^{3/2}\ell^{3/2}q^{(k+\ell)/2}/N$.
For a mixed term, \eqref{L:eq:modecycle} gives
$D_{N,\ell}=q^{-\ell/2}(\mathcal R_{N,\ell}-\E[\mathcal R_{N,\ell}])/\sqrt{2\ell}$.
Proposition~\ref{L:prop:mixedremainder} applies to every $i\mid k$, since
$N\geq8R^2\geq2(i+\ell)$, and yields
\begin{align*}
 |\Cov(Q_{N,k},D_{N,\ell})|
 &\leq\frac{q^{-(k+\ell)/2}}{\sqrt{k\ell}}
       \sum_{\substack{i\mid k\\i\geq3}} i\,
           |\Cov(C_{N,i},\mathcal R_{N,\ell})|\\
 &\leq\frac{C_B\ell^5q^{(\ell-k)/2}}{N\sqrt{k\ell}}
                    \sum_{i\mid k}i q^i\\
 &\leq C_B\frac{k^{1/2}\ell^{9/2}q^{(k+\ell)/2}}{N}.
\end{align*}
Here $\sum_{i\leq k} i q^i\leq C_q kq^k$. Interchanging $k$ and $\ell$ treats the other mixed term. Finally, \eqref{L:eq:Dnorm} and Cauchy--Schwarz give
\[
 |\Cov(D_{N,k},D_{N,\ell})|
 \leq C_B\frac{k^{5/2}\ell^{5/2}q^{(k+\ell)/2}}{N}.
\]
Multiply each of the four normalised bounds by $\sqrt{2k}\sqrt{2\ell}$ and add the results. This proves \eqref{eq:pointcov}.
\end{proof}

The same normalised bounds also imply
\begin{equation}\label{L:eq:covmain}
 \sup_{3\leq k,\ell\leq R}
 |\Cov(X_{N,k},X_{N,\ell})-K_B(k,\ell)|
 \leq C_B R^5q^R/N.
\end{equation}
In particular, for every fixed $L$,
\begin{equation}\label{L:eq:fixedcovlimit}
 \Cov(Y_{N,k},Y_{N,\ell})=\Lambda_B(k,\ell)+O_{B,L}(N^{-1}),
 \qquad 3\leq k,\ell\leq L.
\end{equation}
The growing-time transfer uses the pointwise estimate \eqref{eq:pointcov}, rather than only its uniform cutoff bound. The covariance scale alone can also be recovered from the core-graph method discussed in Section~\ref{sec:main}; the proof above uses only the explicit exploration estimates in this paper.

\section{An oscillatory Chebyshev energy limit}\label{sec:coefficients}
Set
\begin{equation}\label{eq:phase}
 \phi(x)=\arccos(r_d\cos x),\qquad
 f_s(x)=T_s(r_d\cos x)=\cos(s\phi(x)),\qquad x\in\mathbb R/(2\pi\mathbb Z).
\end{equation}
The phase is smooth and periodic because $0<r_d<1$. It is even and strictly increasing on $(0,\pi)$, with
\[
 \phi(0)=\alpha_d,\quad \phi(\pi)=\pi-\alpha_d,\qquad
 \phi'(x)=\frac{r_d\sin x}{\sqrt{1-r_d^2\cos^2x}}.
\]
The coefficients in \eqref{eq:chebcoeff} are the cosine coefficients of $f_s$. Periodic functions in this section have period $2\pi$. We use the unnormalised $L^1$ norm $\|g\|_{L^1}=\int_0^{2\pi}|g(x)|\dd x$, and $\|g\|_\infty=\sup_{x\in[0,2\pi]}|g(x)|$; distances on the circle are taken modulo $2\pi$.

\begin{lemma}[Uniform coefficient bounds]\label{lem:coeffbounds}
For every $s\geq0$,
\begin{equation}\label{eq:parseval}
 a_{s,0}^2+\frac12\sum_{k=1}^{s}a_{s,k}^2
     =\frac1\pi\int_0^\pi |f_s(x)|^2\dd x\leq1.
\end{equation}
In particular, $|a_{s,k}|\leq\sqrt2$ for $k\geq1$. For each fixed $p\geq0$ and $s\leq R$,
\begin{equation}\label{eq:weightedbound}
 \sum_{k=1}^{s}|a_{s,k}|k^p q^{k/2}
       \leq C_q R^p q^{R/2},\qquad
 \sum_{k=1}^{s}k a_{s,k}^2\leq2s.
\end{equation}
Also $a_{s,k}=0$ when $s$ and $k$ have different parity.
\end{lemma}
\begin{proof}
Parseval and $|\cos(s\phi)|\leq1$ prove \eqref{eq:parseval}. Cauchy--Schwarz and
$\sum_{k=1}^R k^{2p}q^k\leq R^{2p}q^R/(1-q^{-1})$ give the first bound in \eqref{eq:weightedbound}; $k\leq s$ gives the second. Finally $T_s(-r_dx)=(-1)^s T_s(r_dx)$ and the uniqueness of the Chebyshev expansion prove the parity assertion.
\end{proof}

The two stationary points $0,\pi$ are nondegenerate, since
\[
 \phi''(x)=\frac{r_d(1-r_d^2)\cos x}{(1-r_d^2\cos^2x)^{3/2}},\qquad
 \phi''(0)=\frac{r_d}{\sqrt{1-r_d^2}},\quad
 \phi''(\pi)=-\frac{r_d}{\sqrt{1-r_d^2}}.
\]
\begin{lemma}[Uniform oscillatory estimate]\label{TA:lem:osc}
For a continuously differentiable periodic complex-valued amplitude $A$ and $\omega\geq1$,
\begin{equation}\label{TA:eq:osc}
 \left|\int_0^{2\pi}e^{i\omega\phi(x)}A(x)\dd x\right|
 \leq C_{r_d}\omega^{-1/2}(\|A\|_\infty+\|A'\|_{L^1}).
\end{equation}
\end{lemma}
\begin{proof}
For large $\omega$, remove intervals of radius $\rho=\omega^{-1/2}$ around $0$ and $\pi$, on the circle. Their contribution is at most $C\rho\|A\|_\infty$. On the complementary intervals, $|\phi'|\geq c_{r_d}\rho$. Moreover $|\phi'|\geq r_d|\sin x|$ and $\phi''$ is bounded. On the two complementary intervals, the integral of $|\sin x|^{-2}$ is $4\cot\rho\leq C/\rho$ for small $\rho$, so
\[
 \int_{\operatorname{dist}(x,\{0,\pi\})\geq\rho}
        \frac{|\phi''(x)|}{|\phi'(x)|^2}\dd x\leq C_{r_d}/\rho.
\]
Integrate by parts on each complementary interval, using
$e^{i\omega\phi}=(i\omega\phi')^{-1}(e^{i\omega\phi})'$.
The boundary terms and the term involving $A\phi''/(\phi')^2$ are bounded by
$C_{r_d}\|A\|_\infty/(\omega\rho)$, and the term involving $A'/\phi'$ is bounded by
$C_{r_d}\|A'\|_1/(\omega\rho)$. Combining these bounds with the contribution from the removed intervals proves the result. For bounded $\omega$, the estimate follows after increasing $C_{r_d}$.
\end{proof}

\subsection{Bilinear energy and observation-time normalisation}
The leading covariance term scales with the geometric mean of the two observation times. We begin with the Fourier identity used in the proof.

Write $f_s(x)=\cos(s\phi(x))$ and $\Delta_u f(x)=f(x+u)-f(x)$. For real smooth periodic functions $f,g$, put
\[
 \widehat f(n)=\frac1{2\pi}\int_0^{2\pi}f(x)e^{-inx}\dd x.
\]
The homogeneous $H^{1/2}$ bilinear form satisfies the identity
\begin{equation}\label{eq:halfenergy}
 \mathcal E(f,g)=\sum_{n\in\mathbb Z}|n|\widehat f(n)\overline{\widehat g(n)}
 =\frac1{16\pi^2}\int_0^{2\pi}\int_{-\pi}^{\pi}
       \frac{\Delta_u f(x)\Delta_u g(x)}{\sin^2(u/2)}\dd u\dd x.
\end{equation}
This identity follows from Fourier orthogonality: it can first be verified for trigonometric polynomials and then extended by approximation. Since each $f_s$ is even,
\begin{equation}\label{eq:energycoeff}
 \mathcal E(f_s,f_{s'})=\frac12\sum_{k\geq1}ka_{s,k}a_{s',k}.
\end{equation}

\begin{lemma}[Uniform absolute-derivative remainder]\label{S:lem:absoluteD}
Let $\varphi\in C^\infty(\mathbb T;\mathbb R)$ have finitely many critical
points, all nondegenerate. With
\[
 D=-i\partial_x,\qquad |D|g=\sum_{k\in\mathbb Z}|k|\widehat g(k)e^{ikx},
\]
there is $C_\varphi<\infty$ such that, for every integer $n\geq1$,
\begin{equation}\label{S:eq:absoluteD}
 \bigl\||D|e^{in\varphi}-n|\varphi'|e^{in\varphi}\bigr\|_{L^1(\mathbb T)}
 \leq C_\varphi.
\end{equation}
The same estimate holds after taking real parts.
\end{lemma}
\begin{proof}
For smooth periodic $g$, the principal-value and symmetrised kernel formulas are
\begin{align}
 |D|g(x)&=\frac1{4\pi}\operatorname{p.v.}\!\int_{-\pi}^{\pi}
 \frac{g(x)-g(x+v)}{\sin^2(v/2)}\dd v,\label{S:eq:DkernelPV}\\
 |D|g(x)&=\frac1{8\pi}\int_{-\pi}^{\pi}
 \frac{2g(x)-g(x+v)-g(x-v)}{\sin^2(v/2)}\dd v.\label{S:eq:Dkernel}
\end{align}
The second integral is absolutely convergent. The identities follow first on Fourier modes from
$\int_{-\pi}^{\pi}(1-\cos(kv))/\sin^2(v/2)\dd v=4\pi|k|$
and then by smooth approximation. The remainder is estimated by separating neighbourhoods of the critical points from their complement.

Let $\mathcal C=\{\varphi'=0\}$ and
$\delta(x)=\min\{\delta_0,\operatorname{dist}_{\mathbb T}(x,\mathcal C)\}$, where
$\delta_0>0$ is fixed and small. Nondegeneracy and compactness give
$|\varphi'(x)|\asymp\delta(x)$. Put $F_n=e^{in\varphi}$ and
$R_n=|D|F_n-n|\varphi'|F_n$. If
$\delta(x)\leq \ell:=n^{-1/2}$, then
$|F_n''|\leq C_\varphi n$ on the relevant $2\ell$-neighbourhood.
Using \eqref{S:eq:Dkernel} for $|v|\leq\ell$ and the trivial numerator
bound outside it gives
\begin{equation}\label{S:eq:criticalR}
 |R_n(x)|\leq C_\varphi\sqrt n.
\end{equation}

Suppose now that $\delta:=\delta(x)>\ell$ and set
$\sigma=\operatorname{sgn}\varphi'(x)$. Choose a fixed
$\chi\in C_c^\infty((-1/4,1/4))$, with $0\leq\chi\leq1$ and
$\chi=1$ on $[-1/8,1/8]$. For fixed $x$ and $\delta$, let
$g(y)=\chi((y-x)/\delta)F_n(y)$ and $b=F_n-g$; the operators below act
in the $y$ variable. Since $g=F_n$ near $x$,
\[
 R_n(x)=(|D|-\sigma D)g(x)+|D|b(x).
\]
The first multiplier vanishes on Fourier modes with $\sigma k\geq0$.
For $\sigma k<0$, after $y=x+\delta v$ write
\[
 Q=n\delta^2+|k|\delta,\qquad
 \psi(v)=\frac{n\{\varphi(x+\delta v)-\varphi(x)\}-k\delta v}{Q}.
\]
On the fixed support of $\chi$,
$|\psi'|\geq c_\varphi$ and
$\|\psi^{(j)}\|_\infty\leq C_{\varphi,j}$ for $1\leq j\leq4$.
Indeed, for $j\geq2$,
$n\delta^j\|\varphi^{(j)}\|_\infty/Q\leq C_\varphi\delta^{j-2}$.
If $\mathcal L_\psi a=-(a/\psi')'$, three integrations by parts give
\[
 \int e^{iQ\psi(v)}\chi(v)\dd v
 =(iQ)^{-3}\int e^{iQ\psi(v)}\mathcal L_\psi^3\chi(v)\dd v,
 \qquad \|\mathcal L_\psi^3\chi\|_1\leq C_\varphi.
\]
Consequently
\begin{equation}\label{S:eq:oppfourier}
 |\widehat g(k)|\leq C_\varphi\delta
       (n\delta^2+|k|\delta)^{-3},\qquad \sigma k<0.
\end{equation}
With $A=n\delta^2\geq1$ and $K=\lfloor A/\delta\rfloor$,
\[
 \delta\sum_{k\geq1}\frac{k}{(A+k\delta)^3}
 \leq \frac{C}{\delta A}.
\]
The contribution from $k\leq K$ is at most
$\delta K(K+1)/(2A^3)\leq C/(\delta A)$, while for $k>K$ it is at most
$\delta^{-2}\sum_{k>K}k^{-2}\leq C/(\delta A)$.
Thus
\begin{equation}\label{S:eq:localR}
 |(|D|-\sigma D)g(x)|\leq C_\varphi(n\delta^3)^{-1}.
\end{equation}

For the nonlocal term, $b$ vanishes near $x$, so
\eqref{S:eq:DkernelPV} gives an ordinary integral,
\[
 |D|b(x)=\int_{\mathbb T}e^{in\varphi(y)}A_{x,\delta}(y)\dd y,
 \qquad
 A_{x,\delta}(y)=-\frac{1-\chi_{x,\delta}(y)}
                         {4\pi\sin^2((y-x)/2)}.
\]
Here $\chi_{x,\delta}$ is the periodic cutoff used to define $g$:
it equals $\chi((y-x)/\delta)$ in a local coordinate about $x$ and
is zero outside that coordinate neighbourhood. The amplitude extends
smoothly across $y=x$, since its numerator vanishes near that point, and
\[
 \|A_{x,\delta}\|_\infty+\|\partial_y A_{x,\delta}\|_{L^1}
 \leq C\delta^{-2}.
\]
Here the kernel and its derivative are $O(\operatorname{dist}(x,y)^{-2})$ and
$O(\operatorname{dist}(x,y)^{-3})$, while the cutoff derivative is
$O(\delta^{-1})$ on a set of length $O(\delta)$. Removing
$n^{-1/2}$-neighbourhoods of the critical points of $\varphi$ and
integrating by parts once on the complement yields the standard estimate
\[
 \left|\int_{\mathbb T}e^{in\varphi(y)}A(y)\dd y\right|
 \leq C_\varphi n^{-1/2}(\|A\|_\infty+\|A'\|_1).
\]
Hence
\begin{equation}\label{S:eq:nonlocalR}
 ||D|b(x)|\leq C_\varphi n^{-1/2}\delta^{-2}.
\end{equation}
Integrating \eqref{S:eq:criticalR} over the $O(n^{-1/2})$ critical
neighbourhoods and \eqref{S:eq:localR}--\eqref{S:eq:nonlocalR} over
their complement gives
\[
 \frac1n\int_{n^{-1/2}}^{\delta_0}v^{-3}\dd v=O(1),\qquad
 \frac1{\sqrt n}\int_{n^{-1/2}}^{\delta_0}v^{-2}\dd v=O(1).
\]
This proves \eqref{S:eq:absoluteD}. Taking real parts is legitimate because
$|D|$ commutes with complex conjugation.
\end{proof}

\begin{theorem}[Uniform bilinear energy remainder]\label{S:thm:energy}
There is a constant $C_d<\infty$ such that, for all positive integers
$s,u$,
\begin{equation}\label{S:eq:energy}
 \left|2\sum_{k\geq1}k a_{s,k}a_{u,k}
             -\sqrt{su}\,\Gamma_d(s-u)\right|\leq C_d.
\end{equation}
Consequently, for every $t\geq1$,
\begin{equation}\label{S:eq:naturalenergy}
 \sup_{t\leq s,u\leq2t}
 \left|\frac2{\sqrt{su}}\sum_{k\geq1}k a_{s,k}a_{u,k}
                 -\Gamma_d(s-u)\right|\leq\frac{C_d}{t}.
\end{equation}
All time indices are integers. The estimate holds for all positive observation times: the proof applies the multiplier estimate only to the $u$-factor and uses $\|f_s\|_\infty\leq1$, so no comparability assumption on $s$ and $u$ is required.
\end{theorem}
\begin{proof}
For $\phi(x)=\arccos(r_d\cos x)$ and $f_n(x)=\cos(n\phi(x))$,
\[
 \phi'(x)=\frac{r_d\sin x}{\sqrt{1-r_d^2\cos^2x}}.
\]
Its only critical points are $0$ and $\pi$, and they are nondegenerate.
Lemma~\ref{S:lem:absoluteD} therefore gives
\[
 |D|f_u=u|\phi'|f_u+r_u,\qquad \|r_u\|_1\leq C_d.
\]
By \eqref{eq:energycoeff} and Fourier orthogonality,
\[
 2\sum_{k\geq1}k a_{s,k}a_{u,k}
   =\frac2\pi\int_0^{2\pi}f_s(x)|D|f_u(x)\dd x.
\]
Since $\|f_s\|_\infty\leq1$, the product-to-sum identity and the change of
variables $\theta=\phi(x)$ on the two monotone arcs give, uniformly in
$s,u\geq1$,
\begin{align*}
 2\sum_{k\geq1}k a_{s,k}a_{u,k}
 &=u\{\Gamma_d(s-u)+\Gamma_d(s+u)\}+O_d(1).
\end{align*}
For $h\ne0$, \eqref{eq:paritykernel} gives
$|h\Gamma_d(h)|\leq4/\pi$. Hence
$|u\Gamma_d(s+u)|\leq4/\pi$, and, if $s\ne u$,
\[
 |u-\sqrt{su}|=\frac{\sqrt u\,|s-u|}{\sqrt s+\sqrt u}\leq|s-u|,
\]
so $|(u-\sqrt{su})\Gamma_d(s-u)|\leq4/\pi$; for $s=u$ this term is
zero. This proves \eqref{S:eq:energy}, and division by
$\sqrt{su}\geq t$ proves \eqref{S:eq:naturalenergy}.
\end{proof}

The bounded remainder permits a direct comparison at growing offsets and
removes the logarithmic loss from the observation-time covariance estimate.
Identifying the bounded remainder would require a separate band-edge analysis.

\subsection{An offset bound as a consequence}
The following form is useful for the quantitative common-time comparisons. It follows by dividing the bilinear estimate by $t$ and comparing the two normalisations.

\begin{theorem}[Growing-offset energy bound]\label{thm:energy}
For integers $t\geq2$ and $0\leq M\leq t/2$,
\begin{equation}\label{TA:eq:uniformenergy}
 \sup_{|a|,|b|\leq M}
 \left|\frac2t\sum_{k\geq1}ka_{t+a,k}a_{t+b,k}-\Gamma_d(a-b)\right|
 \leq C_d\frac{1+M}{t}.
\end{equation}
The supremum is over integer offsets. No average in $t$ is taken. In
particular, for each fixed pair of integer offsets,
\begin{equation}\label{eq:energylimit}
 \frac2t\sum_{k\geq1}k a_{t+a,k}a_{t+b,k}
       \longrightarrow\Gamma_d(a-b),\qquad t\to\infty.
\end{equation}
\end{theorem}
\begin{proof}
Apply Theorem~\ref{S:thm:energy} to $s=t+a$ and $u=t+b$. Since
$t-M\geq t/2$,
\[
 \left|\frac{\sqrt{(t+a)(t+b)}}t-1\right|
 \leq C\frac Mt,
\]
and $|\Gamma_d(h)|\leq c_d<2$. Dividing the bounded energy remainder by
$t$ proves \eqref{TA:eq:uniformenergy}. The fixed-offset limit follows
immediately.
\end{proof}

\section{The ideal Poisson trace field}\label{sec:poisson}
The field in \eqref{eq:idealsignal} is built from independent Poisson variables.
Interchanging the finite sums over divisors in \eqref{eq:idealsignal} gives
\begin{equation}\label{eq:betarep}
 P_s=\sum_{j=3}^{s}\beta_{s,j}(Z_j-\lambda_j),\qquad
 \beta_{s,j}=2j\sum_{m=1}^{\lfloor s/j\rfloor}a_{s,mj}q^{-mj/2}.
\end{equation}
Extend $\beta_{s,j}$ by zero for $j>s$. Lemma~\ref{lem:coeffbounds} implies
\begin{equation}\label{eq:betabound}
 |\beta_{s,j}|\leq B_j:=\frac{2\sqrt2\,j q^{-j/2}}{1-q^{-j/2}},
 \qquad\sup_{s\geq0,\ j\geq3}|\beta_{s,j}|\leq C_q.
\end{equation}

\subsection{Period and covariance}
For independent standard normals $\xi_j$, let
$P_s^{\mathrm G}=\sum_{j\geq3}\beta_{s,j}\sqrt{\lambda_j}\xi_j$.
In Proposition~\ref{F:prop:alltime}, we construct a coupling of this field with $P_s$.
\begin{lemma}[Non-bipartite and bipartite bases]\label{S:lem:period}
The non-backtracking matrix of a finite connected simple $d$-regular graph with $d\geq3$ is irreducible. Its period is $\mathfrak p_B$, and it is primitive exactly when $B$ is not bipartite. For some $\vartheta_B\in(0,1)$,
\begin{equation}\label{S:eq:perron}
 \left|b_kq^{-k}-\eta_B(k)\right|\leq C_B\vartheta_B^k,
 \qquad \eta_B(k)=1+(\mathfrak p_B-1)(-1)^k.
\end{equation}
If $B$ is bipartite, then, at every odd integer time $s$, both $\cA_N^0(s)$ and $\cA_N^*(s)$ vanish identically, as do $P_s$ and $P_s^{\mathrm G}$.
\end{lemma}
\begin{proof}
The directed graph of allowed non-backtracking transitions has indegree and outdegree $q$ at every directed edge. It is weakly connected. Indeed, at each base vertex, its incoming and outgoing arcs are joined by all transitions except immediate reversals. The underlying local bipartite graph is $K_{d,d}$ minus a perfect matching, which is connected for $d\geq3$. These local graphs share arcs along every base edge, and $B$ is connected. In any finite directed graph with equal indegree and outdegree at every vertex, a strongly connected component with no outgoing edges also has no incoming edges: sum the degree differences over that component. Weak connectivity therefore implies strong connectivity. This proves irreducibility.

The Ihara--Bass characteristic identity is
\begin{equation}\label{S:eq:bass}
 \det(zI-H_B)=(z^2-1)^{m_0-n_0}
                 \det(z^2I-zA_B+qI);
\end{equation}
see \cite{Bass}, or the equivalent eigenvalue statement in \cite[Lemma 5.3]{GJR}. If $|z|=q>1$ is an eigenvalue, then $z+q/z$ is a real adjacency eigenvalue. Writing $z=qe^{i\theta}$ gives imaginary part $(q-1)\sin\theta$, so $z=\pm q$. The adjacency eigenvalue $d$ is simple because $B$ is connected, and $-d$ is an eigenvalue exactly when $B$ is bipartite, in which case it is also simple. Both assertions follow from
\[
 x^{\mathsf T}(dI-A_B)x=\sum_{\{v,w\}\in E}(x_v-x_w)^2,
 \qquad
 x^{\mathsf T}(dI+A_B)x=\sum_{\{v,w\}\in E}(x_v+x_w)^2.
\]
In the second identity a nonzero kernel vector changes sign along every edge, which is possible exactly for a bipartition. The root $q$ in \eqref{S:eq:bass}, and the root $-q$ when it is present, are simple since $q>1$ and the corresponding adjacency eigenvalues are simple. Perron--Frobenius theory now identifies the period as one or two. All remaining eigenvalues have modulus less than $q$, and taking traces of powers proves \eqref{S:eq:perron}, whether or not the matrix is diagonalisable.

A lift of a bipartite graph is bipartite. Assign to a directed edge the sign of its initial vertex in the lifted bipartition, and let $J$ be the resulting diagonal sign matrix. A nonzero entry of $U_{G_N}$ connects arcs with opposite signs, so $JU_{G_N}J=-U_{G_N}$. Its odd traces are zero. Also $b_j=C_{N,j}=0$ for odd $j$. The coefficient parity in Lemma~\ref{lem:coeffbounds} shows that $\mathcal D_N(s)=m_N(s)=0$ for odd $s$ in the finite formulas of Section~\ref{TA:sec:centering}: the constant coefficient vanishes, and each nonzero coefficient has odd index. The formula for $P_s$ in \eqref{eq:betarep} has only zero contributions at odd $s$, as does the Gaussian formula in Proposition~\ref{F:prop:alltime}. This proves the last assertion.
\end{proof}

\begin{lemma}[Covariance correction and linear variance bound]\label{lem:summable}
For every base $B$ satisfying the standing assumptions, with $\eta_B$ as in
\eqref{S:eq:perron},
\begin{equation}\label{S:eq:periodsum}
 \sum_{k,\ell\geq1}
 \left|\Lambda_B(k,\ell)-2k\eta_B(k)\ind_{\{k=\ell\}}\right|<\infty.
\end{equation}
Uniformly for nonnegative integer $s,u$,
\begin{equation}\label{eq:poissoncovcompare}
 \left|\Cov(P_s,P_u)-2\sum_{k\geq1}k\eta_B(k)a_{s,k}a_{u,k}\right|\leq C_B,
 \qquad \Var(P_s)\leq8m_0s+C_B.
\end{equation}
For a non-bipartite base, $\eta_B(k)=1$.
\end{lemma}
\begin{proof}
Set
\[
 d_k=2k(b_kq^{-k}-\eta_B(k)),\qquad
 E_{k\ell}=2q^{-(k+\ell)/2}
 \sum_{\substack{j\mid\gcd(k,\ell),\ j\geq3\\(k,\ell)\ne(j,j)}}j b_j\geq0.
\]
Equation~\eqref{eq:idealY} gives
\[
 \Lambda_B(k,\ell)-2k\eta_B(k)\ind_{\{k=\ell\}}
       =d_k\ind_{\{k=\ell\}}+E_{k\ell}.
\]
The diagonal correction is absolutely summable by \eqref{S:eq:perron}.
For $z_j=q^{-j/2}$, Tonelli's theorem and $b_jq^{-j}\leq2m_0$ give
\begin{align*}
 \sum_{k,\ell\geq1}E_{k\ell}
 &=2\sum_{j\geq3}j b_j
       \left\{\left(\frac{z_j}{1-z_j}\right)^2-z_j^2\right\}\\
 &=2\sum_{j\geq3}j(b_jq^{-j})
              \frac{z_j(2-z_j)}{(1-z_j)^2}
 \leq C_B\sum_{j\geq3}j q^{-j/2}<\infty.
\end{align*}
Since all summands are nonnegative, Tonelli's theorem justifies the rearrangement.
The triangle inequality proves \eqref{S:eq:periodsum}; multiplication by the
bounded coefficients gives the covariance assertion.
For the variance, apply $b_k\leq2m_0q^k$ to the exact diagonal term.
Lemma~\ref{lem:coeffbounds} yields
\[
 \Var(P_s)\leq4m_0\sum_k k a_{s,k}^2+2\sum_{k,\ell}E_{k\ell}
             \leq8m_0s+C_B.
\]
\end{proof}

\begin{proposition}[Covariance with deterministic time growth]\label{S:prop:covariance}
Put
\[
 J_B(s,u)=\mathfrak p_B\,
       \ind_{\{\mathfrak p_B\mid s,\ \mathfrak p_B\mid u\}}
                       \Gamma_d(s-u).
\]
For all positive integers $s,u$,
\begin{equation}\label{S:eq:rawcov}
 |\Cov(P_s,P_u)-\sqrt{su}\,J_B(s,u)|\leq C_B.
\end{equation}
Consequently, for all positive integers $s,u$,
\begin{equation}\label{S:eq:cov_envelope}
 \left|\frac{\Cov(P_s,P_u)}{\sqrt{su}}-J_B(s,u)\right|
       \leq\frac{C_B}{\sqrt{su}}.
\end{equation}
In particular, uniformly for $s,u\in[t,2t]$,
\begin{equation}\label{S:eq:natcov}
 \left|\frac{\Cov(P_s,P_u)}{\sqrt{su}}-J_B(s,u)\right|
       \leq C_B/t.
\end{equation}
If $R=\max(s,u)$, $s,u\in[t,2t]$, and $N\geq8R^2$, then
\begin{equation}\label{S:eq:finitecov}
 |\Cov(\Tr U_{G_N}^s,\Tr U_{G_N}^u)-\sqrt{su}\,J_B(s,u)|
 \leq C_B\left\{1+\frac{R^6q^R}{N}\right\}.
\end{equation}
The error after division by $\sqrt{su}$ is bounded by
$C_B\{t^{-1}+R^5q^R/N\}$.
\end{proposition}
\begin{proof}
For a non-bipartite base, Lemma~\ref{lem:summable} and
Theorem~\ref{S:thm:energy} apply directly. For a bipartite base and even
$s,u$, coefficient parity gives
\[
 2\sum_k k\eta_B(k)a_{s,k}a_{u,k}=4\sum_k k a_{s,k}a_{u,k},
\]
which is twice the non-bipartite energy; this agrees with
$\mathfrak p_B=2$. If either time is odd, the ideal trace and the kernel
$J_B(s,u)$ vanish. The absolutely summable correction in
Lemma~\ref{lem:summable} is uniformly bounded, proving
\eqref{S:eq:rawcov}--\eqref{S:eq:natcov}.

For the finite lift use the exact representation \eqref{eq:tracecoef} and
Proposition~\ref{lem:pointcov}. With
$S_a(s)=\sum_{k=3}^s|a_{s,k}|k^a q^{k/2}$, the finite--ideal covariance
difference is bounded by
\begin{equation}\label{C:eq:finiteidealcov}
\begin{aligned}
 &|\Cov(\cA_N(s),\cA_N(u))-\Cov(P_s,P_u)|\\
 &\quad\leq\frac{C_B}{N}\{S_2(s)S_2(u)+S_1(s)S_5(u)
                      +S_5(s)S_1(u)+S_3(s)S_3(u)\}
 \leq C_B R^6q^R/N.
\end{aligned}
\end{equation}
This comparison uses only $s,u\leq R$ and $N\geq8R^2$. Combining it
with \eqref{S:eq:rawcov} proves \eqref{S:eq:finitecov}; division by
$\sqrt{su}\geq t$ and $R/t\leq2$ gives the last assertion.
\end{proof}

For the next lemma, put $T_M=(\Gamma_d(a-b))_{0\leq a,b<M}$ and let $\|\cdot\|_{\op}$ denote the Euclidean operator norm.
\begin{lemma}[Constants for the Grover covariance]\label{F:lem:toeplitz}
For $\gamma=\Gamma_d$ in \eqref{eq:kernel},
\begin{equation}\label{F:eq:gammaconstants}
 \sup_M\|T_M\|_{\op}\leq2,\qquad
 \sum_{h\in\mathbb Z}\Gamma_d(h)^2=2c_d.
\end{equation}
\end{lemma}
\begin{proof}
Let $\mathcal I_d=[\alpha_d,\pi-\alpha_d]\cup[-\pi+\alpha_d,-\alpha_d]$. Then
\[
 \Gamma_d(h)=\frac1\pi\int_{\mathcal I_d}e^{ih\theta}\dd\theta.
\]
For every finite complex vector $z$,
\[
 z^*T_Mz=\frac1\pi\int_{\mathcal I_d}\left|\sum_{a=0}^{M-1}z_ae^{ia\theta}\right|^2\dd\theta
 \leq2\sum_a|z_a|^2.
\]
This also proves positivity. In the normalised Fourier convention, $\Gamma_d(h)$ are the coefficients of $2\ind_{\mathcal I_d}$. Parseval gives
\[
 \sum_h\Gamma_d(h)^2=\frac1{2\pi}\int_{-\pi}^\pi4\ind_{\mathcal I_d}(\theta)\dd\theta
   =\frac{4(\pi-2\alpha_d)}\pi=2c_d.
\]
\end{proof}

\subsection{Fixed-offset fluctuations}
For a non-bipartite base, Theorem~\ref{thm:energy} and Lemma~\ref{lem:summable} give
\begin{equation}\label{eq:idealcovlimit}
 \frac1t\Cov(P_{t+a},P_{t+b})\longrightarrow\Gamma_d(a-b)
 \quad\text{for fixed integer }a,b.
\end{equation}
\begin{proposition}[Gaussian limit for the Poisson field]\label{prop:idealclt}
Under primitivity, for any fixed offsets $h_1,\ldots,h_p$,
\[
 (P_{t+h_a}/\sqrt t)_{a=1}^p\xrightarrow{\mathrm d}(G_{h_a})_{a=1}^p.
\]
The covariance entries converge, and convergence also holds in $\Wtwo$. The convergence holds jointly with any fixed finite vector $(Y_k)_{k=3}^{L}$, and the Gaussian field is independent of that vector in the limit.
\end{proposition}
\begin{proof}
In subscripts on constants, $h$ denotes the fixed list $(h_1,\ldots,h_p)$, including its length. Fix $c=(c_1,\ldots,c_p)\in\mathbb R^p$, and put
$u_{t,j}=t^{-1/2}\sum_a c_a\beta_{t+h_a,j}$. By \eqref{eq:betabound},
$\max_j|u_{t,j}|\leq C_{q,c}t^{-1/2}$. The variance of the resulting compensated Poisson sum
\[
 v_t=\sum_j\lambda_j u_{t,j}^2
\]
is bounded and, by \eqref{eq:idealcovlimit}, tends to
$\sum_{a,b}c_ac_b\Gamma_d(h_a-h_b)$. The characteristic function has the logarithm
\[
 \log\E\left[\exp\left(iz\sum_j u_{t,j}(Z_j-\lambda_j)\right)\right]
   =\sum_j\lambda_j(e^{izu_{t,j}}-1-izu_{t,j}),
\]
where the right-hand side specifies the continuous branch of the logarithm that vanishes at $z=0$. Taylor's formula bounds its difference from $-z^2v_t/2$ by
\begin{equation}\label{eq:charerror}
 \frac{|z|^3}{6}\sum_j\lambda_j|u_{t,j}|^3
 \leq\frac{|z|^3}{6}\max_j|u_{t,j}|\,v_t
 \leq C_{B,c,h}|z|^3t^{-1/2}.
\end{equation}
The characteristic-function argument and the Cram\'er--Wold device prove the Gaussian limit. Poisson variables with zero mean vanish identically.

To prove independence from the fixed low-mode vector, delete the terms $3\leq j\leq L$ in \eqref{eq:betarep}. Their normalised $L^2$ norm is $O_{B,L}(t^{-1/2})$ by \eqref{eq:betabound}. The remaining vector depends only on the independent variables $Z_j$ with $j>L$, and is exactly independent of $(Y_k)_{k=3}^{L}$. Thus the joint weak limit is the product law. Equation \eqref{eq:idealcovlimit} gives convergence of the expected squared Euclidean norms of the vectors at growing times; the low-mode vector has a fixed finite second moment. Weak convergence and convergence of second moments yield $\Wtwo$ convergence, also jointly with the low-mode vector.
\end{proof}

We verify nondegeneracy and then give an explicit representation of the limiting field. For distinct offsets and a nonzero real vector $c$,
\[
 \sum_{a,b}c_ac_b\Gamma_d(h_a-h_b)
  =\frac2\pi\int_{\alpha_d}^{\pi-\alpha_d}
       \left|\sum_a c_ae^{ih_a\theta}\right|^2\dd\theta>0.
\]
A nonzero trigonometric polynomial cannot vanish on this whole interval. Equivalently, with independent real Gaussian white-noise measures $W_1,W_2$ on this interval, each with Lebesgue control measure,
\begin{equation}\label{eq:fieldrepresentation}
 G_h=\sqrt{\frac2\pi}\int_{\alpha_d}^{\pi-\alpha_d}
        \{\cos(h\theta)\,W_1(\dd\theta)+\sin(h\theta)\,W_2(\dd\theta)\}.
\end{equation}
This defines a stationary real Gaussian process indexed by the integers, with covariance \eqref{eq:kernel}.

\section{Transfer to growing times on the finite lift}\label{sec:transfer}
The distributional and second-moment transfers are treated separately, using the distance conventions from Section~\ref{sec:main}.

\begin{proposition}[Finite-lift distributional comparison]\label{prop:BLtransfer}
Fix offsets $h_1,\ldots,h_p$ and set $H=\max_a|h_a|$. For integers $t\geq2H+3$, put $R=t+H$. If $N\geq8R^2$, then
\begin{equation}\label{eq:BLtransfer}
 \BL\left(\Law((\cA_N(t+h_a)/\sqrt t)_a),
           \Law((P_{t+h_a}/\sqrt t)_a)\right)
     \leq C_{B,h}\frac{R^4q^R}{N}.
\end{equation}
The same estimate, with a constant also depending on fixed $L$, holds after appending $(Y_{N,k})_{k=3}^L$ on the first side and $(Y_k)_{k=3}^L$ on the second, for $t$ sufficiently large.
\end{proposition}
\begin{proof}
Substitute \eqref{L:eq:modecycle} into \eqref{eq:tracecoef} and apply Lemma~\ref{L:lem:remainder}. The normalised total $L^1$ error caused by the centred remainder for a time $s\leq R$ is bounded by
\[
 \frac{C_B}{N\sqrt t}\sum_{k=3}^{s}|a_{s,k}|k^4q^{k/2}
   \leq C_B\frac{R^4q^{R/2}}{N\sqrt t}.
\]
Here \eqref{eq:weightedbound} was used. The cycle-count expression, still centred at $\mu_{N,j}$, is a deterministic map of the entire count vector up to $R$. By Theorem~\ref{L:thm:poisson}, its law changes by at most $C_BR^4q^R/N$ in total variation on replacing the cycle counts by independent Poisson variables of means $\lambda_j$. The corresponding bounded Lipschitz distance is at most twice this bound, regardless of the coefficients or the dimension of the original count vector.

For a single unnormalised Chebyshev mode, Proposition~\ref{L:prop:mean} bounds the change of centre by
\[
 2q^{-k/2}\sum_{j\mid k,\ j\geq3}j|\mu_{N,j}-\lambda_j|
 \leq C_B k^2q^{k/2}/N.
\]
After summing with $a_{s,k}$ and dividing by $\sqrt t$, the centre change is at most
$C_BR^2q^{R/2}/(N\sqrt t)$. Since $R/t$ is bounded and $q\geq2$, all these errors are absorbed by the right side of \eqref{eq:BLtransfer}. Appending a fixed low-mode vector adds only $C_{B,L}/N$ to the remainder and centring errors, while using the same Poisson approximation.
\end{proof}

\begin{proposition}[Finite-lift covariance comparison]\label{prop:covtransfer}
Fix integers $a,b$, put $H=\max\{|a|,|b|\}$, and for integers $t\geq2H+3$ set $R=t+H$. If $N\geq8R^2$, then
\begin{equation}\label{eq:covtransfer}
 \left|\frac1t\Cov(\cA_N(t+a),\cA_N(t+b))
       -\frac1t\Cov(P_{t+a},P_{t+b})\right|
     \leq C_{B,a,b}\frac{R^5q^R}{N}.
\end{equation}
Primitivity is not required for this comparison.
\end{proposition}
\begin{proof}
The weighted finite--ideal estimate \eqref{C:eq:finiteidealcov}
uses only $s,u\leq R$ and $N\geq8R^2$.
Apply it with $s=t+a$, $u=t+b$, and divide by $t$.
Since $R/t\leq3/2$, the result is at most
$C_B R^6q^R/(Nt)\leq C_{B,a,b}R^5q^R/N$.
\end{proof}

For the remainder of this section, suppose that the base is non-bipartite. For fixed offsets, define the deterministic quantity
\begin{equation}\label{eq:epsilon}
 \epsilon_{r_d,a,b}(t)=\left|
      \frac2t\sum_{k\geq1}k a_{t+a,k}a_{t+b,k}-\Gamma_d(a-b)\right|.
\end{equation}
Theorem~\ref{thm:energy} proves $\epsilon_{r_d,a,b}(t)\to0$. We obtain the following bound, which separates the two errors:
\begin{equation}\label{eq:covtotal}
 \left|\frac1t\Cov(\cA_N(t+a),\cA_N(t+b))-\Gamma_d(a-b)\right|
 \leq C_{B,a,b}\left(\frac{R^5q^R}{N}+\frac1t\right)
            +\epsilon_{r_d,a,b}(t).
\end{equation}
Theorem~\ref{thm:energy} gives the explicit bound
$\epsilon_{r_d,a,b}(t)\leq C_d(1+\max\{|a|,|b|\})/t$
for all sufficiently large $t$. This deterministic energy error is separate from the finite-lift error.

\begin{proof}[Proof of Theorem~\ref{thm:main} and Corollary~\ref{cor:ff}]
Condition \eqref{eq:weakcondition} ensures $N\geq8R_N^2$ eventually. Proposition~\ref{prop:BLtransfer} and Proposition~\ref{prop:idealclt} give \eqref{eq:mainclt}, including joint convergence with the independent fixed low-mode vector. Under \eqref{eq:momentcondition}, which implies \eqref{eq:weakcondition}, equation \eqref{eq:covtotal} gives convergence of all covariances of the high-time vectors. Their coordinates have mean zero, so their expected squared Euclidean norms converge to that of the Gaussian limit. The second moments of the fixed low-mode vector converge by \eqref{L:eq:fixedcovlimit}. Consequently, the full joint vector converges weakly and its second moments converge.

To see that this suffices for $\Wtwo$ convergence, take a Skorokhod coupling with almost sure vector convergence. The squared norms converge almost surely and their expectations converge. For nonnegative variables, the identity $|x-y|=x+y-2\min(x,y)$ implies convergence of squared norms in $L^1$, and therefore uniform integrability. Squared vector differences are bounded by twice the sum of the two squared norms, so their expectations tend to zero. This coupling proves $\Wtwo$ convergence. In particular, the covariances between the low modes and the growing-time coordinates also converge. Their limiting independence was established by the preceding joint weak convergence.

Taking $a=b=0$ in \eqref{eq:covtotal} proves \eqref{eq:ramplaw}. Formula \eqref{eq:paritykernel} follows by direct integration of \eqref{eq:kernel}.
\end{proof}

\section{Explicit centring and quantitative time averages}\label{sec:timeaverages}
We first give the centring and a transfer estimate valid for both periods
and both time normalisations. We then establish finite-size bounds for common-time
averages on a non-bipartite base. The observation-time law of large numbers is
proved in Section~\ref{S:sec:llnproof}.
\subsection{An explicit centring from the base}\label{TA:sec:centering}
The following centring is computed directly from the base graph and the number of sheets. Recall $\beta_{s,j}$ from \eqref{eq:betarep} and put
\[
 L_N(s)=\sum_{k=3}^s a_{s,k}q^{-k/2}\mathcal R_{N,k}.
\]
Define
\begin{align}
 \mathcal D_N(s)={}&N(m_0-n_0)(1+(-1)^s)+2n_0N a_{s,0}\notag\\
          &-N(m_0-n_0)\sum_{k=1}^s a_{s,k}q^{-k/2}(1+(-1)^k),\label{TA:eq:detpart}\\
 w_k^*&=\sum_{j\mid k,\ j\geq3}b_j,\qquad
 m_N(s)=\mathcal D_N(s)+\sum_{k=1}^s a_{s,k}q^{-k/2}w_k^*.\label{TA:eq:explicitcenter}
\end{align}
These are finite sums computable from $B$, $N$ and $s$. At $s=0$, we use the empty-sum convention, giving $m_N(0)=2m_0N$.

\begin{proposition}[Deterministic and exact mean centrings]\label{TA:prop:centering}
For all $s\geq0$,
\begin{equation}\label{TA:eq:rawtrace}
 \Tr U_{G_N}^s=\mathcal D_N(s)+\sum_{k=1}^s a_{s,k}q^{-k/2}W_{N,k}.
\end{equation}
In particular,
\begin{equation}\label{TA:eq:explicitdecomp}
 \Tr U_{G_N}^s-m_N(s)
    =\sum_{j=3}^s\beta_{s,j}(C_{N,j}-\lambda_j)+L_N(s).
\end{equation}
For $s\leq R$ and $N\geq8R^2$,
\begin{equation}\label{TA:eq:centererror}
 |\E[\Tr U_{G_N}^s]-m_N(s)|\leq C_B\frac{R^4q^{R/2}}{N}.
\end{equation}
\end{proposition}
\begin{proof}
For a regular graph with $n$ vertices and $m$ edges, the determinant identities imply
\[
 \Tr U_G^s=(m-n)(1+(-1)^s)+2\Tr T_s(A_G/d),
\]
\[
 2\Tr T_k(A_G/(2\sqrt q))
   =q^{-k/2}\{\Tr H_G^k-(m-n)(1+(-1)^k)\},\quad k\geq1.
\]
They hold componentwise for disconnected graphs. Substitute \eqref{eq:chebcoeff} and $n=n_0N$, $m=m_0N$. The constant coefficient contributes $2n_0N a_{s,0}$; this gives \eqref{TA:eq:rawtrace}. Since $w_k^*=\sum_{j\mid k}2j\lambda_j$, the cycle decomposition \eqref{L:eq:walkdecomp} gives \eqref{TA:eq:explicitdecomp}.

The difference between the exact mean and the explicit centring is bounded by
\[
 \sum_{k=3}^s|a_{s,k}|q^{-k/2}
 \left\{\sum_{j\mid k,\ j\geq3}2j|\mu_{N,j}-\lambda_j|+\E[\cR_{N,k}]\right\}.
\]
By Proposition~\ref{L:prop:mean}, the divisor sum is at most $C_B k^2q^k/N$; the remainder is at most $C_B k^4q^k/N$. Apply Lemma~\ref{lem:coeffbounds} with $p=4$ to prove \eqref{TA:eq:centererror}.
\end{proof}

The centring is thus an explicit finite sum. For numerical work, the coefficients can be evaluated by recurrence, and the new-sector compression in the first paper \cite{AKM} can be used to compute the fluctuating traces without numerically subtracting inherited eigenvalues. Its numerical conditioning and floating-point stability require separate analysis for large times or sheet numbers.

\subsection{A common transfer estimate}\label{TA:sec:transfer}\label{F:sec:transfer}
For $t\geq4$, $1\leq M\leq t/2$, $R=t+M$ and $N\geq8R^2$, set
$C_N=(C_{N,j})_{3\leq j\leq R}$ and
\[
 \mathcal X_{N,t,h}=t^{-1/2}\sum_{j=3}^R\beta_{t+h,j}(C_{N,j}-\lambda_j),
 \qquad Q^{\mathrm{cyc}}_{N,t,M}=M^{-1}\sum_{h<M}\mathcal X_{N,t,h}^2.
\]
For $3\leq k\leq R$, put
$\Psi_k(C_N)=2q^{-k/2}\sum_{j\mid k,\ j\geq3}j(C_{N,j}-\lambda_j)$.
Theorem~\ref{L:thm:poisson} compares every measurable function of $C_N$
with the same function of $(Z_j)_{3\leq j\leq R}$ in total variation,
with error at most $C_B\tau_{N,R}$. The same comparison applies to the joint map
consisting of a time average, a step or interpolated path, and the fixed vector of polynomials $\Psi_k$.
For each fixed $L$ with $R\geq L$, the first-moment and exact-mean bounds give
\begin{equation}\label{F:eq:lowerror}
 \E\left[\sum_{k=3}^L|Y_{N,k}-\Psi_k(C_N)|\right]\leq C_{B,L}/N.
\end{equation}
Indeed, \eqref{L:eq:modecycle} separates each difference into the centred
remainder and the change from $\mu_{N,j}$ to $\lambda_j$; both have this
bound for fixed $L$.

\begin{lemma}[Weighted transfer of squared traces]\label{TA:lem:maxerror}
Let $V_h=\cA_N^\star(t+h)/\sqrt t$, $X_h=\mathcal X_{N,t,h}$ and
$E_N^\star=\max_{h<M}|V_h-X_h|$, for either trace centring. In the finite
range above,
\begin{equation}\label{TA:eq:maxerror}
 \E[E_N^\star]\leq C_B\kappa_{N,t,M},\qquad
 \kappa_{N,t,M}=\frac{MR^4q^{R/2}}{N\sqrt t}.
\end{equation}
For any deterministic weights $0\leq w_h\leq1$, on each realisation,
\begin{align}
 \frac1M\sum_{h<M}w_h|V_h^2-X_h^2|
    &\leq2E_N^\star\sqrt{Q^{\mathrm{cyc}}_{N,t,M}}+(E_N^\star)^2,
                                          \label{TA:eq:squaredifference}\\
 \sup_{0\leq u\leq1}\left|\frac1{\sqrt M}
       \sum_{h<\lfloor Mu\rfloor}w_h(V_h^2-X_h^2)\right|
    &\leq\sqrt M\{2E_N^\star\sqrt{Q^{\mathrm{cyc}}_{N,t,M}}
                                      +(E_N^\star)^2\}.
                                          \label{F:eq:pathdifference}
\end{align}
The path bound also holds after interpolation. If $t\to\infty$ and
$\tau_{N,R}\to0$, within the same finite range, its right-hand side tends
to zero in probability. All these assertions hold for either period of the base.
\end{lemma}
\begin{proof}
By \eqref{TA:eq:explicitdecomp}, the error for explicit centring is
$L_N(t+h)/\sqrt t$. Nonnegativity of $\cR_{N,k}$ and
Lemmas~\ref{L:lem:remainder} and \ref{lem:coeffbounds} give
\[
 \E[|L_N(s)|]\leq\frac{C_B}{N}\sum_{k=3}^s|a_{s,k}|k^4q^{k/2}
                   \leq C_B R^4q^{R/2}/N,\qquad s\leq R.
\]
For exact-mean centring use also \eqref{TA:eq:centererror}. Bounding the
maximum by the sum proves \eqref{TA:eq:maxerror}.
The identity $V_h^2-X_h^2=2X_h(V_h-X_h)+(V_h-X_h)^2$ and
Cauchy--Schwarz in $h$ give both pathwise bounds. A partial sum is bounded
by the full sum of absolute values; interpolation takes convex combinations.
Moreover,
\begin{equation}\label{F:eq:scaledkappa}
 \frac{\kappa_{N,t,M}}{\tau_{N,R}}
       =\frac{M}{\sqrt t\,q^{R/2}}\leq C\sqrt R\,q^{-R/2},\qquad
 \frac{\sqrt M\kappa_{N,t,M}}{\tau_{N,R}}
       \leq C Rq^{-R/2}.
\end{equation}
Thus $\sqrt M E_N^\star\to0$ in probability. The ideal counterpart of
$Q^{\mathrm{cyc}}_{N,t,M}$ has bounded expectation by
\eqref{eq:poissoncovcompare}; the total variation comparison therefore implies tightness of
$Q^{\mathrm{cyc}}_{N,t,M}$. The first term on the right of
\eqref{F:eq:pathdifference} is a product of a variable converging to zero in probability and a tight variable,
and the second is $(\sqrt M E_N^\star)^2/\sqrt M$.
\end{proof}
The choices $w_h=1$ and $w_h=t/(t+h)$ give respectively common-time and
observation-time normalisation. The total-variation comparison acts on the
whole measurable map, including all deterministic normalisations; the actual trace remainder
is controlled separately by the lemma.

For observation-time normalisation, keeping each observation time in the
remainder estimate removes the restriction on the relative window length.
Put
\[
 \widetilde X_h=\frac1{\sqrt{t+h}}
       \sum_{j=3}^{R}\beta_{t+h,j}(C_{N,j}-\lambda_j),
 \qquad
 \widetilde Q_{N,t,M}=\frac1M\sum_{h<M}\widetilde X_h^2.
\]

\begin{lemma}[Observation-time transfer on unrestricted windows]\label{S:lem:owntransfer}
Let $t,M\geq1$, $R=t+M\geq3$, and $N\geq8R^2$. For either trace
centring set
\[
 \widetilde E_N^\star=
 \max_{h<M}\left|\frac{\cA_N^\star(t+h)}{\sqrt{t+h}}
                         -\widetilde X_h\right|.
\]
Then
\begin{equation}\label{S:eq:ownerror}
 \E[\widetilde E_N^\star]
       \leq C_B\frac{R^{7/2}q^{R/2}}{N},\qquad
 \sqrt M\,\E[\widetilde E_N^\star]
       \leq C_B\tau_{N,R}q^{-R/2}.
\end{equation}
On each realisation, the difference between the squared observation-time
trace process and its cycle counterpart is bounded in supremum norm by
\begin{equation}\label{S:eq:ownpath}
 \sqrt M\left\{2\widetilde E_N^\star\sqrt{\widetilde Q_{N,t,M}}
                       +(\widetilde E_N^\star)^2\right\}.
\end{equation}
The bound also holds for the interpolations. If $M\to\infty$ and
$\tau_{N,R}\to0$, this bound tends to zero in probability, without a
restriction on $M/t$.
\end{lemma}
\begin{proof}
Use the proof of Lemma~\ref{TA:lem:maxerror} at each time $s=t+h$,
with cutoff $s$ rather than $R$. Lemmas~\ref{L:lem:remainder} and
\ref{lem:coeffbounds}, together with \eqref{TA:eq:centererror}, give
\[
 \E\left[\left|\frac{\cA_N^\star(s)}{\sqrt s}
       -\frac1{\sqrt s}\sum_{j\geq3}\beta_{s,j}(C_{N,j}-\lambda_j)
                  \right|\right]
 \leq\frac{C_Bs^{7/2}q^{s/2}}{N}.
\]
For $s=1,2$, the cycle sum and the centred trace are both zero.
Bounding the maximum by the sum and using the geometric decay from the
largest time yields
\[
 \sum_{s=t}^{t+M-1}s^{7/2}q^{s/2}
       \leq C_qR^{7/2}q^{R/2}.
\]
Since $M\leq R$, this proves both inequalities in
\eqref{S:eq:ownerror}. The difference-of-squares identity and
Cauchy--Schwarz prove \eqref{S:eq:ownpath}, exactly as in
\eqref{F:eq:pathdifference}.

The ideal counterpart $Q_{t,M}^{\natural}=M^{-1}\sum_{h<M}P_{t+h}^2/(t+h)$
has uniformly bounded expectation because $\Var(P_s)\leq C_Bs$ for
$s\geq1$ by Lemma~\ref{lem:summable}. For every $A>0$, the
complete-cycle-vector comparison and Markov's inequality give
\[
 \Prob(\widetilde Q_{N,t,M}>A)
 \leq\Prob(Q_{t,M}^{\natural}>A)+C_B\tau_{N,R}
 \leq C_B/A+C_B\tau_{N,R}.
\]
Thus $\widetilde Q_{N,t,M}$ is tight along the stated sequences.
This comparison uses only probabilities of measurable events, not
expectations of unbounded functions.
Equation~\eqref{S:eq:ownerror} gives
$\sqrt M\widetilde E_N^\star\to0$ in probability, and hence both
terms in \eqref{S:eq:ownpath} tend to zero. Interpolation takes convex
combinations of neighbouring grid values.
\end{proof}

\subsection{Ideal moments and quantitative common-time averages}\label{TA:sec:ideal}
In this subsection, the common-time bounds require $\mathfrak p_B=1$; the cumulant identity holds for both periods.
Write
\begin{equation}\label{TA:eq:error}
 e(t,M)=\frac{M}{t},
 \qquad Q_{t,M}=\frac1{Mt}\sum_{h<M}P_{t+h}^2.
\end{equation}
For a non-bipartite base, Theorem~\ref{thm:energy} and
Lemma~\ref{lem:summable} imply, with one fixed constant $C_B$,
\begin{equation}\label{TA:eq:idealerror}
 \max_{a,b<M}|K_t(a,b)-\Gamma_d(a-b)|\leq D_B(t,M),
 \qquad D_B(t,M)=C_B\{e(t,M)+t^{-1}\},
\end{equation}
where $K_t(a,b)=t^{-1}\Cov(P_{t+a},P_{t+b})$.

\begin{proposition}[Ideal second and fourth moments]\label{TA:prop:idealLLN}
For every base and all nonnegative $s,u$,
\begin{equation}\label{TA:eq:fourthidentity}
\begin{aligned}
 \Cov(P_s^2,P_u^2)&=2\Cov(P_s,P_u)^2+\rho_B(s,u),\\
 \rho_B(s,u)&=\sum_{j\geq3}\lambda_j\beta_{s,j}^2\beta_{u,j}^2\in[0,C_B].
\end{aligned}
\end{equation}
For $\mathfrak p_B=1$, $t\geq4$ and $1\leq M\leq t/2$,
\begin{equation}\label{TA:eq:idealmse}
 |\E[Q_{t,M}]-c_d|\leq D_B(t,M),\qquad
 \E[|Q_{t,M}-c_d|^2]\leq C_B\{M^{-1}+t^{-2}+e(t,M)^2\}.
\end{equation}
If also $t,M\to\infty$ and
\begin{equation}\label{TA:eq:smallwindow}
 \sqrt M\{e(t,M)+t^{-1}\}\to0,
\end{equation}
then $M\Var(Q_{t,M})\to4c_d$. In particular, this holds for
$M=\lfloor t^\gamma\rfloor$ with fixed $0<\gamma<2/3$.
\end{proposition}
\begin{proof}
The joint cumulant-generating function is
$\sum_j\lambda_j\{\exp(\sum_a z_a\beta_{s_a,j})-1-\sum_a z_a\beta_{s_a,j}\}$.
Differentiating at zero gives \eqref{TA:eq:fourthidentity}.
The envelope $B_j$ in \eqref{eq:betabound} and
$\lambda_j\leq m_0q^j/j$ imply
$\sum_j\lambda_jB_j^4\leq C_B\sum_j j^3q^{-j}<\infty$.
For a non-bipartite base, averaging the diagonal of
\eqref{TA:eq:idealerror} bounds the bias, and
\begin{equation}\label{TA:eq:idealvar}
 \Var(Q_{t,M})=\frac2{M^2}\sum_{a,b<M}K_t(a,b)^2+O_B(t^{-2}).
\end{equation}
By Lemma~\ref{F:lem:toeplitz}, the row sums of squared entries of
$(\Gamma_d(a-b))$ are bounded. Apply
$K_t(a,b)^2\leq2\Gamma_d(a-b)^2+2D_B(t,M)^2$ and add the squared bias
to obtain \eqref{TA:eq:idealmse}.
For the scaled variance,
\[
 \frac2M\sum_{a,b<M}\Gamma_d(a-b)^2
 =2\sum_{|h|<M}(1-|h|/M)\Gamma_d(h)^2\longrightarrow4c_d.
\]
Replacing $\Gamma_d$ by $K_t$ introduces an error of at most
$C_B\{D_B(t,M)\log(2+M)+M D_B(t,M)^2\}$, using
$|\Gamma_d(h)|\leq C_d/(1+|h|)$. This vanishes under
\eqref{TA:eq:smallwindow}, as does the fourth-cumulant contribution
$O_B(M/t^2)$. For a power window $M=\lfloor t^\gamma\rfloor$, the condition
$\sqrt M\{e(t,M)+t^{-1}\}\to0$ holds for every fixed $0<\gamma<2/3$.
\end{proof}

\begin{proposition}[Finite-size probability and mean bounds]\label{C:prop:quantitative}
Suppose $\mathfrak p_B=1$. For $t\geq4$, $1\leq M\leq t/2$,
$R=t+M$, $N\geq8R^2$, $0<\eta\leq1$, either centring satisfies
\begin{equation}\label{TA:eq:probability}
 \Prob(|\cF_{N,t,M}^\star-c_d|>\eta)
 \leq C_B\left\{(1+\eta^{-1})\tau_{N,R}
     +\eta^{-2}\bigl(M^{-1}+t^{-2}+e(t,M)^2\bigr)\right\}.
\end{equation}
Uniformly for $0\leq h<M$,
\begin{equation}\label{F:eq:meaninput}
 \left|\E\left[\frac{\cA_N^\star(t+h)^2}{t}\right]-c_d\right|
 \leq C_B\{e(t,M)+t^{-1}+\delta_{N,R}
       +\ind_{\{\star=*\}}R^8q^R/(N^2t)\}.
\end{equation}
The same bound holds for $|\E[\cF_{N,t,M}^\star]-c_d|$.
The last term is at most $C_B\delta_{N,R}$ in this finite range.
\end{proposition}
\begin{proof}
For common-time weights, \eqref{TA:eq:squaredifference} gives
$|\cF^\star-Q^{\mathrm{cyc}}|\leq2E_N^\star\sqrt{Q^{\mathrm{cyc}}}+(E_N^\star)^2$.
Since $c_d<2$, on $|Q^{\mathrm{cyc}}-c_d|\leq\eta/2$ and
$E_N^\star\leq\eta/8$, this bound is smaller than $\eta/2$:
$\sqrt{5/2}/4+1/64<1/2$. Hence
\[
 \Prob(|\cF^\star-c_d|>\eta)
 \leq\frac{8\E[E_N^\star]}{\eta}+C_B\tau_{N,R}
                    +\Prob(|Q_{t,M}-c_d|>\eta/2).
\]
Markov's inequality, \eqref{TA:eq:maxerror}, \eqref{F:eq:scaledkappa}
and \eqref{TA:eq:idealmse} prove \eqref{TA:eq:probability}.
For the mean bound, the error in \eqref{C:eq:finiteidealcov} with $s=u=t+h$, divided by $t$,
is at most $C_B\delta_{N,R}$ because $R/t\leq3/2$.
Combining with \eqref{TA:eq:idealerror} gives \eqref{F:eq:meaninput} for
exact-mean centring. For explicit centring use
\[
 \E[\cA_N^*(s)^2]=\Var(\Tr U_{G_N}^s)
                  +(\E[\Tr U_{G_N}^s]-m_N(s))^2
\]
and \eqref{TA:eq:centererror}. The extra term is absorbed since
$R^3/(Nt)\leq3/16$. Averaging proves the time-average bound.
\end{proof}
As $t,M\to\infty$, $M/t\to0$ and $\tau_{N,R}\to0$ imply common-time self-averaging; adding $\delta_{N,R}\to0$ gives $L^1$ convergence by nonnegativity. Theorem~\ref{S:thm:lln} includes both periods and separates this common-time limit from the unrestricted observation-time limit. The variance limit in
Proposition~\ref{TA:prop:idealLLN} concerns the ideal field only.
\section{Functional limits and temporal averaging}\label{sec:functional}\label{sec:natural}
We prove the observation-time limit using a time-decaying Gaussian coupling
and an invariance principle for quadratic arrays with a product bound on
the covariance error. These estimates apply to arbitrary growing windows
of positive times. We then transfer the limit to the finite lift, identify
the common-time growth profile, and prove the law of large numbers directly.
\subsection{A simultaneous Gaussian coupling}\label{F:sec:coupling}
Write $\|V\|_{L^p}=(\E[|V|^p])^{1/p}$ for $1\leq p<\infty$.
For $\psi(x)=e^{|x|}-1$, let
$\|V\|_\psi=\inf\{a>0:\E[\psi(V/a)]\leq1\}$, and let $W_\psi$
be the infimum of this norm of the difference over couplings.
Rio's Corollary 2.2 \cite{Rio} gives
\[
 \sup_{\lambda\geq0}W_\psi(\Poi(\lambda),\mathcal N(\lambda,\lambda))<\infty.
\]
For each prescribed $\lambda_j$, choose a coupling whose error norm is at most one plus this
uniform bound, and standardise its Gaussian coordinate. When $\lambda_j=0$, take
$Z_j=0$ and any independent standard normal. The countable product gives
independent pairs $(Z_j,\xi_j)$ such that
\begin{equation}\label{C:eq:rio}
 \|Z_j-\lambda_j-\sqrt{\lambda_j}\xi_j\|_\psi\leq C,
 \qquad\|Z_j-\lambda_j-\sqrt{\lambda_j}\xi_j\|_{L^p}\leq C_p
 \quad(1\leq p<\infty).
\end{equation}
The second bound follows from $x^p\leq C_pe^x$ for $x\geq0$.
This single product coupling is used for all finite $p$ and all positive
integer times.

\begin{proposition}[Time-decaying ideal trace coupling]\label{F:prop:alltime}
Under the coupling in \eqref{C:eq:rio}, let
$P_s^{\mathrm G}=\sum_{j\geq3}\beta_{s,j}\sqrt{\lambda_j}\xi_j$.
For every fixed $1\leq p<\infty$,
\begin{equation}\label{F:eq:alltime}
 \left\|\sup_{s\geq1}\sqrt s\,|P_s-P_s^{\mathrm G}|\right\|_{L^p}
       \leq C_{p,d},\qquad
 \Cov(P_s^{\mathrm G},P_u^{\mathrm G})=\Cov(P_s,P_u).
\end{equation}
For all positive integers $t,M$, and positive deterministic numbers $d_h$,
\begin{equation}\label{F:eq:squarecoupling}
 \left\|\sup_{0\leq v\leq1}\left|
 \frac1{\sqrt M}\sum_{h<\lfloor Mv\rfloor}
       \frac{P_{t+h}^2-(P_{t+h}^{\mathrm G})^2}{d_h}
                         \right|\right\|_{L^p}
 \leq\frac{C_{p,B}}{\sqrt M}\sum_{h<M}\frac1{d_h}.
\end{equation}
The same bound holds after interpolation. In particular, with
\begin{equation}\label{S:eq:harmonic}
 H_{t,M}=\sum_{h=0}^{M-1}\frac1{t+h}\leq1+\log M,
\end{equation}
the observation-time error is at most $C_{p,B}H_{t,M}/\sqrt M$,
uniformly for $t\geq1$. The common-time error is at most
$C_{p,B}\sqrt M/t$.
\end{proposition}
\begin{proof}
Apply Lemma~\ref{TA:lem:osc} with amplitude $A(x)=\cos(kx)$.
The cosine coefficient formula for $f_s$ and
$\|A\|_\infty+\|A'\|_{L^1}\leq1+4k$ give
\begin{equation}\label{S:eq:coeffdecay}
 |a_{s,k}|\leq C_d(1+k)s^{-1/2},\qquad s,k\geq1.
\end{equation}
For $k>s$, the coefficient is zero. Put $z_j=q^{-j/2}$.
Equation~\eqref{eq:betarep} therefore implies
\begin{equation}\label{S:eq:betadecay}
 \sqrt s\,|\beta_{s,j}|\leq
 D_j:=2C_dj\left\{\frac{z_j}{1-z_j}
                    +\frac{jz_j}{(1-z_j)^2}\right\},
 \qquad \sum_{j\geq3}D_j<\infty.
\end{equation}
Writing $\zeta_j=Z_j-\lambda_j-\sqrt{\lambda_j}\xi_j$, we have
\[
 \sup_{s\geq1}\sqrt s\,|P_s-P_s^{\mathrm G}|
       \leq\sum_{j\geq3}D_j|\zeta_j|.
\]
Minkowski's inequality, \eqref{C:eq:rio} and monotone convergence
bound the right-hand side in every finite $L^p$. The supremum is
measurable because the time set is countable; the majorant is finite
almost surely. Both covariances in \eqref{F:eq:alltime} equal
$\sum_j\lambda_j\beta_{s,j}\beta_{u,j}$.

By Lemma~\ref{lem:summable}, Gaussian moments satisfy
$\|P_s^{\mathrm G}\|_{L^{2p}}\leq C_{p,B}\sqrt s$ for $s\geq1$.
The coupling gives the same bound for $P_s$. Hence H\"older's inequality yields
\[
 \|P_s^2-(P_s^{\mathrm G})^2\|_{L^p}
 \leq\|P_s-P_s^{\mathrm G}\|_{L^{2p}}
               \|P_s+P_s^{\mathrm G}\|_{L^{2p}}
 \leq C_{p,B}.
\]
Bounding a partial sum by the full sum of absolute values and applying
Minkowski proves \eqref{F:eq:squarecoupling}. Interpolation preserves
the bound. Finally, $t+h\geq h+1$ gives \eqref{S:eq:harmonic}.
\end{proof}

The same time decay controls the removal of fixed cycle lengths.
For fixed $L$, put
$P_s^{\leq L}=\sum_{3\leq j\leq L}\beta_{s,j}(Z_j-\lambda_j)$.
Equation~\eqref{S:eq:betadecay} gives
\[
 \left\|\sup_{s\geq1}\sqrt s\,|P_s^{\leq L}|\right\|_{L^p}
       \leq C_{p,B,L}.
\]
In particular, H\"older's inequality gives, uniformly for $s\geq1$,
\[
 \|P_s^2-(P_s-P_s^{\leq L})^2\|_{L^p}
 \leq\|P_s^{\leq L}\|_{L^{2p}}
       \bigl(2\|P_s\|_{L^{2p}}+\|P_s^{\leq L}\|_{L^{2p}}\bigr)
 \leq C_{p,B,L}.
\]
The $L^p$ norm of the path supremum of the deletion error is at most
$C_{p,B,L}H_{t,M}/\sqrt M$, also after interpolation. The remaining field depends only on $(Z_j)_{j>L}$, whereas
$(Y_k)_{k=3}^L$ depends only on $(Z_j)_{3\leq j\leq L}$.
They are independent for every window before passage to the limit.

\subsection{An invariance principle for Gaussian quadratic forms}\label{F:sec:gaussian}
For a matrix $A$, write $\|A\|_{\op}=\sup_{\|x\|_2=1}\|Ax\|_2$ and $\|A\|_F=(\sum_{a,b}|A(a,b)|^2)^{1/2}$ for the Euclidean operator norm and Frobenius norm. Matrix indices in this subsection run from $0$ to $M-1$ when the matrix size is $M$; the square root of a positive semidefinite matrix always denotes its positive semidefinite square root. The following array result applies to the normalised field in \eqref{F:eq:alltime}, which is not exactly stationary. The covariance matrices may be singular, and no uniform lower bound on their eigenvalues is needed. The row length may tend to infinity along any sequence. No coupling between different rows is assumed.

\begin{theorem}[Gaussian arrays with a product covariance-error bound]\label{F:thm:gaussian}
Let $\gamma:\mathbb Z\to\mathbb R$ be even, and suppose that
$T_M=(\gamma(a-b))_{0\leq a,b<M}$ is positive semidefinite with
\[
 \sup_M\|T_M\|_{\op}<\infty,\qquad
 g_2:=\sum_{h\in\mathbb Z}\gamma(h)^2<\infty.
\]
Let $X_M$ be a centred real Gaussian vector with covariance $K_M$.
Suppose there are nonnegative weights $w_{M,a}$ and a constant $C$
such that
\begin{equation}\label{F:eq:arraycondition}
 \begin{gathered}
 |K_M(a,b)-\gamma(a-b)|\leq Cw_{M,a}w_{M,b},\\
 H_M:=\sum_{a<M}w_{M,a}^2=o(\sqrt M),\qquad
 \sup_M\sum_{a<M}w_{M,a}^4<\infty.
 \end{gathered}
\end{equation}
Then
\[
 S_M(v)=\frac1{\sqrt M}\sum_{h<\lfloor Mv\rfloor}
                  \{X_{M,h}^2-K_M(h,h)\}
 \Longrightarrow\sqrt{2g_2}\,W(v)
\]
in $D([0,1])$ with $J_1$, and the polygonal interpolations converge
in $C([0,1])$. Both conclusions remain valid when $K_M(h,h)$ is
replaced by $\gamma(0)$.
\end{theorem}
\begin{proof}
\emph{Matrix estimates.}
Write $E_M=K_M-T_M$. For unit vectors $x,y$, the product bound gives
\[
 |x^{\mathsf T}E_My|
 \leq C\Bigl(\sum_a|x_a|w_{M,a}\Bigr)
          \Bigl(\sum_b|y_b|w_{M,b}\Bigr)\leq CH_M.
\]
It follows that
\begin{equation}\label{F:eq:operator}
 \|E_M\|_{\op}\leq CH_M,\qquad
 \|E_M\|_F\leq CH_M,\qquad
 \frac{\|K_M\|_{\op}}{\sqrt M}\longrightarrow0.
\end{equation}
Since $\|T_M\|_F^2\leq Mg_2$, Cauchy--Schwarz also yields
\begin{equation}\label{F:eq:squareserror}
 \frac1M\sum_{a,b<M}|K_M(a,b)^2-T_M(a,b)^2|
 \leq C\frac{H_M}{\sqrt M}+C\frac{H_M^2}{M}\longrightarrow0.
\end{equation}
Here the squares are entrywise squares, not matrix products.

\emph{Covariance of partial sums.}
Wick's formula gives
\[
 \Cov(S_M(v),S_M(w))
 =\frac2M\sum_{a<\lfloor Mv\rfloor}\sum_{b<\lfloor Mw\rfloor}K_M(a,b)^2.
\]
Replacing each entrywise square by $T_M(a,b)^2$ changes this by $o(1)$.
For a fixed lag $h$, the number of pairs in the rectangle with $a-b=h$,
divided by $M$, tends to $\min(v,w)$ and is at most one.
Dominated convergence against $\gamma(h)^2$ proves
\begin{equation}\label{F:eq:browncov}
 \Cov(S_M(v),S_M(w))\longrightarrow2g_2\min(v,w).
\end{equation}

\emph{Finite-dimensional convergence.}
A fixed linear combination of partial sums has the form
\[
 M^{-1/2}\{X_M^{\mathsf T}D_MX_M-\Tr(D_MK_M)\}
       \overset{\mathrm d}=\sum_r\theta_{M,r}(\xi_r^2-1),
\]
where $D_M$ is a real diagonal matrix with uniformly bounded entries,
$\xi_r$ are independent standard normals, and $\theta_{M,r}$ are the
eigenvalues of $M^{-1/2}K_M^{1/2}D_MK_M^{1/2}$.
This representation also holds for singular $K_M$.
Equation~\eqref{F:eq:operator} gives
$\max_r|\theta_{M,r}|\to0$, while
$2\sum_r\theta_{M,r}^2$ converges by \eqref{F:eq:browncov}.
For fixed real $z$, the logarithm of the characteristic function is
\[
 \sum_r\{-iz\theta_{M,r}-\tfrac12\log(1-2iz\theta_{M,r})\}
 =-z^2\sum_r\theta_{M,r}^2
   +O_z\Bigl(\max_r|\theta_{M,r}|\sum_r\theta_{M,r}^2\Bigr).
\]
Each logarithm is taken continuously from value zero at $z=0$.
The remainder tends to zero. If the limiting variance is zero,
convergence follows in $L^2$. The Cram\'er--Wold device gives the
asserted finite-dimensional convergence.

\emph{Tightness of the interpolations.}
For a centred Gaussian quadratic form $Q=\sum_r\theta_r(\xi_r^2-1)$,
\begin{equation}\label{F:eq:fourthquadratic}
 \E[Q^4]=12\Bigl(\sum_r\theta_r^2\Bigr)^2+48\sum_r\theta_r^4
       \leq15\Var(Q)^2.
\end{equation}
For $0\leq v<w\leq1$, an interpolated increment is
\[
 Q=\frac1{\sqrt M}\sum_{a<M}\omega_a\{X_{M,a}^2-K_M(a,a)\},
 \qquad 0\leq\omega_a\leq1,\quad \sum_a\omega_a=M(w-v).
\]
The weight $\omega_a$ is $M$ times the overlap length with the $a$th
grid interval. Using $K_M(a,b)^2\leq2T_M(a,b)^2+
2C^2w_{M,a}^2w_{M,b}^2$, we obtain
\begin{align*}
 \Var(Q)&\leq\frac C M\left\{\sum_a\omega_a+
                  \Bigl(\sum_a\omega_aw_{M,a}^2\Bigr)^2\right\}\\
 &\leq\frac C M\sum_a\omega_a=C(w-v).
\end{align*}
The first term uses $\sum_b\gamma(a-b)^2\leq g_2$ and $\omega_b\leq1$.
For the second, put $W_4=\sup_M\sum_{a<M}w_{M,a}^4<\infty$.
Weighted Cauchy--Schwarz gives
\[
 \Bigl(\sum_a\omega_aw_{M,a}^2\Bigr)^2
 \leq\Bigl(\sum_a\omega_a\Bigr)
       \Bigl(\sum_a\omega_aw_{M,a}^4\Bigr)
 \leq W_4\sum_a\omega_a.
\]
More explicitly, $\Var(Q)\leq4(g_2+C^2W_4)(w-v)$, where $C$ is
the constant in \eqref{F:eq:arraycondition}. This estimate uses the
fourth-power weight bound, whereas the covariance limit and the
vanishing largest eigenvalue use $H_M=o(\sqrt M)$.
Consequently
\begin{equation}\label{F:eq:tightness}
 \E[|\overline S_M(w)-\overline S_M(v)|^4]\leq C|w-v|^2.
\end{equation}
This holds also when either endpoint lies between grid times, so the
continuous-process tightness criterion applies.

\emph{Step processes and diagonal centring.}
The fourth-power weight bound implies $\sup_{M,a}w_{M,a}^2<\infty$,
and hence $\sup_{M,a}K_M(a,a)<\infty$.
Each jump has fourth moment $O(M^{-2})$, so a union bound gives
\[
 \Prob\left(\max_{a<M}\frac{|X_{M,a}^2-K_M(a,a)|}{\sqrt M}>\eta\right)
       \leq\frac{C}{M\eta^4}\longrightarrow0.
\]
Thus the step process and its interpolation have the same continuous
limit. Finally, replacing $K_M(a,a)$ by $\gamma(0)$ changes either
process by at most $CH_M/\sqrt M\to0$ in supremum norm.
\end{proof}

The uniform-error assumption
$\max_{a,b}|K_M(a,b)-\gamma(a-b)|\leq\epsilon_M$ with
$\sqrt M\epsilon_M\to0$ is included by taking
$w_{M,a}=\sqrt{\epsilon_M}$. The product bound also allows a covariance
error that is not uniformly small near the beginning of a row, provided both weight conditions in \eqref{F:eq:arraycondition} hold.

\subsection{The ideal observation-time limit}\label{F:sec:ideal}
Define
$S^{P,\natural}_{t,M}(u)=M^{-1/2}\sum_{h<\lfloor Mu\rfloor}
(P_{t+h}^2/(t+h)-c_d)$.
\begin{proposition}[Ideal limit for either period]\label{S:prop:ideal}
For any sequence of positive integer pairs $(t,M)$ with $M\to\infty$,
$S_{t,M}^{P,\natural}\Longrightarrow2\sqrt{\mathfrak p_Bc_d}\,W$
in $D([0,1])$ with $J_1$, and its polygonal interpolation converges in $C([0,1])$. The convergence holds jointly with each fixed vector $(Y_k)_{k=3}^L$, which is independent of the limiting Brownian motion.
\end{proposition}
\begin{proof}
Proposition~\ref{F:prop:alltime}, with $d_h=t+h$, bounds the uniform
$L^p$ distance between the Poisson and Gaussian squared processes by
$C_{p,B}H_{t,M}/\sqrt M$, which tends to zero uniformly in $t\geq1$.
The same bound holds after interpolation.

Suppose first that $B$ is non-bipartite. For the Gaussian row
$X_{M,a}=P_{t+a}^{\mathrm G}/\sqrt{t+a}$,
\eqref{S:eq:cov_envelope} gives
\[
 |K_M(a,b)-\Gamma_d(a-b)|\leq
              \frac{C_B}{\sqrt{(t+a)(t+b)}}.
\]
Take $w_{M,a}=(t+a)^{-1/2}$. Then
\[
 \sum_{a<M}w_{M,a}^2=H_{t,M}\leq1+\log M=o(\sqrt M),\qquad
 \sum_{a<M}w_{M,a}^4\leq\sum_{n\geq1}n^{-2}<\infty.
\]
Theorem~\ref{F:thm:gaussian} and Lemma~\ref{F:lem:toeplitz} therefore
give the Brownian limit with variance rate $4c_d$.

For a bipartite base, let $r_t=0$ when $t$ is even and $r_t=1$ otherwise. The even times in the window are
\[
 s_j=t+r_t+2j,\qquad 0\leq j<m,
 \qquad m=\#\{0\leq h<M:t+h\text{ is even}\}.
\]
Here $m/M\to1/2$, regardless of the parity of $t$ or $M$. Their normalised Gaussian covariance differs from the following Toeplitz kernel by at most $C_B/\sqrt{s_js_k}$:
\[
 \gamma_2(j-k)=2\Gamma_d(2(j-k)).
\]
The corresponding Toeplitz matrices have operator norm at most $4$: before multiplication by $2$ they are principal submatrices of the matrices in Lemma~\ref{F:lem:toeplitz}. Since the Grover covariance kernel vanishes at odd lags,
\[
 \sum_{j\in\mathbb Z}\gamma_2(j)^2
       =4\sum_{j\in\mathbb Z}\Gamma_d(2j)^2=8c_d.
\]
With $w_{m,j}=s_j^{-1/2}$, the sum of squared weights is at most
$1+\log(2m)$ and the sum of fourth powers is bounded by
$\sum_{n\geq1}n^{-2}$. Theorem~\ref{F:thm:gaussian} applied to this row gives
\[
 V_m(v)=\frac1{\sqrt m}\sum_{j<\lfloor mv\rfloor}
       \{(P_{s_j}^{\mathrm G})^2/s_j-2c_d\}
              \Longrightarrow4\sqrt{c_d}\,W(v).
\]

For $n=\lfloor Mu\rfloor$, let $j(n)$ count the even times among the first $n$ observations. The full Gaussian process is exactly
\begin{equation}\label{S:eq:evenembed}
 \sqrt{m/M}\,V_m(j(n)/m)+\frac{c_d}{\sqrt M}\{2j(n)-n\}.
\end{equation}
The second term is bounded uniformly by $c_d/\sqrt M$ because
$|2j(n)-n|\leq1$. Moreover, $|2m-M|\leq1$, so for $m\geq1$,
\[
 \sup_{u\in[0,1]}|j(\lfloor Mu\rfloor)/m-u|
 \leq m^{-1}+M^{-1}\leq C/M.
\]
Since $M\to\infty$, this applies to all sufficiently large windows.
The interpolations of $V_m$ are tight with continuous limits, and their maximum jump sizes tend to zero, by the proof of Theorem~\ref{F:thm:gaussian}. Their moduli of continuity at this deterministic time scale therefore tend to zero in probability. Thus composition in \eqref{S:eq:evenembed} preserves the limit. The factor $\sqrt{m/M}$ tends to $1/\sqrt2$, giving $2\sqrt{2c_d}\,W$. The full process also has vanishing maximum jump: its odd increments are deterministic of size $c_d/\sqrt M$, and its even increments are the active-row increments with the preceding bounded deterministic adjustment. Hence its polygonal interpolation has the same continuous limit. The parity of the initial time need not be fixed for this averaged limit.

Equation~\eqref{F:eq:squarecoupling} transfers both cases to the Poisson
field. Delete the finitely many cycle terms with $j\leq L$.
The bound following Proposition~\ref{F:prop:alltime} makes the resulting
process error at most $C_{B,L}H_{t,M}/\sqrt M$ in $L^1$.
The remaining process is exactly independent of $(Y_k)_{k=3}^L$.
Since the deletion error tends to zero, the full process has the same
independent joint limit.
\end{proof}

\subsection{Transfer to the finite lift}
\begin{proof}[Proof of Theorem~\ref{S:thm:natural}]
Write $t=t_N$, $M=M_N$, and $R=R_N$. The assumptions imply
$R\to\infty$, and eventually $N\geq8R^2$ and $R\geq L$.
Using the notation of Lemma~\ref{S:lem:owntransfer}, define
\[
 S_{N,t,M}^{\mathrm{cyc},\natural}(u)
  =\frac1{\sqrt M}\sum_{h<\lfloor Mu\rfloor}
                         \{\widetilde X_h^2-c_d\}.
\]
The complete-cycle-vector comparison shows that its law differs from
that of $S_{t,M}^{P,\natural}$ by at most $C_B\tau_{N,R}$ in total
variation. This also holds jointly with the fixed cycle polynomials
$\Psi_k$ and with the polygonal interpolation.
Lemma~\ref{S:lem:owntransfer} makes both trace-centred processes
uniformly close in probability to this same cycle process.
Proposition~\ref{S:prop:ideal} and Slutsky's theorem prove the process
limits and joint convergence of the two centrings to the same path.
Equation~\eqref{F:eq:lowerror}, whose fixed-degree proof needs no
comparison of $t$ and $M$, transfers the joint low-mode limit and its
independence. Evaluation at $u=1$ gives the scalar conclusion.
\end{proof}

\subsection{The common-time growth profile}
\begin{proof}[Proof of Corollary~\ref{S:cor:macro}]
Suppress the subscripts and let $U_j=\cZ^{\natural}(j/M)$.
For $n=\lfloor Mu\rfloor$, the profile-centred process equals
\begin{align}
 \sum_{h=0}^{n-1}(1+h/t)(U_{h+1}-U_h)
   &=(1+(n-1)/t)U_n-\frac1t\sum_{j=1}^{n-1}U_j.
 \label{S:eq:macroabel}
\end{align}
Both sides are zero when $n=0$. Since $M/t\to\rho<\infty$ and
$M\to\infty$, we also have $t\to\infty$.
Theorem~\ref{S:thm:natural} has a continuous limit, so convergence in
$J_1$ allows the use of uniform convergence after a Skorokhod representation.
The Riemann sums in \eqref{S:eq:macroabel} then converge uniformly to
$\rho\sigma_B\int_0^uW(v)\dd v$, and the first term converges to
$(1+\rho u)\sigma_BW(u)$. Deterministic integration by parts gives
\[
 (1+\rho u)W(u)-\rho\int_0^uW(v)\dd v
       =\int_0^u(1+\rho v)\dd W(v).
\]
This proves the step-process limit jointly with the natural process.
The profile-process jumps are the natural-process jumps multiplied by
weights bounded by $1+M/t$. Their maximum tends to zero in probability,
so the interpolations have the same limit. The mapping preserves the
joint conclusions for the two centrings and the low modes.

The limiting stochastic integral is a centred continuous Gaussian
process with independent increments and variance
$\int_0^u(1+\rho v)^2\dd v=u+\rho u^2+\rho^2u^3/3$.
This proves the time-change identity in law. At $u=1$ the deterministic
profile average is $c_d(1+(M-1)/(2t))$, giving the scalar limit.
\end{proof}

\subsection{The drift under common-time normalisation}
\begin{proof}[Proof of Corollary~\ref{S:cor:drift}]
Fix one centring and suppress $N,t,M,\star$ from the notation. Put
$U_j=\cZ^{\natural}(j/M)$ and $n=\lfloor Mu\rfloor$. Discrete summation by parts gives the exact identity
\begin{equation}\label{S:eq:abel}
 \cZ(u)-\cZ^{\natural}(u)
 =\frac{c_dn(n-1)}{2t\sqrt M}
   +\frac1t\left\{(n-1)U_n-\sum_{j=1}^{n-1}U_j\right\}.
\end{equation}
For $n=0$, the braces are interpreted as zero. Indeed, the left side equals
$M^{-1/2}\sum_{h<n}(h/t)\cA_N^\star(t+h)^2/(t+h)$, and the centred summand divided by $\sqrt M$ is $U_{h+1}-U_h$.

The expression in braces in \eqref{S:eq:abel} is bounded in absolute value by
$2M\|\cZ^{\natural}\|_\infty$. This norm is tight by Theorem~\ref{S:thm:natural}, and $M/t\to0$. Uniformly in $u$,
\[
 \left|\frac{n(n-1)}{2t\sqrt M}
                -\frac{M^{3/2}}{2t}u^2\right|
           \leq C\sqrt M/t\longrightarrow0.
\]
These two estimates prove \eqref{S:eq:drift}. Polygonal interpolation preserves the random-error bound; the uniform error in interpolating the quadratic drift is at most $C/(t\sqrt M)$. If $M^{3/2}/t\to\kappa<\infty$, addition of the deterministic continuous drift proves the asserted joint limits. If it tends to infinity, \eqref{S:eq:abel} at $u=1$ and tightness of $\cZ^{\natural}(1)$ imply $\cZ(1)\to+\infty$ in probability. The terminal marginals, and hence the processes, are not tight. The argument is pathwise and requires no asymptotic estimate of finite-lift expectations.
\end{proof}

\subsection{Self-averaging and averaged moment errors}\label{S:sec:llnproof}
The law of large numbers follows from direct ideal second moments and
local-in-time finite-lift moment bounds. In particular, the expectation
argument does not transfer unbounded moments through total variation.

\begin{proof}[Proof of Theorem~\ref{S:thm:lln}]
Write $t=t_N$, $M=M_N$, $R=t+M$ and $H=H_{t,M}$.
First consider $Q_{t,M}^{\natural}=M^{-1}\sum_{h<M}P_{t+h}^2/(t+h)$.
The diagonal case of \eqref{S:eq:rawcov} gives
\[
 |\E[Q_{t,M}^{\natural}]-c_d|\leq C_B(H+1)/M.
\]
Here the average of $\eta_B(t+h)$ differs from one by at most $1/M$.
The leading kernel $J_B(t+a,t+b)$ has bounded row sums of squared
entries. The product error in \eqref{S:eq:cov_envelope}, and the
fourth-cumulant identity \eqref{TA:eq:fourthidentity}, therefore imply
\begin{equation}\label{S:eq:naturalmse}
 \E[|Q_{t,M}^{\natural}-c_d|^2]
       \leq C_B\left\{\frac1M+\frac{(H+1)^2}{M^2}\right\}
       \longrightarrow0.
\end{equation}
To see this directly, write
$K_{ab}=\Cov(P_{t+a},P_{t+b})/\sqrt{(t+a)(t+b)}$ and
$J_{ab}=J_B(t+a,t+b)$. Then
\[
 \Var(Q_{t,M}^{\natural})
 =\frac2{M^2}\sum_{a,b<M}K_{ab}^2
  +\frac1{M^2}\sum_{a,b<M}
       \frac{\rho_B(t+a,t+b)}{(t+a)(t+b)}.
\]
The row sums $\sum_bJ_{ab}^2$ are bounded for either period.
Since $|K_{ab}-J_{ab}|\leq C_B/\sqrt{(t+a)(t+b)}$, the inequality
$K_{ab}^2\leq2J_{ab}^2+2(K_{ab}-J_{ab})^2$ gives
$\sum_{a,b<M}K_{ab}^2\leq C_B(M+H^2)$.
The fourth-cumulant contribution is at most $C_BH^2/M^2$.
Adding the squared bias gives the asserted mean-square estimate.
The complete-cycle-vector comparison transfers convergence in
probability to the cycle statistic. Dividing \eqref{S:eq:ownpath} by
$\sqrt M$ then transfers it to $\cF_{N,t,M}^{\natural,\star}$.

For expectations, use \eqref{C:eq:finiteidealcov} with $s=u$ and
cutoff $s$. For $s\geq3$, the variance error after division by $s$
is at most $C_Bs^5q^s/N$; for $s=1,2$ it is zero. Summing the geometric
weights before bounding them gives
\[
 \frac1M\sum_{s=t}^{t+M-1}
    \frac{|\Var(\cA_N^0(s))-\Var(P_s)|}{s}
       \leq\frac{C_B}{NM}\sum_{s=t}^{t+M-1}s^5q^s
       \leq C_B\frac{\delta_{N,R}}M.
\]
For explicit centring, \eqref{TA:eq:centererror} with cutoff $s$ adds
at most
\[
 \frac{C_B}{N^2M}\sum_{s=t}^{t+M-1}s^7q^s
       \leq C_B\frac{R^7q^R}{N^2M}
       = C_B\frac{\delta_{N,R}}M\frac{R^2}{N}.
\]
Consequently
\begin{equation}\label{S:eq:averagedmean}
 |\E[\cF_{N,t,M}^{\natural,\star}]-c_d|
 \leq C_B\left\{\frac{H+1}{M}+\frac{\delta_{N,R}}M
                +\ind_{\{\star=*\}}\frac{R^7q^R}{N^2M}\right\}.
\end{equation}
The cutoff assumption implies $R^2/N\to0$. Thus the last term is
$\delta_{N,R}/M$ times the vanishing factor $R^2/N$, and
the expectation error tends to zero under
\eqref{S:eq:averagedmoment}.
For nonnegative $F_N$ converging in probability to a constant $c\geq0$,
\[
 \E[|F_N-c|]=\E[F_N]+c-2\E[\min(F_N,c)].
\]
Since $\min(F_N,c)$ is bounded and converges to $c$ in probability,
$\E[\min(F_N,c)]\to c$. Thus convergence of expectations
proves the asserted $L^1$ limit.

For the common-time statistic assume $M/t\to\rho<\infty$ and put
$v_h=1+h/t$ and $a_{t,M}=M^{-1}\sum_{h<M}v_h=1+(M-1)/(2t)$.
The same variance argument, with the weights retained, gives
\[
 \Var\left(\frac1M\sum_{h<M}v_h\frac{P_{t+h}^2}{t+h}\right)
 \leq C_B(1+M/t)^2\left\{\frac1M+\frac{H^2}{M^2}\right\}
 \longrightarrow0.
\]
Here $1+M/t$ is bounded along the sequence under consideration.
Its expectation satisfies
\[
 \left|\frac1M\sum_{h<M}\frac{\E[P_{t+h}^2]}t-c_da_{t,M}\right|
       \leq C_B(t^{-1}+M^{-1}).
\]
For a bipartite base, the additional main term is an alternating sum
of $v_h$; its absolute value is bounded by $C(1+M/t)$ before division
by $M$. This follows by pairing successive weights, with at most one
unpaired endpoint. The bounded covariance remainder contributes
$C_B/t$.
Thus the ideal common-time average converges in probability to
$c_d(1+\rho/2)$. Total variation and the unscaled version of
\eqref{S:eq:ownpath}, multiplied by $\max_hv_h$, transfer this limit
to either finite-lift statistic. The finite--ideal expectation
errors are likewise at most $\max_hv_h$ times the preceding locally
summed bounds, so they tend to zero under
\eqref{S:eq:averagedmoment}. Nonnegativity proves the common-time
$L^1$ limit as well.
\end{proof}

\subsection{Fixed-offset fluctuations for bipartite bases}\label{S:sec:bipartite}
Time averages over consecutive observations have the same limit regardless of the initial-time parity. For fixed offsets, the odd-time coordinates are identically zero, and the initial-time parity determines the limiting vector.

\begin{theorem}[Fixed-offset limit for a bipartite base]\label{S:thm:bipoffset}
Suppose $B$ is bipartite. Fix $r\geq1$ integer offsets $h_1,\ldots,h_r$, put $H=\max_a|h_a|$, and let $t_N\to\infty$ have a fixed parity $\varepsilon\in\{0,1\}$. Set $R_N=t_N+H$. Under $R_N^4q^{R_N}/N\to0$, for either trace centring,
\[
 \left(\cA_N^\star(t_N+h_a)/\sqrt{t_N}\right)_{a=1}^r
 \Longrightarrow
 \left(\sqrt2\,\ind_{\{\varepsilon+h_a\text{ even}\}}G_{h_a}\right)_{a=1}^r,
\]
where $G$ has covariance \eqref{eq:kernel}. The convergence holds jointly with each fixed low-mode vector, and the two limiting vectors are independent. If $R_N^5q^{R_N}/N\to0$ as well, it holds in $\Wtwo$, also jointly with the low-mode vector, and all covariance entries converge. In particular, along even times with $t_N^5q^{t_N}/N\to0$,
$\mathcal F_N^{\mathrm c}(t_N)/t_N\to2c_d$; at odd times the form factor is identically zero.
\end{theorem}
\begin{proof}
The odd-time coordinates vanish by Lemma~\ref{S:lem:period}. For the even coordinates, apply \eqref{S:eq:rawcov} at reference time $\min_a(t_N+h_a)$, which is comparable to $t_N$. Their normalised ideal covariance tends to $2\Gamma_d(h_a-h_b)$. The compensated-Poisson representation \eqref{eq:betarep} has maximum coefficient $O(t_N^{-1/2})$ in each fixed linear combination. The characteristic-exponent proof of Proposition~\ref{prop:idealclt} therefore applies with this covariance, including degenerate coordinates. Deleting finitely many $Z_j$ gives independence from the low modes exactly as in that proof.

Proposition~\ref{prop:BLtransfer} transfers the joint law; its proof uses no primitivity. The deterministic difference between the two trace centrings is negligible by \eqref{TA:eq:centererror} and its first-moment proof. Under the additional condition, Proposition~\ref{prop:covtransfer} and \eqref{L:eq:fixedcovlimit} give convergence of second moments. Weak convergence with convergent second moments yields the claimed $\Wtwo$ convergence by the argument in the proof of Theorem~\ref{thm:main}. The even-time form-factor conclusion is the diagonal case.
\end{proof}

\subsection{Common-time and expectation-centred consequences}\label{F:sec:common}\label{F:sec:mean}
\begin{corollary}[Common-time consequences]\label{F:thm:main}
Suppose $\mathfrak p_B=1$, $t_N,M_N\to\infty$, $R_N=t_N+M_N$ and
\begin{equation}\label{F:eq:conditions}
 \sqrt{M_N}\{e(t_N,M_N)+t_N^{-1}\}\to0,
 \qquad\tau_{N,R_N}\to0.
\end{equation}
For either trace centring, $\cZ_{N,t_N,M_N}^\star\Longrightarrow2\sqrt{c_d}W$
in $D([0,1])$ with $J_1$, and the interpolations converge in $C([0,1])$.
The processes corresponding to the two centrings are uniformly close in probability and converge jointly
to the same path. Convergence holds jointly with each fixed low-mode vector,
and the Brownian motion is independent of the limiting low-mode vector.
In particular,
\[
 \sqrt{M_N}(\cF_{N,t_N,M_N}^\star-c_d)
                \Longrightarrow\mathcal N(0,4c_d).
\]
If also $\delta_{N,R_N}/\sqrt{M_N}\to0$, the same conclusions hold when
each summand $\cA_N^\star(t_N+h)^2/t_N$ is centred by its own expectation.
The scalar statistic is then centred at $\E[\cF_{N,t_N,M_N}^\star]$.
\end{corollary}
\begin{proof}
Since $e(t,M)=M/t$, the temporal condition entails both $M/t\to0$ and $M^{3/2}/t\to0$.
Theorem~\ref{S:thm:natural} and Corollary~\ref{S:cor:drift} with $\kappa=0$
give the process and joint limits; evaluation at $u=1$ gives the scalar one.
For expectation centring, Corollary~\ref{F:cor:expectation} below
applies with $\rho=0$ and gives the asserted conclusions under
$\delta_{N,R}/\sqrt M\to0$.
\end{proof}
For a non-bipartite base, the ideal common-time process
$S^P_{t,M}(u)=M^{-1/2}\sum_{h<\lfloor Mu\rfloor}(P_{t+h}^2/t-c_d)$
has the same functional and low-mode joint limits as $t,M\to\infty$
under the temporal condition in \eqref{F:eq:conditions}, without a condition on $N$.
Indeed, apply Proposition~\ref{S:prop:ideal} and \eqref{S:eq:abel} to the
ideal squared traces to obtain
\[
 \|S^P_{t,M}-S^{P,\natural}_{t,M}\|_\infty
 \leq\frac{c_dM^{3/2}}{2t}+
           \frac{2M}{t}\|S^{P,\natural}_{t,M}\|_\infty\to0
 \quad\text{in probability}.
\]
Interpolation preserves this bound. The finite-lift cycle comparison,
also jointly with a time average and fixed low modes, and its fluctuation-scale
remainder follow from Section~\ref{TA:sec:transfer} with $w_h=1$.

\subsection{Expectation-centred squared traces}\label{F:sec:expectation}
The two trace centrings need not centre their squares. To subtract the
expectation of every squared observation, define
\begin{align}
 \cZ_{N,t,M}^{\natural,\star,\mathrm E}(u)
   &=\cZ_{N,t,M}^{\natural,\star}(u)
        -\E[\cZ_{N,t,M}^{\natural,\star}(u)],
       \label{F:eq:expectnatural}\\
 \cZ_{N,t,M}^{\star,\mathrm E}(u)
   &=\cZ_{N,t,M}^{\star}(u)
        -\E[\cZ_{N,t,M}^{\star}(u)].
       \label{F:eq:expectcommon}
\end{align}
Thus the summands in \eqref{F:eq:expectnatural} are
$\{\cA_N^\star(t+h)^2-\E[\cA_N^\star(t+h)^2]\}/(t+h)$,
with the same $M^{-1/2}$ factor as before; in
\eqref{F:eq:expectcommon} the denominator is $t$.
Rather than apply the uniform bound \eqref{F:eq:meaninput} separately
at every observation, we sum the local moment errors first.

\begin{corollary}[Expectation-centred functional limits]\label{F:cor:expectation}
Assume \eqref{S:eq:conditions} and, in addition,
\begin{equation}\label{F:eq:expectcondition}
 \frac{\delta_{N,R_N}}{\sqrt{M_N}}\longrightarrow0.
\end{equation}
Set $\sigma_B=2\sqrt{\mathfrak p_Bc_d}$. For either trace centring,
\[
 \cZ_{N,t_N,M_N}^{\natural,\star,\mathrm E}
       \Longrightarrow\sigma_BW
       \quad\text{in }D([0,1])\text{ with }J_1,
\]
and the polygonal interpolations converge in $C([0,1])$.
The expectation-centred processes and the processes in
Theorem~\ref{S:thm:natural}, for both choices of trace centring,
are uniformly close in probability and converge jointly to the same
Brownian motion. This convergence holds jointly with every fixed vector
$(Y_{N,k})_{k=3}^L$, whose limiting vector is independent of $W$.

If also $M_N/t_N\to\rho\in[0,\infty)$, then
\[
 \cZ_{N,t_N,M_N}^{\star,\mathrm E}(u)
       \Longrightarrow\sigma_B\int_0^u(1+\rho v)\dd W(v)
\]
in $D([0,1])$ with $J_1$, and the polygonal interpolations converge in
$C([0,1])$. This holds jointly with the preceding natural processes and
the profile-centred processes of Corollary~\ref{S:cor:macro}, with the
same Brownian motion, for both trace centrings and with the same
independent limiting low-mode vector. In particular,
\begin{align*}
 \sqrt{M_N}\bigl(\cF_{N,t_N,M_N}^{\natural,\star}
             -\E[\cF_{N,t_N,M_N}^{\natural,\star}]\bigr)
     &\Longrightarrow\mathcal N(0,4\mathfrak p_Bc_d),\\
 \sqrt{M_N}\bigl(\cF_{N,t_N,M_N}^{\star}
             -\E[\cF_{N,t_N,M_N}^{\star}]\bigr)
     &\Longrightarrow\mathcal N\!\left(
       0,4\mathfrak p_Bc_d(1+\rho+\rho^2/3)\right),
\end{align*}
where the second conclusion uses $M_N/t_N\to\rho$.
\end{corollary}
\begin{proof}
Write $t=t_N$, $M=M_N$, $R=t+M$, and $H=H_{t,M}$.
Eventually $N\geq8R^2$. For each $3\leq s<R$,
\eqref{C:eq:finiteidealcov} with cutoff $s$, together with
\eqref{TA:eq:centererror}, gives
\begin{equation}\label{F:eq:localexpectation}
 \left|\frac{\E[\cA_N^\star(s)^2]}s
                  -\frac{\E[P_s^2]}s\right|
 \leq C_B\left\{\frac{s^5q^s}{N}
        +\ind_{\{\star=*\}}\frac{s^7q^s}{N^2}\right\}.
\end{equation}
Indeed, exact trace centring leaves only the variance difference. For
explicit trace centring, the extra term is
$\{\E[\Tr U_{G_N}^s]-m_N(s)\}^2/s$.
At $s=1,2$, the centred traces and $P_s$ are zero for both choices.
The ideal diagonal estimate \eqref{S:eq:rawcov} implies
\[
 \left|\frac{\E[P_s^2]}s-c_d\eta_B(s)\right|\leq C_B/s,
 \qquad s\geq1.
\]
Here $\eta_B$ is defined in \eqref{S:eq:perron}. For either period and
every integer $0\leq n\leq M$,
\[
 \left|\sum_{h=0}^{n-1}\{\eta_B(t+h)-1\}\right|\leq1.
\]
In the bipartite case this is a partial alternating sum; its cancellation
is retained before taking absolute values.

Put $b_{N,t,M}^\star(u)=\E[\cZ_{N,t,M}^{\natural,\star}(u)]$.
For the other errors in \eqref{F:eq:localexpectation}, bounding every
partial sum by the full sum of absolute values gives
\begin{align}
 \|b_{N,t,M}^\star\|_\infty
 &\leq\frac{C_B(1+H)}{\sqrt M}
      +\frac{C_B}{N\sqrt M}\sum_{s=t}^{t+M-1}s^5q^s
      +\frac{C_B\ind_{\{\star=*\}}}{N^2\sqrt M}
           \sum_{s=t}^{t+M-1}s^7q^s\notag\\
 &\leq C_B\left\{\frac{1+H}{\sqrt M}
       +\frac{\delta_{N,R}}{\sqrt M}
            \left(1+\ind_{\{\star=*\}}\frac{R^2}{N}\right)\right\}.
       \label{F:eq:expectdrift}
\end{align}
The geometric sums are taken before their terminal-cutoff bounds are
used. Since $H\leq1+\log M$ and $\tau_{N,R}\to0$ implies $R^2/N\to0$,
the last expression tends to zero under \eqref{F:eq:expectcondition}.
Subtracting this deterministic path from the natural process proves
all its conclusions by Theorem~\ref{S:thm:natural} and Slutsky's theorem.
Interpolation is linear and preserves the uniform error bound.

For the common-time process suppose $M/t\to\rho<\infty$, and let
$B_j=b_{N,t,M}^\star(j/M)$. The deterministic shift to subtract from the profile-centred process in
\eqref{S:eq:macro} to obtain \eqref{F:eq:expectcommon} is
\[
 \sum_{h=0}^{n-1}(1+h/t)(B_{h+1}-B_h)
   =(1+(n-1)/t)B_n-\frac1t\sum_{j=1}^{n-1}B_j,
 \qquad n=\lfloor Mu\rfloor,
\]
with both sides zero for $n=0$. Its supremum norm is at most
$(1+2M/t)\|b_{N,t,M}^\star\|_\infty$, which tends to zero by
\eqref{F:eq:expectdrift} and boundedness of $M/t$.
Summation by parts retains the cancellation in the weighted bipartite
sum. Corollary~\ref{S:cor:macro} now proves the
common-time limit and all the joint conclusions. The interpolated
centring error has the same bound. Evaluation at $u=1$ proves the
scalar limits.
\end{proof}

Condition~\eqref{F:eq:expectcondition} is sufficient and is not claimed
necessary. For example, along the sheet numbers
$N_n=(n^2+n)^5q^{n^2+n}$, take $t_n=n^2$, $M_n=n$, and $R_n=n^2+n$.
Then $\delta_{N_n,R_n}=1$ and $\tau_{N_n,R_n}=1/R_n$, while
$\delta_{N_n,R_n}/\sqrt{M_n}=n^{-1/2}\to0$.
Thus expectation centring does not require $\delta_{N,R}$ itself to
vanish. The common-time expectation-centred result is established when
$M/t\to\rho<\infty$, without the additional assumption $M^{3/2}/t\to0$.

\section{Conclusions and spectral interpretation}\label{sec:interpretation}
A single-time observation and a time average have different limits.
For a non-bipartite base, Theorem~\ref{thm:main} and the continuous mapping
theorem give
\[
 \frac{|\cA_N(t_N)|^2}{t_N}\xrightarrow{\mathrm d}c_d\chi_1^2.
\]
Thus the normalised squared trace retains randomness at a single growing
time. In contrast, Theorem~\ref{S:thm:lln} recovers $c_d$ by averaging
observation-time normalised measurements on one lift, and
Theorem~\ref{S:thm:natural} describes the cumulative error by
$2\sqrt{\mathfrak p_Bc_d}W$. These conclusions allow a fixed initial time
and impose no additional restriction on $M/t$ beyond the cycle cutoff.
They hold with both exact mean centring and the explicit centring $m_N$.
For bipartite bases, odd-time traces vanish and the leading even-time
variance is doubled, giving Brownian variance $8c_d$ rather than $4c_d$.
Fixed low-degree statistics retain the base-dependent Poisson structure
and are independent of the temporal limits in the joint limiting law.

The common-time normalisation separates a different deterministic
profile from the fluctuations. For $M/t\to\rho<\infty$,
Corollary~\ref{S:cor:macro} gives the variance clock
$u+\rho u^2+\rho^2u^3/3$ after subtracting that profile.
On windows with $M/t\to0$, leaving it unsubtracted produces the quadratic
drift in Corollary~\ref{S:cor:drift}. Table~\ref{tab:normalisations}
summarises the corresponding centring and scaling choices.

\paragraph{Covariance and the Grover band.}
The covariance is determined by the band $[\alpha_d,\pi-\alpha_d]$, rather than by the Fourier transform of the Kesten--McKay pushforward density. In the bilinear energy calculation, $|\phi'|$ is the Jacobian of $\theta=\phi(x)$ and yields Lebesgue measure on that band. Equivalently, $\Gamma_d$ is the Fourier transform of
$\pi^{-1}\ind_{\mathcal I_d}(\theta)\dd\theta$, with the two-band set
$\mathcal I_d$ of Lemma~\ref{F:lem:toeplitz}.
Reflection symmetry makes this transform real; invariance under a shift by $\pi$ gives zero odd-lag covariances. The even and odd subsequences of the non-bipartite limiting Gaussian field are independent, although correlations remain within each subsequence. For example, $c_3\approx1.567306208$, $\Gamma_3(2)\approx-0.400140585$, while $c_4=4/3$ and $\Gamma_4(2)=-\sqrt3/\pi$.

\paragraph{Relation to the orthogonal covariance form.}
For real polynomials $F,G$, write $F(2\cos x)=c_0(F)+\sum_{k\geq1}c_k(F)\cos(kx)$. In the normalisation of \cite[Section 1.3.1 and Theorem 2.5]{CES}, the leading orthogonal covariance form is
\[
 \mathcal Q_{\mathrm O}(F,G)=\frac12\sum_{k\geq1}k c_k(F)c_k(G).
\]
For $F_s(y)=2T_s(r_dy/2)$, one has $c_k(F_s)=2a_{s,k}$ and hence
\[
 \mathcal Q_{\mathrm O}(F_s,F_u)=2\sum_{k\geq1}k a_{s,k}a_{u,k}.
\]
Theorem~\ref{S:thm:energy} evaluates this form for all positive observation times. The period correction and finite-lift estimates justify the use of this form for the present ensemble. Transferring a Wigner-matrix limit theorem would require more than equality of the covariance forms. The Grover determinant and old/new decomposition remain established algebraic inputs \cite{KS,SKKM,AKM}.

\paragraph{Time scale and scope.}
The evolution operator acts on a space of dimension $D_N=2m_0N$. The logarithmic range satisfies $t_N/D_N\to0$; under the moment condition of Corollary~\ref{cor:ff},
\[
 \mathcal F_N^{\mathrm c}(t_N)/D_N=(c_d+o(1))t_N/D_N\to0.
\]
The linear law concerns original eigenphases without unfolding, at times small compared with the matrix dimension. All finite-lift limits are taken under the unconditioned covering-permutation law. They neither assert almost-sure convergence under a coupling across sheet numbers nor use the sampling conditional on connectedness in \cite{AKM}.

The squared trace is also distinct from an average return probability. With $\delta_e$ the coordinate unit vector at the directed edge $e$ in $\ell^2(D(G))$, in general
\[
 \left|\sum_e\langle\delta_e,U_G^s\delta_e\rangle\right|^2
 \ne\sum_e|\langle\delta_e,U_G^s\delta_e\rangle|^2.
\]
These conclusions concern spectral statistics rather than position-space
transport, localisation, mixing, or local eigenvalue universality.
The cutoff conditions are sufficient, not asserted optimal.
Proposition~\ref{TA:prop:idealLLN} gives an ideal-model variance limit;
the finite-lift weak limits do not by themselves yield convergence of
the scaled variance or fourth moments. Expectation-centred versions of both squared-trace processes are
established in Corollary~\ref{F:cor:expectation} under the additional
sufficient condition $\delta_{N,R}/\sqrt M\to0$. The common-time version
also assumes $M/t\to\rho<\infty$, without the short-window drift
condition. These moment conditions are not asserted necessary. Higher-order energy asymptotics, finite-lift
moment limits and path-space Wasserstein bounds remain separate
questions.

\section*{Statements and Declarations}
\paragraph{Funding.} This work was supported by JSPS KAKENHI Grant Number 23K12507.

\paragraph{Competing interests.} The authors declare no conflicts of interest.

\paragraph{Data availability.} The results are theoretical and are established by the proofs in this article. No external datasets are required for these results.

\begingroup\footnotesize
\bibliographystyle{abbrvurl}
\bibliography{references}
\endgroup
\end{document}